\pdfoutput=1
\documentclass[11pt,a4paper]{article}

\usepackage{latexsym,amssymb,amsmath}
\usepackage{MnSymbol, mathtools}
\usepackage[english]{babel}
\usepackage[latin1]{inputenc}
\usepackage[pdftex]{color}

\newcommand{\efe}{\mathbb F}
\newcommand{\FF}{\mathbb F}

\newcommand{\M}{\mathbb{M}}

\def\rank{\mathop{\rm rank}\nolimits}

\def\diag{\mathop{\rm Diag}\nolimits}

\def\max{\mathop{\rm max}\nolimits}

\def\halmos{\rule{6pt}{6pt}}

\newcommand{\wh}{\widehat}

\newcommand{\la}{\lambda}

\def\col{\mathop{\rm col}\nolimits}

\newcommand{\be}{\begin{equation}}
\newcommand{\ee}{\end{equation}}
\newtheorem{theorem}{Theorem}[section]
\newtheorem{definition}[theorem]{Definition}
\newtheorem{proposition}[theorem]{Proposition}
\newtheorem{corollary}[theorem]{Corollary}
\newtheorem{lemma}[theorem]{Lemma}
\newtheorem{rem}[theorem]{Remark}
\newtheorem{example}[theorem]{Example}

\title{On minimal bases and indices of rational matrices and their linearizations}

\author{A. Amparan\footnotemark[1], F. M. Dopico\footnotemark[2], S.
Marcaida\footnotemark[1],
I. Zaballa\footnotemark[1] \\ \\
{\small \em Dedicated to Paul Van Dooren on the occasion of his 70th birthday}}

\date{}

\begin{document}

\maketitle

\renewcommand{\thefootnote}{\fnsymbol{footnote}}

\footnotetext[1]{Departamento de Matem\'{a}tica Aplicada y EIO,
Universidad del Pa\'{\i}s Vasco UPV/EHU, Apdo. Correos 644, Bilbao 48080,
Spain. E-mail addresses: agurtzane.amparan@ehu.eus (A. Amparan),
silvia.marcaida@ehu.eus (S. Marcaida), ion.zaballa@ehu.eus (I. Zaballa).
Supported by ``Ministerio de Ciencia, Innovaci\'on y Universidades'' of Spain and ``Fondo Europeo de Desarrollo
Regional (FEDER)'' of EU through grants MTM2017-83624-P and MTM2017-90682-REDT,
and by UPV/EHU through grant GIU16/42.}

\footnotetext[2]{Departamento de Matem\'aticas,
Universidad Carlos III de Madrid, Avda. Universidad 30, 28911 Legan\'es, Spain. E-mail address:
dopico@math.uc3m.es (F. M. Dopico). Supported
by ``Ministerio de Ciencia, Innovaci\'on y Universidades'' of Spain and ``Fondo Europeo de Desarrollo Regional (FEDER)'' of EU through grants MTM2015-65798-P and MTM2017-90682-REDT {and by the ``Proyecto financiado por la Agencia Estatal de Investigaci\'on (PID2019-106362GB-I00/AEI/10.13039/501100011033)".}}

\renewcommand{\thefootnote}{\arabic{footnote}}
\begin{abstract}
  {A complete theory of the relationship between the minimal bases and indices of rational matrices and those of their strong linearizations is presented}. Such theory is based on establishing first the relationships between the minimal bases and indices of rational matrices and those of their polynomial system matrices under the classical minimality condition and certain additional conditions of properness. {This is related to pioneer results obtained by Verghese, Van Dooren and Kailath in 1979-80, which were the first proving results of this type under different nonequivalent conditions}. It is shown that the definitions of linearizations and strong linearizations do not guarantee any relationship between the minimal bases and indices of the linearizations and the rational matrices in general. {In contrast, simple relationships are obtained for the family of strong block minimal bases linearizations, which can be used to compute minimal bases and indices of any rational matrix, including rectangular ones, via algorithms for pencils. These results extend the corresponding ones for other families of linearizations available in recent literature for square rational matrices.}
  %	 is in contrast with other families of linearizations available in recent literature that only deal with square matrices. Moreover, the results obtained for strong block minimal bases linearizations allow to deduce the corresponding results for other classes of linearizations.}
\end{abstract}



%\input{00README.tex}

%---------------------------------------------------------------------------------------------
\section{Introduction}\label{sec.Intro}
%-------------------------------------------------------------------------------------------

{Rational matrices, i.e., matrices whose entries are rational functions have been studied intensively since the 1950s in connection with linear systems and control theory \cite{Kailath80,Rosen70}. The pioneering work of Rosenbrock \cite{Rosen70} showed that some properties of a rational matrix $G(\la)$ can be studied from its polynomial system matrices,
that is, polynomial matrices
of the form
\begin{equation}\label{s0.eqpolsysmat}
P(\la)=\begin{bmatrix}
A(\la) & B(\la)\\
-C(\la) & D(\la)
\end{bmatrix}
\end{equation}
with $A(\la)$ regular, i.e., $\det A(\la)\not\equiv 0$, and with transfer function matrix $G(\la)=D(\la)+C(\la)A(\la)^{-1}B(\la)$. A classical result of Rosenbrock states that if $P(\la)$ has least order (equivalently, is minimal or irreducible), i.e., $\begin{bsmallmatrix} A(\la) &B(\la)\end{bsmallmatrix}$ and $\begin{bsmallmatrix} A(\la) \\-C(\la) \end{bsmallmatrix}$ have no {\em finite} zeros, then the {\em finite pole} structure of $G(\la)$ is given by the finite zero structure of $A(\la)$ and the {\em finite zero} structure of $G(\la)$ is given by the finite zero structure of $P(\la)$. Thus, the finite pole and zero structures of a rational matrix can be obtained from any of its minimal polynomial system matrices. However, rational matrices have other structural data that are interesting in applications and that cannot be obtained, in general, from its minimal polynomial system matrices. Actually, rational matrices have {\em poles and zeros at infinity} \cite[pp. 449-450]{Kailath80} (see also \cite{silverman1983system} for applications and some classical algorithms for computing the structure at infinity). Moreover, singular, i.e., nonregular, rational matrices have {\em minimal bases and indices} \cite{Fo75}, \cite[Section 6.5.4]{Kailath80}, which form their null-space or singular structure. Minimal bases of rational matrices (in general rectangular) have a number of important applications, as, for instance, in the solution of minimal design problems, originally proposed in \cite{wang1973}, and further studied in \cite[Section 8]{Fo75} and \cite{kung1980}. In this type of problems one must solve $G_1 (\la) X(\la) = G_2 (\la)$ for a rational unknown $X(\la)$ with special properties, where $G_1 (\la)$ and $G_2 (\la)$ are given rational matrices. The solution can be constructed, when possible, from a minimal basis of $[G_1 (\la) ,\; -G_2 (\la)]$.

Motivated by the previous discussion, the pioneer work of Verghese, Van Dooren and Kailath established in \cite{VeDoka79,Ver80} sufficient conditions on a polynomial system matrix as in \eqref{s0.eqpolsysmat} that allow to obtain from $P(\la)$ the pole and zero structures at infinity and the minimal bases and indices of its transfer function matrix $G(\la)$, in addition to the finite pole and zero structures. In \cite{VeDoka79}, the authors considered polynomial system matrices in generalized state-space form (GSSF), that is, $A(\la)=\la E-A$ and $B(\la)=B$, $C(\la)=C$ and $D(\la)=D$ are constant matrices. These
polynomial system matrices were defined to be  strongly irreducible if
$\begin{bsmallmatrix}\la E-A &B\end{bsmallmatrix}$ and $\begin{bsmallmatrix}\la E-A \\-C\end{bsmallmatrix}$
have no finite or infinite zeros. Under these conditions, it was shown in \cite{VeDoka79} that the pole and zero structures at infinity and the minimal bases and indices of $P(\la)$ and $G(\la)$ are very easily related to each other. One year later, Verghese extended in \cite{Ver80} the results of \cite{VeDoka79} to general polynomial system matrices. However, the definition of strong irreducibility in this case is more complicated since it involves checking not only the finite zeros of
 $\begin{bsmallmatrix}A(\la) &B(\la)\end{bsmallmatrix}$ and $\begin{bsmallmatrix}A(\la) \\-C(\la)\end{bsmallmatrix}$
but also the infinite zeros of two extensions of $P(\la)$.

The results in \cite{VeDoka79} were fundamental for the development of efficient and stable numerical algorithms for computing the complete list of structural data of any rational matrix $G(\la)$, i.e., its finite and infinite pole and zero structures and its minimal indices. The reason is that Van Dooren presented in \cite[Section V]{vd1981} stable algorithms for constructing strongly irreducible polynomial system matrices in GSSF, starting from an easily constructible one in GSSF that is not strongly irreducible. Since a polynomial system matrix in GSSF is a pencil, the staircase algorithm in \cite{vd1979} can be applied to compute its structural data, which are related to those of $G(\la)$ as established in \cite{VeDoka79}. The minimal bases of this pencil can be computed by combining the output of \cite{vd1979} with the method in \cite[Sections 3 and 4]{beelen1987} and, finally, one can recover the minimal bases of $G(\la)$ as explained in \cite{VeDoka79}.

Despite of its importance and major impact, the results in \cite{vd1981,VeDoka79} are not always the most convenient tools for computing the complete list of structural data of a rational matrix $G(\la)$ and its minimal bases. One reason of this is that in recent years a number of new {\em linearizations} of rational matrices arising in applications have been developed for computing their structural data via algorithms for pencils and that these linearizations are neither in GSSF nor are necessarily strongly irreducible in the sense of \cite{Ver80}. In plain words, these new linearizations are minimal linear polynomial system matrices of a rational matrix $\widehat{G} (\la)$, not necessarily equal to $G(\la)$, but related to $G(\la)$. This motivates to look for some conditions that can replace Verghese's strong irreducibility conditions at infinity in \cite{Ver80}, that are useful in this new setting, and that guarantee that a minimal polynomial system matrix as in \eqref{s0.eqpolsysmat} allows to obtain the pole and zero structures at infinity and the minimal bases and indices of its transfer function matrix, in addition to the finite pole and zero structures. As a first contribution, we will prove in this paper that $C(\la) A(\la)^{-1}$ and $A(\la)^{-1} B(\la)$ being proper rational matrices are such conditions, which have the advantage of being directly checked on the polynomial system matrix itself. Moreover, we will show that these conditions do not imply Verghese's conditions and vice versa. These properness conditions were in fact introduced in \cite[Lemma 2.4 and Corollary 2.5]{AmDoMaZa18}, where it is proved that they allow to recover the pole and zero structures at infinity of a rational matrix from its polynomial system matrices satisfying these conditions. We will prove in Section \ref{Sec_PSM} that they also allow to recover the minimal indices and bases and will use this fact for proving other results in this paper.

In recent years, rational matrices have received considerable attention in the context of nonlinear eigenvalue problems (NLEPs), either because they arise directly in some of these problems or because they are used for approximating matrices of nonlinear functions in order to solve numerically NLEPs \cite{guttel-tisseur,mehrmann-voss-04}. Thus, the numerical solution of a NLEP reduces in practice to computing the finite zeros that are not poles of a rational matrix $G(\la)$. One of the most efficient approaches for doing this is to apply an eigenvalue algorithm to a linearization of $G(\la)$, i.e., a pencil that allows to recover the information of interest of $G(\la)$, in the spirit of Van Doreen's approach \cite{vd1981}. In NLEPs, this approach started in  \cite{SuBai11}, continued in \cite{nleigs}, and influenced the search for rigorous definitions of linearizations of rational matrices and the development of new families of linearizations, as it has been done before for polynomial matrices (see the references in \cite{BuDoPe18,MMMM06}).

The first formal definition of linearization of a rational matrix was given in \cite{AlBe16}. A different definition of linearization was introduced in \cite{AmDoMaZa18}, together with the first definition of strong linearization of a rational matrix (i.e., a linearization that allows to recover the infinite pole and zero structures, in addition to the finite ones). The definition of linearization in \cite{AmDoMaZa18} includes the one in \cite{AlBe16}. A comparison between both definitions can be found in \cite[Section 5]{AmDoMaZa16}. One of the differences between both definitions is that both are linear polynomial system matrices, but the off-diagonal blocks in \cite{AlBe16} are forced to be constant matrices, while this is not the case in \cite{AmDoMaZa18}. Another definition of strong linearization was introduced in \cite{DaAl19_2}, that is equivalent to the one in \cite{AmDoMaZa18}, except for the fact that the off-diagonal matrices are constant. Since the definitions mentioned above do not capture some of the pencils that have been used recently in the numerical solution of NLEPs \cite{nleigs,automatic}, the concept of local linearizations of rational matrices was introduced in \cite{DoMaQuVD19}.

Simultaneously to the definitions of linearizations of rational matrices, different families of linearizations that can be constructed without any numerical operation have been developed. Some of them are restricted to square rational matrices, as the different Fiedler-like linearizations presented in \cite{AlBe16,AlBe18,DasAl19laa,DaAl20}, and the affine-spaces and $\mathbb{M}_1$ and  $\mathbb{M}_2$-strong linearizations introduced in \cite{DaAl19_2} and in \cite{DoMaQu19}, respectively, which are closely connected to each other. In contrast, the family of strong block minimal bases linearizations introduced in \cite{AmDoMaZa18} is valid for general rational matrices, i.e., rectangular or square. An important property of this wide family is that it contains modulo permutations the Fiedler-like linearizations mentioned above. This is proved in Section \ref{sec.Fiedler} and extends to rational matrices a well-known result for polynomial matrices \cite{BuDoPe18}. Moreover, it is known that the $\mathbb{M}_1$ and  $\mathbb{M}_2$-strong linearizations are very simply connected to strong block minimal bases linearizations \cite{DoMaQu19}. Thus, the properties of strong block minimal bases linearizations imply results for other families of linearizations and for more general rational matrices.

The main focus of this paper is on the relationship of the minimal bases and indices of a rational matrix and its linearizations. Previous works on this problem deal with strongly irreducible pencils in GSSF \cite{VeDoka79}, strongly irreducible polynomial system matrices when they are linear \cite{Ver80} (the pencils in \cite{VeDoka79,Ver80} are particular cases of linearizarizations in the sense of \cite{AmDoMaZa18}, though not strong) and, recently, Fiedler-like linearizations \cite{DasAl19laa,DaAl20} and affine spaces of linearizations \cite{DaAl19_2}. The results in \cite{DasAl19laa,DaAl19_2,DaAl20} have the advantage with respect to those in \cite{VeDoka79,Ver80} of dealing with linearizations that can be constructed without computational cost from certain representations of the rational matrix. However, the results in \cite{DasAl19laa,DaAl19_2,DaAl20} have a disadvantage for computing minimal indices and bases, since they only work for square rational matrices while the applications where minimal indices and bases are of interest very often include rectangular matrices. This is the case, for example, of transfer function matrices of systems with different numbers of inputs and outputs and in the minimal design problems mentioned at the beginning of this section. In this setting, a second contribution of this paper is to establish for the first time simple relationships between the minimal bases and indices of a rational matrix and its strong block minimal bases linearizations. The proof of these results is presented in Section \ref{Sec_blockminimalbases} based on those in Section \ref{Sec_PSM} and allow to use strong block minimal bases linearizations for computing minimal bases and indices via algorithms for pencils \cite{vd1979,beelen1987}. Since many strong block minimal bases linearizations can be constructed without computational cost and are valid for rectangular matrices, this approach is advantageous with respect to the previous ones. In addition, due to the relation of this family of linearizations with other families of linearizations, the results in Section \ref{Sec_blockminimalbases} allow to obtain as corollaries the corresponding ones for Fiedler-like pencils, already presented in \cite{DasAl19laa,DaAl20}, in Subsection \ref{subsec.consec} and the ones for $\mathbb{M}_1$ and  $\mathbb{M}_2$-strong linearizations in Section \ref{Sec_M12minimalbases}. This approach constitutes a unified treatment of the relationships between minimal bases and indices of rational matrices and those of most classes of linearizations developed in the recent literature.

The paper is completed with the study in Sections \ref{Sec_linearizations} and \ref{Sec_minimalindicesstrong} of a question of a fundamental nature: what is the relationship between the minimal bases and indices of a rational matrix and those of any of its linearizations or strong linearizations? Here and in what follows by linearizations and strong linearizations we will mean those defined in \cite{AmDoMaZa18}. For linearizations, our conclusions are that we can recover polynomial bases, but not
minimal ones, and that the dimension of the left (right, respectively) nullspace of any rational matrix and that of its
linearizations coincide. Using  Van Dooren's index sum theorem  \cite{VeDoka79},  we can obtain the sums of the right and left minimal indices of any rational matrix
from any of its strong linearizations. However, as far as the minimal indices themselves are concerned, we will
show that  the minimal indices of a rational matrix may differ arbitrarily from those of its strong linearizations.

The paper is organized as follows: Section \ref{sec.prelim} contains the notation and some preliminary results. In Section \ref{Sec_PSM} the relationship between the minimal bases and indices of a rational
matrix and its minimal polynomial system matrices satisfying some conditions of properness is given.
Section \ref{Sec_linearizations} relates polynomial bases of rational matrices and their linearizations in general.
Section \ref{Sec_minimalindicesstrong} shows that the minimal indices of a rational matrix and of its strong linearizations
may differ arbitrarily, but that there is a connection between the sums of their left and right minimal indices.
Section \ref{Sec_blockminimalbases} is devoted to obtain minimal bases and indices of any rational matrix from its strong
block minimal bases linearizations and vice versa. The same goal is pursued in Section \ref{Sec_M12minimalbases} for $\mathbb{M}_1$ and
$\mathbb{M}_2$-strong linearizations. Connections between strong block minimal bases linearizations and Fiedler-like linearizations are investigated in Section \ref{sec.Fiedler}. Some brief remarks on eigenvectors are included for completeness in Section \ref{sec.eigenvectors}. Finally, some conclusions are discussed in Section \ref{Sec.Conclusions}.
}

%---------------------------------------------------------------------------------------------
\section{Preliminaries}\label{sec.prelim}
%-------------------------------------------------------------------------------------------

 Most of the results included in this section are classic and can be found in standard references as \cite{Kailath80,Rosen70,Vard91}, together with much more information on rational matrices.

Let
$\FF$ be an arbitrary field and $\overline{\FF}$ its algebraic closure. Let $\FF[\la]$ be the ring of polynomials with
coefficients in $\FF$ and $\FF(\la)$ the field of fractions of
$\FF[\la]$, i.e., the field of rational functions over $\FF$.
%Denote by $\overline{\FF}$ the algebraic closure of $\FF$.
%If
%$t(\la)=\frac{p(\la)}{q(\la)}\in\FF(s)$ where $p(\la), q(\la)\in\FF[\la]$,
%$q(\la)\neq 0$, we consider the map
%$\delta_\infty(\cdot):{\efe}(\la)\longrightarrow{\mathbb{Z}}\cup\{+\infty\}$
%defined via:
%$$
%\delta_\infty(t(\la))=\begin{cases}
%\deg(q(\la))-\deg(p(\la)), &\text{if $t(\la)\neq 0$,}\\
%+\infty, &\text{if $t(\la)=0,$}
%\end{cases}
%$$
%where $\deg(\cdot)$ stands for ``degree of''.
The elements of $\FF(\la)$ with the degree of the numerator at most the degree of the denominator are called proper rational functions. The set of proper rational functions over $\efe$ form a ring denoted $\FF_{pr}(\la)$.
%, that is, let
% $\FF_{pr}(\la)=\{t(\la)
%\in {\efe}(\la): \delta_\infty(t(\la))\geq0\}$ denotes
If the degree of the denominator of a rational function is strictly larger than the degree of its numerator then the rational function is called strictly proper.
%If  $\delta_\infty(t(\la))>0$ then
%$t(\la)$ is called strictly proper. For all non-zero elements
%$t_1(\la), t_2(\la)$ of $\FF(\la)$, the map $\delta_\infty(\cdot)$
%satisfies the following two properties (see \cite[Sec. 3.2]{Vard91})

%\begin{equation}\label{valpro}
%\delta_\infty(t_1(\la)t_2(\la))=\delta_\infty(t_1(\la))+\delta_\infty(t_2(\la)),
%\end{equation}

%\begin{equation}\label{valsum}
%\delta_\infty(t_1(\la)+t_2(\la))\geq \min\{\delta_\infty(t_1(\la)),\delta_\infty(t_2(\la))\}.
%\end{equation}

Vectors with entries in $\FF[\la]$ are called vector polynomials.  $\FF(\la)^p$ stands for the vector space  of $p$-tuples
of rational functions.
We denote by $\efe[\la]^{p\times m}$ (resp., $\efe(\la)^{p\times m}$, $\efe_{pr}(\la)^{p\times m}$) the set of $p\times m$ matrices with entries in
$\efe[\la]$ (resp., $\efe(\la)$, $\efe_{pr}(\la)$). Matrices in $\efe[\la]^{p\times m}$ are called polynomial matrices or matrix polynomials indistinctly. The degree of a polynomial matrix is the highest degree of all its entries. The square polynomial matrices whose inverses are polynomial matrices are called unimodular matrices. Matrices  in $\efe(\la)^{p\times m}$ are known as rational matrices and
matrices  with entries in $\efe_{pr}(\la)$ are termed as proper rational matrices. In particular, if the entries are all strictly proper then they are called strictly proper rational matrices. Invertible matrices in $\efe_{pr}(\la)^{p\times p}$, that is, square proper rational matrices whose inverses are also proper, are called biproper.  Equivalently, biproper matrices are square proper rational matrices whose determinants are biproprer rational functions. { $I_n$ denotes the $n\times n$ identity matrix and $0_{p\times m}$ (or $0_{p,m}$) the $p\times m$ zero matrix, which reduces to $0_p$ when $p=m$. We will write just $I$ or $0$ when the size is clear from the context.}

%____________________________________
{\color{black}\subsection{Spectral structure of rational matrices}}
%________________________________________

We introduce now the spectral structure (both finite and infinite) of rational matrices. Recall that two rational matrices $G_1(\la), G_2(\la)\in\efe(\la)^{p\times m}$ are unimodularly equivalent if there exist unimodular matrices $U_1(\la)\in\efe[\la]^{p\times p}$ and $U_2(\la)\in\efe[\la]^{m\times m}$ such that $G_2(\la)=U_1(\la)G_1(\la)U_2(\la)$. Any rational matrix is unimodularly equivalent
to its finite Smith--McMillan form (see, for example, \cite[Chapter 3, Section 4]{Rosen70}
or \cite[Section 6.5.2]{Kailath80}). That is to say, if
$G(\la)\in\FF(\la)^{p\times m}$ then there are unimodular matrices
$U_1(\la)\in\FF[\la]^{p\times p}$ and $U_2(\la)\in\FF[\la]^{m\times m}$
such that
\begin{equation}\label{s1.eqMcMillan}
U_1(\la)G(\la)U_2(\la)=\left[\begin{array}{cc}\diag\left(\frac{\epsilon_1(\la)}{\psi_1(\la)},\ldots,
\frac{\epsilon_r(\la)}{\psi_r(\la)}\right)&0\\0&0\end{array}\right]
\end{equation}
where  $\epsilon_1(\la),\ldots,\epsilon_r(\la),\psi_1(\la),\ldots,\psi_r(\la)$ are nonzero
monic (leading coefficient equal to 1) polynomials, $\epsilon_i(\la), \psi_i(\la)$ are pairwise coprime for all $i=1,\ldots,r$, and
$\epsilon_1(\la)\mid\cdots\mid\epsilon_r(\la)$ while $\psi_r(\la)\mid\cdots\mid\psi_1(\la)$,
where $\mid$ stands for divisibility. {\color{black}The index $r$ is the normal rank of $G(\la)$, which we denote by $\rank G(\la)$.}
The finite zeros of $G(\la)$
are the roots in $\overline{\FF}$  of
$\epsilon_r(\la)$ and its finite poles
are the roots in  $\overline{\FF}$ of $\psi_1(\la)$.
If $\la_0\in\overline{\FF}$ is a zero of $G(\la)$ then, for $i=1,\ldots, r$, we can write
$\epsilon_i(\la)=(\la-\la_0)^{m_i}\wh{\epsilon}_i(\la)$ with $\wh{\epsilon}_i(\la_0)\neq 0$
and $m_i\geq 0$. The nonzero elements in $(m_1,\ldots, m_r)$ are called the
partial multiplicities of $\la_0$ as a zero of $G(\la)$. In the same way, if $\la_0\in\overline{\FF}$ is a pole of $G(\la)$ then, for $i=1,\ldots, r$, we can write
$\psi_i(\la)=(\la-\la_0)^{n_i}\wh{\psi}_i(\la)$ with $\wh{\psi}_i(\la_0)\neq 0$
and $n_i\geq 0$. The nonzero elements in $(n_1,\ldots, n_r)$ are called the
partial multiplicities of $\la_0$ as a pole of $G(\la)$.
We understand by finite zero structure of $G(\la)$ its finite zeros together with their respective partial multiplicities. Analogously, the finite pole structure of $G(\la)$ consists of its finite poles each with its partial multiplicities.

Rational matrices may have structure at infinity as well. Recall (see, for example, \cite{Vard91}) that two rational matrices of the same size $G_1(\la), G_2(\la)\in\efe(\la)^{p\times m}$ are equivalent at infinity if there exist biproper matrices $B_1(\la)\in\efe_{pr}(\la)^{p\times p}$ and $B_2(\la)\in\efe_{pr}(\la)^{m\times m}$ such that $G_2(\la)=B_1(\la)G_1(\la)B_2(\la)$. Any rational matrix is equivalent at infinity
to its Smith--McMillan form at infinity. That is to say, if
$G(\la)\in\FF(\la)^{p\times m}$ then there are biproper matrices
$B_1(\la)\in\FF_{pr}(\la)^{p\times p}$ and $B_2(\la)\in\FF_{pr}(\la)^{m\times m}$
such that
$$ %\begin{equation}
B_1(\la)G(\la)B_2(\la)=\left[\begin{array}{cc}\diag\left((\frac{1}{\la})^{q_1},\ldots,
(\frac{1}{\la})^{q_r}\right)&0\\0&0\end{array}\right]
$$ %\end{equation}
where $r=\rank G(\la)$ and $q_1\leq\cdots\leq q_r$ are integers. These are called the invariant orders at infinity of $G(\la)$.
They determine the zeros and poles at infinity of $G(\la)$, also called infinite zeros and poles. Namely, if $q_1\leq\cdots\leq q_k<0=q_{k+1}=\cdots=q_{u-1}<q_u\leq\cdots\leq q_r$ are the invariant orders at infinity of $G(\la)$ then $G(\la)$ has $r-u+1$ zeros at
infinity each one of order $q_u,\ldots,q_r$ and $k$ poles at infinity each one of order $-q_k,\ldots,-q_1$.
%Thus, the total number of zeros at infinity is $\sum_{q_i>0}q_i$ while the total number of poles at infinity is $\sum_{q_i<0}(-q_i)$.
Notice that proper rational matrices have all its invariant orders at infinity nonnegative, that is, they do not have poles at infinity. Moreover, all the
invariant orders at infinity of strictly proper rational matrices are positive.
%	 if $G(\la)$ is proper then all its invariant orders at infinity are nonnegative. Moreover, if $G(\la)$ is strictly proper then all its invariant orders at infinity are positive.

Note that any rational matrix can be decomposed uniquely as
$
G(\la)=D(\la)+G_{sp}(\la)
$
with $D(\la)$ a polynomial matrix and $G_{sp}(\la)$ a strictly proper matrix. When $G(\la)$ is not strictly proper, that is, when $D(\la)\neq 0$, the first invariant order at infinity of $G(\la)$, $q_1$, turns out to be minus the degree of the polynomial part of $G(\la)$, i.e., $q_1=-\deg(D(\la))$ (see \cite[Section 2]{AmDoMaZa18}),  where $\deg(\cdot)$ stands for ``degree of''.

%\bigskip
%_______________________________________________________
\subsection{Polynomial system matrices}
%_________________________________________________________

Any rational matrix $G(\la)\in\efe(\la)^{p\times m}$ can be written as $G(\la)=D(\la)+C(\la)A( \la)^{-1}B(\la)$ where $A(\la)\in\FF[\la]^ {n\times n}$ is {\color{black}regular, i.e., $\det A(\la)\not\equiv 0$}, $B(\la)\in\FF[\la]^{n\times m}$, $C(\la)\in\FF[\la]^{p\times n}$ and $D(\la)\in\FF[\la]^{p\times m}$. %and $n\geq \deg(\det A(\la))$.
The polynomial matrix formed with these matrices
\begin{equation}\label{s1.eqpolsysmat}
P(\la)=\begin{bmatrix}
A(\la) & B(\la)\\
-C(\la) & D(\la)
\end{bmatrix}
\end{equation}
is called a polynomial system matrix of (or giving rise to)
$G(\la)$ (see \cite{Rosen70}). The rational matrix
$G(\la)$ is called the transfer function matrix of
$P(\la)$ and  $\deg(\det A(\la))$ is known as the order of
$P(\la)$.
We allow $n$ to be equal to 0 in the definition of polynomial system matrix. In this case we say that $P(\la)=D(\la)$ is a polynomial system matrix giving rise to $G(\la)=D(\la)$, that is, $A(\la)$, $B(\la)$ and $C(\la)$ are empty matrices.
Besides,  when $A(\la)$ is a monic linear matrix
polynomial, say $A(\la)=\la I_n -A$, and $B(\la)=B$ and $C(\la)=C$ are constant matrices,
$P(\la)$ is said to be a polynomial system matrix of $G(\la)$ in
state-space form.

Different polynomial system matrices may exist with different orders
giving rise to the same transfer function matrix. A polynomial system matrix of $G(\la)$ is said to have least order, or
to be minimal, if its order is the smallest integer for which matrix polynomials $A(\la)$ ({\color{black}regular}),
%$n\geq \deg(\det A(\la))$)
$B(\la)$, $C(\la)$ and $D(\la)$ giving rise to  $G(\la)=D(\la)+C(\la)A( \la)^{-1}B(\la)$ exist (\cite[Chapter 3, Section 5.1]{Rosen70}
or \cite[Section 1.10]{Vard91}). In consequence, associated with any rational matrix $G(\la)$ there is a unique least order,
which is the order of any minimal polynomial system matrix giving rise to $G(\la)$, and is denoted by $\nu(G(\la))$. Interested
readers can find  in \cite[Chapter 3, Section 5.1]{Rosen70} three algorithms to compute $\nu(G(\la))$ without going to the
length of finding a least order polynomial system matrix giving  rise to $G(\la)$.

One of the many characterizations of when a polynomial system matrix has least order is given in terms of coprimeness.
Two polynomial matrices
$A(\la)\in\FF[\la]^{p\times m}$, $C(\la)\in\FF[\la]^{q\times m}$ with $p+q\geq m$ are
called right coprime if their only right common divisors
are unimodular matrices.
That is to say, if there exist $\wh{A}(\la)\in\FF[\la]^{p\times m}, \wh{C}(\la)\in\FF[\la]^{q\times
m},
X(\la)\in\FF[\la]^{m\times m}$ such that $A(\la)=\wh{A}(\la) X(\la)$ and  $C(\la)=\wh{C}(\la)X(\la)$, then $X(\la)$ is unimodular.
%\[
%\begin{array}{l}
%A(\la)=\bar{A}(\la) X(\la)\\
%C(\la)=\bar{C}(\la)X(\la)
%\end{array}
%%(\wh{A}(\la)\in\FF[\la]^{p\times n}, %\wh{B}(\la)\in\FF[\la]^{q\times n},  %X(\la)\in\FF[\la]^{n\times %n})
%\Rightarrow X(\la) \text{  unimodular}.
%\]
Let us recall some equivalent conditions that characterize when two polynomial matrices
are right coprime   (see, for example, \cite[Chapter 2, Section 6]{Rosen70}, \cite[Chapter 1]{Vard91}, \cite{AmMaZa16}):

\begin{proposition}\label{prop_coprim}
Let $A(\la)\in\FF[\la]^{p\times m}$and $C(\la)\in\FF[\la]^{q\times m}$ with $p+q\geq m$. The following conditions are equivalent:
\begin{itemize}
    \item[\rm{(i)}] $A(\la)$ and $C(\la)$ are right coprime.
%    \item[\rm{(ii)}] The Smith normal form of $\begin{bmatrix} A(\la) \\ C(\la) \end{bmatrix}$ is $\begin{bmatrix} I_m \\ 0 \end{bmatrix}$.
    %\item[(ii)] There exists unimodular matrices  $U(\la)\in\FF[\la]^{m\times m} $,  $V(\la)\in\FF[\la]^{(p+q)\times (p+q)} $ such that $$U(\la)\left[\begin{array}{cc}A(\la) & B(\la)\end{array}\right]V(\la)=\left[\begin{array}{cc}I_m & 0\end{array}\right]$$
    \item[\rm{(ii)}] There exist matrices  $X(\la)\in\FF[\la]^{m\times p}$, $Y(\la)\in\FF[\la]^{m\times q}$ such that $X(\la)A(\la)+Y(\la)C(\la)=I_m$.
    \item[\rm{(iii)}] $\rank \begin{bsmallmatrix} A(\la_0) \\ C(\la_0) \end{bsmallmatrix}=m$ for all $\la_0\in \overline{\FF}$.
    %\item[(iii)]  $\rank \left[\begin{array}{cc}A(\la_0) & B(\la_0)\end{array}\right]=m$ for all $\la_0\in \overline{\FF}$
\end{itemize}
\end{proposition}

On the other hand, $A(\la)\in\FF[\la]^{m\times p}$ and $C(\la)\in\FF[\la]^{m\times q}$, $p+q\ge m$, are
left coprime if their transposes $A(\la)^T$ and $C(\la)^T$ are right coprime.

It turns out that the polynomial system matrix in (\ref{s1.eqpolsysmat}) has least order if and only if $A(\la)$ and $B(\la)$ are left coprime and $A(\la)$ and
$C(\la)$ are right coprime (\cite[Chapter 3]{Rosen70}).

A celebrated result by Rosenbrock \cite[Chapter 3, Theorem 4.1]{Rosen70} relates the finite structure (zero and pole structure) of a rational matrix with the finite structure of its minimal polynomial system matrices. Namely,
 when the polynomial system  matrix in (\ref{s1.eqpolsysmat}) giving rise to $G(\la)$
has least order, the finite zero structure of $G(\la)$ is the finite zero structure of $P(\la)$  and the finite pole structure of $G(\la)$ is the finite zero structure
of $A(\la)$. A consequence of this fact is that the least order of $G(\la)$, $\nu(G(\la))$, which is the degree of the determinant of $A(\la)$, is equal to the sum of the partial multiplicities of the finite poles of $G(\la)$.  In other words, $\nu(G(\la))$ is the sum of the degrees of the denominators in the finite Smith--McMillan form of $G(\la)$.

{\color{black}
%____________________________________
\subsection{Singular structure of rational matrices}
%________________________________________
}

Let us introduce now the singular structure of a rational matrix. Denote by $\mathcal{N}_\ell (G(\la))$ and $\mathcal{N}_r (G(\la))$ the left and
right null-spaces over $\FF(\la)$ of $G(\la)$, respectively, i.e.,
if $G(\la)\in\FF(\la)^{p\times m}$,
\[
\begin{array}{l}
\mathcal{N}_\ell (G(\la))=\{x(\la)\in\FF(\la)^{p}: x(\la)^TG(\la)=0\},\\
\mathcal{N}_r (G(\la))=\{x(\la)\in\FF(\la)^{m}: G(\la)x(\la)=0\}.
\end{array}
\]
These sets are vector subspaces
of $\FF(\la)^p$ and $\FF(\la)^m$, respectively. For any subspace of
$\FF(\la)^p$, it is always possible to find a basis consisting of vector
polynomials; simply take an arbitrary basis and multiply each vector
by the  least common multiple of the denominators of its  entries.
The order of a polynomial basis is defined as the sum of
the degrees of its vectors  (see \cite{Fo75}). If $\mathcal{V}$ is a
subspace of $\FF(\la)^p$, a minimal basis of  $\mathcal{V}$
is  a polynomial basis of  $\mathcal{V}$ with least order among all
polynomial bases of $\mathcal{V}$. The fundamental result in this setting is that the non-decreasing ordered list of degrees of
the vector polynomials in any minimal basis of $\mathcal{V}$  is  always the same (see \cite{Fo75}).  These degrees are called the minimal indices of $\mathcal{V}$.

%A right  (resp., left) minimal bases of a rational matrix $G(\la)$ is defined as a minimal basis of the rational subspace $\mathcal{N}_r (G(\la))$ (resp., $\mathcal{N}_\ell (G(\la))$).
We refer to a polynomial matrix $N(\la)\in\efe[\la]^{m\times l}$ itself as a right polynomial basis of  a rational matrix $G(\la)\in\efe(\la)^{p\times m}$ if the columns of $N(\la)$ form a basis of
$\mathcal{N}_r (G(\la))$. If the  columns of $N(\la)$ form a minimal basis of $\mathcal{N}_r (G(\la))$ then $N(\la)$ is referred  to as a right minimal basis of $G(\la)$. Notice that $l=\dim \mathcal{N}_r (G(\la))\leq m$. Moreover, $l=m$ if and only if $G(\la)=0$.

Analogously, a polynomial matrix $N(\la)\in\efe[\la]^{p\times q}$ is a left polynomial (resp., minimal) basis of a rational matrix $G(\la)\in\efe(\la)^{p\times m}$ if the columns of $N(\la)$ form a polynomial (resp., minimal) basis of $\mathcal{N}_{\ell} (G(\la))$. As above, $q=\dim \mathcal{N}_{\ell}(G(\la))\leq p$, and $q=p$ if and only if $G(\la)=0$.

The right  (resp., left) minimal indices of a rational matrix $G(\la)$ are the minimal indices of
$\mathcal{N}_r (G(\la))$ (resp., $\mathcal{N}_\ell (G(\la))$). If $N(\la)$ is a right (resp., left) minimal basis of $G(\la)$ then the right (resp., left) minimal indices of $G(\la)$ are the degrees of the columns of $N(\la)$ when ordered non-decreasingly.

One of the most usual characterizations of minimal bases is a slightly modified version of the Main Theorem given in \cite{Fo75}, which can be also found in \cite[Theorem 2.14]{DeDoVa15}. Before presenting this theorem let us recall what a column proper or column reduced matrix is. Let $N(\la)\in \FF[\la]^{m \times l}$. We denote by $\deg(\col_j(N(\la)))$ the degree of the $j$-th column of $N(\la)$, that is, the degree of the highest degree entry in column $j$. Put $d_j=\deg(\col_j(N(\la)))$. The matrix $N(\la)$ can always be written  (see \cite[Section 6.3.2]{Kailath80}) as
\begin{equation}\label{reduced}
N(\la)=N_h \diag(\la^{d_1}, \ldots, \la^{d_l})+ L(\la)
\end{equation}
where  $N_h$ is the  highest column degree coefficient matrix of $N(\la)$, and $L(\la)$ is a polynomial matrix collecting the remaining terms,
which
 has lower column degrees than the corresponding ones of $N(\la)$. The polynomial  matrix $N(\la)$ is called  column proper or   column reduced if
 $\rank N_h=l$.

\begin{theorem}{({\cite[Theorem 2.14]{DeDoVa15}})}\label{mincharact}
The columns of a matrix polynomial $N(\la)$ over a field $\efe$
%\mnote{Este resultado est\'a tambi\'en en \cite{VeDoka79} when $\efe=\CC$. Por otra parte %yo pondr�a $N(\la)$ over the field $\FF(s)$}
are a minimal basis of the subspace they span if and only if
$N(\la_0)$ has full column rank for all $\la_0\in\overline{\FF}$ and $N(\la)$ is column reduced.
\end{theorem}

\section{Minimal bases and indices of polynomial system matrices}\label{Sec_PSM}
%-------------------------------------------------------------------------------------------

{The goal of this section is to show that if the blocks of a minimal polynomial system matrix $P(\la)$ of a rational matrix $G(\la)$ satisfy certain properness conditions, then the right (resp. left) minimal indices of $P(\la)$ and $G(\la)$ are the same and the right (resp. left) minimal bases of $P(\la)$ and $G(\la)$ are easily related to each other through a simple bijection. This is stated in Corollary \ref{col_PGPmin}, which is the final result in this section and is obtained as a consequence of a number of intermediate results that deal independently with the left and the right minimal indices and bases. The properness conditions mentioned above are satisfied, in particular, by many  classes of linearizations of rational matrices introduced recently \cite{AlBe16,AlBe18,AmDoMaZa18,DasAl19laa,DaAl19_2,DaAl20,DoMaQu19}, which are minimal linear polynomial system matrices of some rational matrices $\widehat{G} (\la)$ different from $G(\la)$, in general, but related to $G(\la)$. The results of this section are used in Section \ref{Sec_blockminimalbases} to establish the relationships between the minimal indices and bases of $G(\la)$ and those of its strong block minimal bases linearizations, which in turn imply the corresponding results for $\mathbb{M}_1$ and $\mathbb{M}_2$-strong linearizations and Fiedler-like linearizations of $G(\la)$ in Sections \ref{Sec_M12minimalbases} and \ref{sec.Fiedler}. This provides a unified treatment of all these results and alternative proofs to the results in \cite{DasAl19laa,DaAl20}.

The properness conditions used in this section have been introduced in \cite[Lemma 2.4 and Corollary 2.5]{AmDoMaZa18} with the purpose of obtaining the invariant orders at infinity of $G(\la)$ from those of $P(\la)$. The new result is that they also allow to recover the minimal indices and bases of $G(\la)$ from those of $P(\la)$. Thus the coprimeness conditions of the blocks of $P(\la)$ guaranteeing its minimality  and the properness conditions together allow to recover the complete finite eigenstructure of $G(\la)$ from $P(\la)$, as a consequence of the classical result of Rosenbrock \cite{Rosen70} (see the comments following Proposition \ref{prop_coprim}), the complete infinite eigenstructure of $G(\la)$ from $P(\la)$, as a consequence of the results in \cite{AmDoMaZa18}, and the minimal bases and indices of $G(\la)$ from $P(\la)$, as a consequence of the results in this section.

As outlined in the introduction, another set of conditions on a polynomial system matrix $P(\la)$ that allow to recover the complete finite and infinite spectral structures and the minimal indices and bases of its transfer function matrix $G(\la)$ from $P(\la)$ is the strong irreducibility introduced in \cite{Ver80}. The polynomial system matrix $P(\la)$ in \eqref{s1.eqpolsysmat} is strongly irreducible if $\begin{bsmallmatrix} A(\la) & B(\la)
\end{bsmallmatrix}$ and $\begin{bsmallmatrix} A(\la) \\ -C(\la)\end{bsmallmatrix}$ have no finite zeros, and
$\begin{bsmallmatrix} A(\la)&B(\la)  &0\\ -C(\la)&D(\la)&-I_p\end{bsmallmatrix}$ and
$\begin{bsmallmatrix} A(\la)&B(\la)\\ -C(\la)&D(\la)\\0&I_m\end{bsmallmatrix}$
have no infinite zeros\footnote{In the case $B(\la), C(\la)$ and $D(\la)$ are constant matrices these conditions are equivalent to the strong irreducibility originally introduced in \cite{VeDoka79}.}.
The former two conditions are equivalent to  $A(\la)$ and $B(\la)$ be left coprime and $A(\la)$ and $C(\la)$ be right coprime, i.e., the classical Rosenbrock's minimality of $P(\la)$ (see Proposition
\ref{prop_coprim}). And the latter two conditions are used in \cite{Ver80} to relate the invariant orders at infinity and the minimal bases and indices of $P(\la)$ and $G(\la)$. We substitute in this section the
conditions on the infinite zeros by the easier to check conditions that $C(\la)A(\la)^{-1}$ and $A(\la)^{-1}B(\la)$ are proper rational functions. However, we remark that these pairs of conditions, the one on the infinite zeros and the one on properness, are not equivalent, and, even more, that none of them implies the other one, as we show at the end of this section via two examples. This is consistent with the fact that the rules for recovering the invariant orders at infinity are different under the properness conditions \cite[Corollary 2.5]{AmDoMaZa18} that under the strong irreducibility conditions \cite[Result 1]{Ver80}.
}

 %guarantee that the invariant orders at infinity of a polynomial system matrix allow us to recover those of its transfer function matrix, according to \cite[Lemma 2.4 and Corollary 2.5]{AmDoMaZa18}.
%Moreover, recall also that the coprimeness conditions allow us to recover the zero and pole finite structure of
% the transfer function from the polynomial system matrix, as we have explained after Proposition \ref{prop_coprim}.

%{\blue Before getting into the details, we want to acknowledge that the results of this section are inspired in pioneer results of \mnote{Yo quitar\'ia este par\'agrafo porque ya se ha mencionado en la introducci�n de manera muy relevante} Paul Van Dooren and co-workers (see \cite[Theorem 2]{VeDoka79} and \cite[Result 2]{Ver80}), though the properness assumptions used here do not appear in \cite{Ver80,VeDoka79}. Instead, a condition called `strong irreducibility' of the polynomial system matrix is employed in \cite{Ver80,VeDoka79}. }

The proofs of the following two lemmas follow the same pattern as the first part of the proofs of Theorem 2 in
\cite{VeDoka79} and Result 2 in \cite{Ver80} and they are omitted.

\begin{lemma}\label{rankPG}
    Let $P(\la)$ of \eqref{s1.eqpolsysmat}
%   \begin{equation}
 %   $
 %   P(\la)=\begin{bmatrix}
%    A(\la) & B(\la)\\-C(\la) &D(\la)
%    \end{bmatrix}\in\FF[\la]^{(n+p)\times (n+m)}
 %   $
%   \end{equation}
be a polynomial system matrix
%whose transfer function rational matrix is
    of a rational matrix $G(\la)$. Then $\rank P(\la)=n+\rank G(\la)$, $\dim \mathcal{N}_\ell (G(\la))=\dim \mathcal{N}_\ell (P(\la))$ and
    $\dim \mathcal{N}_r (G(\la))=\dim \mathcal{N}_r (P(\la))$.
\end{lemma}

%\textbf{Proof}.- Notice that
%$$
%\begin{bmatrix}I_n&0\\C(\la)A(\la)^{-1}&I_p\end{bmatrix}P(\la)\begin{bmatrix}I_n&-A(\la)^{-1}B(\la)\\0&I_m\end{bmatrix}=
%\begin{bmatrix}
%A(\la)&0\\0&G(\la)
%\end{bmatrix}.
%$$
%Thus, $\rank P(\la)=\rank \begin{bmatrix}
%A(\la)&0\\0&G(\la)
%\end{bmatrix}=n+\rank G(\la)$. Now, using the rank-nullity theorem, $\dim \mathcal{N}_\ell (G(\la))=p-\rank G(\la)=p+n-\rank P(\la)=\dim \mathcal{N}_\ell (P(\la))$ and $\dim \mathcal{N}_r (G(\la))=m-\rank G(\la)=m+n-\rank P(\la)=\dim \mathcal{N}_r (P(\la))$.
%\hfill\halmos

%\bigskip
\begin{lemma}\label{basisPG}
Under the assumptions of Lemma \ref{rankPG},
%     Let
%   \begin{equation}
%     $
%     P(\la)=\begin{bsmallmatrix}
%     A(\la) & B(\la)\\-C(\la) &D(\la)
%     \end{bsmallmatrix}\in\FF[\la]^{(n+p)\times (n+m)}
%     $
%   \end{equation}
 %    be a polynomial system matrix
%whose transfer function rational matrix is
%     of  a rational matrix $G(\la)$.
 if  $
    \begin{bsmallmatrix}
H_1(\la) \\H_2(\la)
    \end{bsmallmatrix}
%    \in\FF[\la]^{(n+m)\times l}
    $ is a right polynomial basis of $P(\la)$ then $H_2(\la)$ is a right polynomial basis of $G(\la)$ and $H_1(\la)=-A(\la)^{-1}  B(\la)  H_2(\la)$.
\end{lemma}

%\textbf{Proof}.- {\blue If  $\label{}
%    \begin{bmatrix}
%H_1(\la) \\H_2(\la)
%    \end{bmatrix}
%    \in\FF[\la]^{(n+m)\times l}
%    $ is a right polynomial basis of $P(\la)$ then}
%$P(\la)\begin{bmatrix}
%H_1(\la) \\H_2(\la)
%\end{bmatrix}=0$.
%Then
%$$
%\begin{bmatrix}
%I_n & 0\\C(\la)A(\la)^{-1} & I_p
%\end{bmatrix}\begin{bmatrix}
%A(\la) & B(\la)\\-C(\la) &D(\la)
%\end{bmatrix}\begin{bmatrix}
%H_1(\la) \\H_2(\la)
%\end{bmatrix}=0.
%$$
%We get, via a direct multiplication, that $A(\la)H_1(\la)+B(\la)H_2(\la)=0$ (or  $H_1(\la)=-A(\la)^{-1} B(\la)  H_2(\la)$) and
%$G(\la)H_2(\la)=0$. Let us see that the columns of $H_2(\la)$ form a right polynomial basis of $G(\la)$. Suppose that
%the columns of $H_2(\la)$ are Linearly dependent, that is
%there exists a  vector $v(\la)\neq 0$
% such that $H_2(\la)v(\la)=0$. Then $\begin{bmatrix}
% H_1(\la) \\H_2(\la)
% \end{bmatrix}v(\la)=\begin{bmatrix}
%-A(\la)^{-1}    B(\la)  H_2(\la)  \\H_2(\la)
% \end{bmatrix}v(\la)=0$, which is in contradiction {\blue with
% $\label{}
 %   \begin{bmatrix}
%H_1(\la) \\H_2(\la)
%    \end{bmatrix}
%    $ being a right polynomial basis of $P(\la)$}. Thus, the columns of $H_2(\la)$ are linearly independent in $\efe(\la)^m$. By Lemma \ref{rankPG}, $\dim \mathcal{N}_r (G(\la))=\dim \mathcal{N}_r (P(\la))$. Hence, $H_2(\la)$ is a right polynomial basis of $G(\la)$.
% \hfill\halmos

%\bigskip
{Lemma \ref{rankPG} means that any polynomial system matrix and its transfer function have the same number of right minimal
indices and the same number of left minimal indices. In turns, Lemma \ref{basisPG} shows how to obtain a right polynomial basis
of a rational matrix from a right polynomial basis of any of its polynomial system matrices. The extension of this result to right
minimal bases under coprimeness and properness conditions requires some technical preliminary lemmas.}

\begin{lemma}\label{coldeg}
Let $N_1(\la) \in\FF[\la]^{n\times l}$ and   $N_2(\la)  \in\FF[\la]^{m\times l}$.
\begin{itemize}
\item [\rm{(i)}] If $N_1(\la)=R(\la)N_2(\la)$ with $R(\la) \in\FF_{pr}(\la)^{n\times m}$ then
$\deg(\col_j(N_1(\la))) \leq \deg(\col_j(N_2(\la)))$ for $j=1, \ldots, l$.
\item [\rm{(ii)}] If $N_1(\la)=R(\la)N_2(\la)$ with $R(\la)$  strictly proper then $\deg(\col_j(N_1(\la))) < \deg(\col_j(N_2(\la)))$ for $j=1, \ldots, l$.
\end{itemize}
\end{lemma}

\textbf{Proof}.-
Let $\frac{p(\la)}{q(\la)}$ be  a proper rational function and let $n(\la)$ be a
polynomial. Notice that
%From (\ref{valpro}) we have
\begin{equation}\label{val}
\deg(p(\la))+\deg(n(\la))-\deg(q(\la))\leq\deg(n(\la)).
%\delta_\infty(r(\la)n(\la))=\delta_\infty(r(\la))+\delta_\infty(n(\la))\geq \delta_\infty(n(\la))=-\deg(n(\la)).
\end{equation}
Let $n_{ij}^{(1)}(\la)$ be an arbitrary element of the $j$-th column of $N_1(\la)$. As $N_1(\la)=R(\la)N_2(\la)$, we can write
$n_{ij}^{(1)}(\la)=\sum_{k=1}^m \frac{p_{ik}(\la)}{q_{ik}(\la)}n_{kj}^{(2)}(\la)$ where $\frac{p_{ik}(\la)}{q_{ik}(\la)}$ is the element in position $(i,k)$ of
$R(\la)$ and $n_{kj}^{(2)}(\la)$ is the element in position $(k,j)$ of $N_2(\la)$. It follows from
%(\ref{valsum}) and
(\ref{val}) that
for each element of the $j$-th column of $N_1(\la)$
$$
\begin{array}{l}\deg\left(n_{ij}^{(1)}(\la)\right)=\deg\left(\sum_{k=1}^m\frac{p_{ik}(\la)}{q_{ik}(\la)}n_{kj}^{(2)}(\la)\right)
\\ \leq
\max_{k}\{\deg(p_{ik}(\la))+\deg(n_{kj}^{(2)}(\la))-\deg(q_{ik}(\la))\}\\
\leq \max_k\{\deg(n_{kj}^{(2)}(\la))\}=\deg(\col_j(N_2(\la)),\qquad\text{for all } i=1,\ldots,n.
\end{array}
$$
%$$\delta_\infty(n_{ij}^1(\la))
%=\delta_\infty(\sum_{k=1}^m r_{ik}(\la)n_{kj}^2(\la))
%\geq \min \{\delta_\infty(r_{ik}(\la)n_{kj}^2(\la))\} \geq \min \{- \deg(n_{kj}^2(\la))\}
%=-\max \{\deg(n_{kj}^2(\la))\}
%.$$
%We conclude that
%$$- \deg(n_{ij}^{(1)}(\la))\geq- \max \{\deg(n_{kj}^2(\la))\}=
%-\deg (n_j^2(\la))
%- \deg(col_j(N_2(\la)))
%\text{ for all } i.$$
Then
$\deg(\col_j(N_1(\la)))= \max_i \{\deg(n_{ij}^{(1)}(\la))\} \leq \deg(\col_j(N_2(\la)))$ and (i) follows.  If  $R(\la)$
is strictly proper the previous inequality is strict.
\hfill\halmos

\bigskip
The following corollary is an immediate consequence of Lemma \ref{coldeg}.
\begin{corollary}\label{degcol}
With the same assumptions and notation of Lemma \ref{basisPG},
%\begin{itemize}
%\item [(i)]
if $A(\la)^{-1} B(\la)$ is proper then $\deg(\col_j(H_1(\la))) \leq \deg(\col_j(H_2(\la)))$ for all $j$. The inequality is strict if
$A(\la)^{-1} B(\la)$ is strictly proper.
%\item [(ii)] If $A(\la)=\la E-A$ with $E$ nonsingular and
%$B(\la)=B$ then $\deg(col_j(N_1(\la))) < \deg(col_j(N_2(\la)))$ for $j=1, \ldots, l$.
%\end{itemize}
\end{corollary}

%\textbf{Proof}.- It follows immediately on applying Lemma
%\ref{coldeg} to the matrix $R(\la)=A(\la)^{-1}   B(\la)$. %and  $R(\la)=(\la
%E-A)^{-1}B$.

%\hfill\halmos

%This corollary allows us to give the following theorem.
%  that relates minimal bases of a polynomial system matrix with  minimal bases of its transfer function rational matrix.
The following lemma  relates the minimal bases of a rational matrix and its transpose as well as their minimal indices. It also states that the transpose of a polynomial system matrix gives rise to the transpose of its transfer function. It can be proved straightforwardly and, therefore, the proof is omitted.
\begin{lemma}\label{lem_trans}
\begin{itemize}
\item[\rm{(a)}] For any rational matrix $G(\la)$, $\mathcal{N}_\ell(G(\la))=\mathcal{N}_r(G(\la)^T)$ and $\mathcal{N}_r(G(\la))=\mathcal{N}_\ell(G(\la)^T)$. Moreover, $H(\la)$ is a left minimal basis of $G(\la)$ if and only if it is a right minimal basis of $G(\la)^T$. Also the left minimal indices of $G(\la)$ and the right minimal indices of $G(\la)^T$ coincide.
\item[\rm{(b)}] If $P(\la)$ is a (minimal) polynomial system matrix giving rise to $G(\la)$ then $P(\la)^T$ is a (minimal) polynomial system matrix giving rise to $G(\la)^T$.
\end{itemize}
\end{lemma}

{As announced, the following result shows how to obtain a minimal basis of a rational matrix from a minimal basis of those of its polynomial system matrices that satisfy coprimeness and properness conditions, and relates their minimal indices.}

\begin{theorem}\label{minimalbasisPG}
%   Let $G(\la)\in\FF(\la)^{p\times m}$ be a rational matrix  and
 Let $P(\la)$ of \eqref{s1.eqpolsysmat}  be a polynomial system matrix
of a rational matrix $G(\la)$.
\begin{itemize}
\item [\rm{(a)}] If $A(\la)$ and $C(\la)$ are right coprime,  $A(\la)^{-1}B(\la)$ is proper and $
    \begin{bsmallmatrix}
    H_1(\la) \\H_2(\la)
    \end{bsmallmatrix}
    \in\FF[\la]^{(n+m)\times l}
    $ is a right minimal basis of $P(\la)$, then $H_2(\la)$ is a   right minimal basis of $G(\la)$ and $H_1(\la)=-A(\la)^{-1}    B(\la)  H_2(\la)$.
    Moreover,
the right minimal indices of $P(\la)$ and $G(\la)$ are the same.

\item [\rm{(b)}] If $A(\la)$ and $B(\la)$ are left coprime,  $C(\la)A(\la)^{-1}$ is proper and $
    \begin{bsmallmatrix}
    H_1(\la) \\H_2(\la)
    \end{bsmallmatrix}
    \in\FF[\la]^{(n+p)\times q}
    $ is a left minimal basis of $P(\la)$, then $H_2(\la)$ is a left minimal basis of $G(\la)$ and $H_1(\la)=(C(\la)A(\la)^{-1})^T  H_2(\la)$.
    Moreover,
the left minimal indices of $P(\la)$ and $G(\la)$ are the same.
\end{itemize}
% such that  $A(\la)$, $C(\la)$ are right coprime and  $A(\la)^{-1}   B(\la)$ is strictly proper.
%    If  $\label{}
%    \begin{bmatrix}
%    H_1(\la) \\H_2(\la)
%    \end{bmatrix}
%    \in\FF[\la]^{(n+m)\times l}
%    $ is a right minimal basis of $P(\la)$, then $H_2(\la)$ is a   right minimal basis of $G(\la)$ and $H_1(\la)=-A(\la)^{-1}    B(\la)  H_2(\la)$.
%    Moreover,
%the right minimal indices of $P(\la)$ and $G(\la)$ are the same.
\end{theorem}

\textbf{Proof}.-
We prove part (a). By Lemma \ref{basisPG}, $H_2(\la)$ is a right  polynomial basis of $G(\la)$ %(i.e., of $\mathcal{N}_r (G(\la))$)
and
%\begin{equation}\label{eq}
$H_1(\la)=-A(\la)^{-1}   B(\la)  H_2(\la)$.
%\end{equation}
We show that $H_2(\la)$ is a minimal basis of $G(\la)$ by applying Theorem \ref{mincharact}. Let us prove first that $H_2(\la_0)$ has full column rank for all
 $\la_0\in\overline{\FF}$. If this were not true, there would exist  $\la_1 \in \overline{\FF}$ and a vector, $v \neq 0$, such that $H_2(\la_1)v=0$. But since
$\begin{bsmallmatrix}
    H_1(\la) \\H_2(\la)
\end{bsmallmatrix}
$
is a right minimal basis for $P(\la)$,
$\begin{bsmallmatrix}
H_1(\la_1) \\H_2(\la_1)
\end{bsmallmatrix}v=\begin{bsmallmatrix}
w\\0
\end{bsmallmatrix}
$
with $w\neq 0$ and
$$P(\la_1)\begin{bmatrix}
H_1(\la_1) \\H_2(\la_1)
\end{bmatrix}v=
\begin{bmatrix}
A(\la_1) & B(\la_1)\\-C(\la_1) &D(\la_1)
\end{bmatrix}
\begin{bmatrix}
w\\0
\end{bmatrix}=\begin{bmatrix}
A(\la_1)\\- C(\la_1)
\end{bmatrix}w=0.
$$
This would be a contradiction because $A(\la)$ and $C(\la)$ are right coprime, i.e.,
$\begin{bsmallmatrix}
A(\la_1)\\- C(\la_1)
\end{bsmallmatrix}$ has full column rank (see Proposition \ref{prop_coprim}).

Next, let us see that $H_2(\la)$ is  column reduced. By hypothesis and Theorem \ref{mincharact}, we know that
$\begin{bsmallmatrix}
    H_1(\la) \\H_2(\la)
\end{bsmallmatrix}$
is column reduced. Our goal is to express the highest column degree coefficient matrix of $\begin{bsmallmatrix}
    H_1(\la) \\H_2(\la)
\end{bsmallmatrix}$
in terms of the highest column degree coefficient matrix of $H_2 (\la)$, which is denoted by $H_{2h}$. For this purpose, note that the assumption that $A(\la)^{-1} B(\la)$ is proper implies that
$-A(\la)^{-1} B(\la) = J + R(\la)$, where $J$ is a constant matrix and $R(\la)$ is strictly proper. Thus, $H_1 (\la) = J H_2(\la) + R(\la) H_2(\la)$ and
\begin{equation} \label{eq:aux1}
\col_j \left( \begin{bmatrix}
    H_1(\la) \\H_2(\la)
\end{bmatrix} \right) =
\begin{bmatrix}
    J \col_j (H_2(\la)) + R(\la) \col_j(H_2(\la)) \\ \col_j (H_2(\la))
\end{bmatrix}.
\end{equation}
Bear in mind that $\col_j (H_2(\la)) \ne 0$ since $H_2(\la)$ is a right polynomial basis of $G(\la)$. Moreover, $R(\la) \col_j(H_2(\la))$ is a vector polynomial, because $\col_j(H_1(\la))$ and $J \col_j (H_2(\la))$ are both vector polynomials. Then, Lemma \ref{coldeg} (ii) guarantees that $\deg (R(\la) \col_j(H_2(\la))) <  \deg (\col_j(H_2(\la)))$. Therefore, the highest degree coefficient of \eqref{eq:aux1} is
$\begin{bsmallmatrix}
    J \col_j (H_{2h}) \\ \col_j (H_{2h})
\end{bsmallmatrix}$,
the degree $d_j$ of \eqref{eq:aux1} is $d_j = \deg (\col_j (H_2(\la)))$ and the highest column degree coefficient matrix of $\begin{bsmallmatrix}
    H_1(\la) \\H_2(\la)
\end{bsmallmatrix}$ is $\begin{bsmallmatrix}
    J H_{2h} \\H_{2h}
\end{bsmallmatrix}$. This latter matrix has full column rank, which implies that $H_{2h}$ has also full column rank, since otherwise there would exist a nonzero constant vector $v$ such that $H_{2h} v = 0$ and $\begin{bsmallmatrix}
    J H_{2h} \\H_{2h}
\end{bsmallmatrix} v =0$, which is a contradiction. This proves that $H_2 (\la)$ is column reduced and, so, a right minimal basis of $G(\la)$. Since the degrees of the corresponding columns of $H_2 (\la)$ and $\begin{bsmallmatrix}
    H_1(\la) \\H_2(\la)
\end{bsmallmatrix}$ coincide, the right minimal indices of $G(\la)$ and $P(\la)$ are the same.

Part (b) is a consequence of part (a) and Lemma \ref{lem_trans}.
\hfill\halmos
%To prove part (b), as $\begin{bmatrix}
%    J_1(\la) \\J_2(\la)
%    \end{bmatrix}$ is a left minimal basis of $P(\la)$ we have that
%$
%\begin{bmatrix}
%    J_1(\la)^T&J_2(\la)^T
%    \end{bmatrix}P(\la)=0
%$, $J_1(\la)^T=J_2(\la)^TC(\la)A(\la)^{-1}$ and, therefore, $J_1(\la)=(C(\la)A(\la)^{-1})^TJ_2(\la)$.
%\hfill\halmos

%\begin{rem}
%   {\rm In the previos theorem, we only use that the matrices  $A(\la)$ and %$C(\la)$ are right coprime. Therefore, we can replace the minimality condition %of the polynomial system matrix by a more weak condition: the right %coprimeness condition between thw matrices
%   $A(\la)$ and $C(\la)$.      }
%\end{rem}

%\begin{rem}
%    {\rm
%        If the polynomial system matrix in Theorem \ref{minimalbasisPG} is minimal (i.e.,  $A(\la)$, $B(\la)$ are left coprime and %$A(\la)$, $C(\la)$ are right
%        coprime)    the previous result can be applied if $A(\la)^{-1}   B(\la)$ is strictly proper.
%    }
%\end{rem}

\bigskip
We have shown so far how to obtain a minimal basis of the transfer function matrix of a polynomial system matrix out of a minimal basis of the latter,  which is the most interesting scenario in applications. For completeness, we consider now
the reciprocal problem. In this respect, Lemma \ref{basisPG} motivates the following result.

\begin{lemma}\label{basisGP}
        Let
        $
        P(\la)
        $
        of \eqref{s1.eqpolsysmat}
        be a polynomial system matrix
%   whose transfer  function rational matrix is
        of a rational matrix $G(\la)$ where $A(\la)$ and  $C(\la)$ are right coprime. Let $H_2(\la)$ be  a
        right polynomial basis of $G(\la)$ and let $H_1(\la)=-A(\la)^{-1}    B(\la)  H_2(\la)$. Then
        $
        \begin{bsmallmatrix}
        H_1(\la) \\H_2(\la)
        \end{bsmallmatrix}
        %\in\FF[\la]^{(n+m)\times (l)}
        $ is a  right polynomial basis of $P(\la)$.
\end{lemma}

\textbf{Proof}.-
Note that
$$P(\la)\begin{bmatrix}
H_1(\la) \\H_2(\la)
\end{bmatrix}=\begin{bmatrix}
A(\la) & B(\la)\\-C(\la) &D(\la)
\end{bmatrix}\begin{bmatrix}
-A(\la)^{-1}    B(\la)  H_2(\la) \\H_2(\la)
\end{bmatrix}=\begin{bmatrix}
0\\G(\la)H_2(\la)
\end{bmatrix}=0.$$
%It only remains to prove that $H_1(\la)$ is polynomial.
Let us see first that $H_1(\la)$ is polynomial. As $A(\la)$ and $C(\la)$ are right coprime, by Bezout's
identity (see Proposition \ref{prop_coprim}), there exist polynomial matrices $X(\la)$ and $Y(\la)$ of
appropriate sizes such that
% $X(\la)A(\la)+ Y(\la)C(\la)=I_n$,  or
 $$\begin{bmatrix}
X(\la) & -Y(\la)
\end{bmatrix}\begin{bmatrix}
    A(\la) \\ -C(\la)
\end{bmatrix}=I_n.$$
Put $H(\la)=X(\la)B(\la)- Y(\la)D(\la)$. Then,
$$\begin{bmatrix}
X(\la) & -Y(\la)
\end{bmatrix}\begin{bmatrix}
A(\la) & B(\la)\\-C(\la) &D(\la)
\end{bmatrix}\begin{bmatrix}
H_1(\la) \\H_2(\la)
\end{bmatrix}=\begin{bmatrix}
I_n & H(\la)
\end{bmatrix}\begin{bmatrix}
H_1(\la) \\H_2(\la)
\end{bmatrix}=0.$$
Hence $H_1(\la)=-H(\la)H_2(\la)$ is a matrix polynomial.
Moreover,
$\begin{bsmallmatrix}
H_1(\la) \\H_2(\la)
\end{bsmallmatrix}$ is a right polynomial basis of $P(\la)$, because its columns belong to $\mathcal{N}_r (P(\la))$, its columns are linearly independent, since $H_2(\la)$ is a basis of $\mathcal{N}_r (G(\la))$, and  $\dim \mathcal{N}_r (G(\la)) \allowbreak =\dim \mathcal{N}_r (P(\la)).$
\hfill\halmos

\bigskip
We can prove now the reciprocal of Theorem \ref{minimalbasisPG}, which shows that, under certain assumptions, minimal bases of polynomial system matrices can be obtained from minimal bases of their transfer functions.
\begin{theorem}\label{minimalbasisPGrec}
Let $G(\la)\in\FF(\la)^{p\times m}$ be a rational matrix  and
let
    $
    P(\la)$ of \eqref{s1.eqpolsysmat}
    be a polynomial system matrix of $G(\la)$.
\begin{itemize}
\item [\rm{(a)}] If $A(\la)$ and $C(\la)$ are right coprime, $A(\la)^{-1}B(\la)$ is proper, $H_2(\la)$ is a  right minimal basis of
$G(\la)$ and $H_1(\la)=-A(\la)^{-1} B(\la)H_2(\la)$ then
$
    \begin{bsmallmatrix}
    H_1(\la) \\H_2(\la)
    \end{bsmallmatrix}
    %\in\FF[\la]^{(n+m)\times (l)}
    $ is a right minimal basis of $P(\la)$.
 Moreover, the right minimal indices of $P(\la)$ and $G(\la)$ are the same.

\item [\rm{(b)}] If $A(\la)$ and $B(\la)$ are left coprime, $C(\la)A(\la)^{-1}$ is proper, $H_2(\la)$ is a left minimal basis of
$G(\la)$ and $H_1(\la)=(C(\la)A(\la)^{-1})^T H_2(\la)$ then
$
    \begin{bsmallmatrix}
    H_1(\la) \\H_2(\la)
    \end{bsmallmatrix}
    $ is a left minimal basis of $P(\la)$.
 Moreover, the left minimal indices of $P(\la)$ and $G(\la)$ are the same.
\end{itemize}
%such that $A(\la)$, $C(\la)$ are right coprime and  $A(\la)^{-1}B(\la)$ is proper. If $H_2(\la)$ is a   right minimal basis of
%$G(\la)$ and $H_1(\la)=-A(\la)^{-1} B(\la)  H_2(\la)$ then
%$\label{}
%    \begin{bmatrix}
%    H_1(\la) \\H_2(\la)
%    \end{bmatrix}
    %\in\FF[\la]^{(n+m)\times (l)}
%    $ is a right minimal basis of $P(\la)$.
% Moreover, the right minimal indices of $P(\la)$ and $G(\la)$ are the same.
\end{theorem}

\textbf{Proof}.- We prove part (a). By Lemma \ref{basisGP} and  Theorem \ref{mincharact}, we just need to prove that   $\begin{bsmallmatrix}
    H_1(\la_0) \\H_2(\la_0)
\end{bsmallmatrix}$ has full column rank for all $\la_0\in\overline{\FF}$ and $\begin{bsmallmatrix}
H_1(\la) \\H_2(\la)
\end{bsmallmatrix}$ is column reduced.
As  $H_2(\la)$ is a  right minimal basis of $G(\la)$,  $H_2(\la_0)$ has  full column rank  for all $\la_0\in\overline{\FF}$, which implies that the matrix $\begin{bsmallmatrix}
H_1(\la_0) \\H_2(\la_0)
\end{bsmallmatrix}$ has full column rank as well. Moreover, $H_2(\la)$ is column reduced. Write
$H_2(\la)=H_{2h}\diag(\la^{d_1}, \ldots, \la^{d_l})+L_2(\la)$
with $H_{2h}$ of full column rank, $d_1,\ldots, d_l$ the right minimal indices of $G(\la)$ and the degree of the $j$-th column of $L_2(\la)$ less than $d_j$ for each $j$. Since $H_1(\la)= -A(\la)^{-1} B(\la)  H_2(\la)$, with  $A(\la)^{-1}    B(\la)$ proper, it follows from Corollary \ref{degcol}
that each column of $H_1(\la)$ has degree less than or equal to the same column of $H_2(\la)$. Therefore,  there is a matrix $H_{1h}$ such that the
highest column degree coefficient matrix of $\begin{bsmallmatrix}
    H_1(\la) \\H_2(\la)
\end{bsmallmatrix}$ is $\begin{bsmallmatrix}
 H_{1h} \\H_{2h}
\end{bsmallmatrix}$, a
full column rank matrix. Moreover, its column degrees are those of $H_2(\la)$.

Part (b) follows from (a) and Lemma \ref{lem_trans}.
\hfill\halmos

\bigskip
Theorems \ref{minimalbasisPG} and \ref{minimalbasisPGrec} together provide our next result.
\begin{corollary}\label{col_PGPmin}
Let $G(\la)\in\FF(\la)^{p\times m}$ be a rational matrix  and
let
    $
    P(\la)
    $ of \eqref{s1.eqpolsysmat}
    be a minimal polynomial system matrix of $G(\la)$.
If both $A(\la)^{-1}B(\la)$ and $C(\la)A(\la)^{-1}$ are proper matrices then
$\begin{bsmallmatrix}H_1(\la) \\H_2(\la)\end{bsmallmatrix}$ is a right (resp., left) minimal basis of $P(\la)$ if and only if
$H_2(\la)$ is a  right (resp., left) minimal basis of
$G(\la)$ and $H_1(\la)=-A(\la)^{-1} B(\la)H_2(\la)$ (resp., $H_1(\la)=(C(\la)A(\la)^{-1})^T H_2(\la)$). Moreover, the right (resp., left) minimal indices of $P(\la)$ and $G(\la)$ are the same.
\end{corollary}
%\begin{rem}
%	{\rm
%		If the polynomial system matrix in Theorem \ref{minimalbasisPGrec} is minimal (i.e.,  $A(\la)$, $B(\la)$ are left coprime and $A(\la)$, $C(\la)$ are right
%		coprime)    the previous result can be applied if $A(\la)^{-1}   B(\la)$ is strictly proper.
%	}
%\end{rem}

{
\begin{rem}\rm
As announced at the beginning of this section, we show that the conditions on properness in Theorems 3.6 and 3.8 and Verghese's conditions on infinite zeros in \cite{Ver80} are not equivalent. More precisely, we prove by means of two examples that for a polynomial system matrix $P(\la)$  as in (\ref{s1.eqpolsysmat}), $C(\la)A(\la)^{-1}$ being proper does not imply $\begin{bsmallmatrix} A(\la)&B(\la)&0\\-C(\la)&D(\la)&-I\end{bsmallmatrix}$ having no infinite zeros an vice versa, even when $A(\la)$ and $B(\la)$ are left coprime. The same two examples prove that $A(\la)^{-1} B(\la)$ being proper does not imply $\begin{bsmallmatrix} A(\la)&B(\la)\\-C(\la)&D(\la)\\ 0 & I\end{bsmallmatrix}$ having no infinite zeros an vice versa, even when $A(\la)$ and $C(\la)$ are right coprime.
\begin{example}\rm
Let $A(\la)=\begin{bsmallmatrix}\la+1&\la^2\\1&\la\end{bsmallmatrix}, B(\la)=\begin{bsmallmatrix}1\\0\end{bsmallmatrix}, C(\la)=\begin{bsmallmatrix}0&1\end{bsmallmatrix}, D(\la)=0$. It is easy to check that $G(\la)=-\frac{1}{\la}$, $A(\la)$ and $B(\la)$ are left coprime and $C(\la)A(\la)^{-1}=\begin{bsmallmatrix}-\frac{1}{\la},&1+\frac{1}{\la}\end{bsmallmatrix}$ is proper. However, $\begin{bsmallmatrix} A(\la)&B(\la)&0\\-C(\la)&D(\la)&-1\end{bsmallmatrix}$
%$\begin{bsmallmatrix}\la+1&\la^2&-1&0\\1&\la&0&0\\0&1&0&-1\end{bsmallmatrix}
%$
has an infinite zero, because its invariant orders at infinity are $q_1=-2$, $q_2=0$ and $q_3=1$. Moreover, $A(\la)$ and $C(\la)$ are right coprime, $A(\la)^{-1} B(\la)$ is proper, but $\begin{bsmallmatrix} A(\la)&B(\la)\\ -C(\la)&D(\la)\\0&1\end{bsmallmatrix}$
has again an infinite zero, because its invariant orders at infinity are also $q_1=-2$, $q_2=0$ and $q_3=1$. Thus, $P(\la)$ of \eqref{s1.eqpolsysmat} is not strongly irreducible, but allows to recover the complete finite and infinite eigenstructures of $G(\la)$ as well as its minimal indices and bases, since it satisfies the coprimeness and properness conditions.
%There are several ways to prove this. By using the fact that $-(q_1+\cdots+q_k)=\max\{\deg(\text{minor of order }k)\}$, it turns out that $-q_1=2$, $-(q_1+q_2)=2$ and $-(q_1+q_2+q_3)=1$ and, therefore, $q_1=-2$, $q_2=0$ and $q_3=1$.
\end{example}

\begin{example}\rm
Let $A(\la)=\begin{bsmallmatrix} \la & 0 & 0 \\ 0&1&0\\0&1&1\end{bsmallmatrix}, B(\la)=\begin{bsmallmatrix}1\\ \la\\1\end{bsmallmatrix}, C(\la)=\begin{bsmallmatrix} 1& 0&\la\end{bsmallmatrix}, D(\la)=0$. It is easy to check that $G(\la)=\frac{1}{\la}-\la^2+\la$, and $A(\la)$ and $B(\la)$ are left coprime. The matrix
%$
%\begin{bsmallmatrix}1&0&-\la&0\\1&1&-1&0\\0&\la&0&-1\end{bsmallmatrix}
%$
$\begin{bsmallmatrix} A(\la)&B(\la)&0\\-C(\la)&D(\la)&-1\end{bsmallmatrix}$
has no infinite zeros since %$-q_1=1$, $-(q_1+q_2)=2$ and $-(q_1+q_2+q_3)=2$, that is, $q_1=-1$, $q_2=-1$ and $q_3=0$.
its invariant orders at infinity are $q_1=-1$, $q_2=-1$, $q_3=-1$ and $q_4 = 0$.
Nevertheless, $C(\la)A(\la)^{-1}=\begin{bsmallmatrix} \frac{1}{\la},&-\la ,&\la\end{bsmallmatrix}$ is not proper. Moreover, $A(\la)$ and $C(\la)$ are right coprime, but the matrix $\begin{bsmallmatrix} A(\la)&B(\la)\\ -C(\la)&D(\la)\\0&1\end{bsmallmatrix}$ has no infinite zeros, since its invariant orders at infinity are again $q_1=-1$, $q_2=-1$, $q_3=-1$ and $q_4 = 0$. However, the matrix $A(\la)^{-1} B(\la) = \begin{bsmallmatrix} \frac{1}{\la},&\la ,&-\la+1\end{bsmallmatrix}$ is not proper. Thus, $P(\la)$ of \eqref{s1.eqpolsysmat} is strongly irreducible and, so, allows to recover the complete finite and infinite eigenstructures of $G(\la)$ as well as its minimal indices and bases, but none of the properness conditions hold.
\end{example}
\end{rem}
}

%---------------------------------------------------------------------------------------------
\section{Polynomial bases of linearizations of rational matrices}\label{Sec_linearizations}
%-------------------------------------------------------------------------------------------

The aim of this section is to study the relationship between  the
polynomial  bases of a rational matrix and the polynomial bases
of its linearizations. It is not possible to extend this relationship to minimal bases because it was already proved
in \cite[Theorem 4.10 (b)]{DeDoMa14} that the minimal bases and indices of a polynomial matrix can not be obtained from the minimal bases and indices of its linearizations in general, and polynomial matrices are particular cases of rational matrices.

A linear pencil
\be\label{eq.deflinar}
L(\la)=\begin{bmatrix}
A_1 \la +A_0 & B_1 \la +B_0\\-(C_1 \la +C_0) & D_1 \la +D_0
\end{bmatrix}
\ee
is said to be a
linearization of a rational matrix $G(\la)$ (see \cite[Definition  3.2]{AmDoMaZa18}) if it is a  minimal polynomial system matrix of a rational
matrix $\wh{G}(\la)$ such that, for some nonnegative integers $s_1,s_2$, $\diag(\wh{G}(\la), I_{s_2})$ and
$\diag(G(\la), I_{s_1})$ are unimodularly equivalent. We can assume without loss of generality that $s_1=s$ and $s_2=0$.
This assumption will be adopted in the rest of the paper every time we deal with linearizations.

A first consequence of this definition is that, by the rank-nullity theorem, $\dim \mathcal{N}_r(\wh{G}(\la))=\dim \mathcal{N}_r(G(\la))$ and $\dim \mathcal{N}_\ell(\wh{G}(\la))=\dim \mathcal{N}_\ell(G(\la))$. Therefore, $G(\la)$ and $\wh{G}(\la)$ have the same number of right minimal indices and the same number of left minimal indices. Furthermore, by Lemma \ref{rankPG}, $\dim \mathcal{N}_r(\wh{G}(\la))=\dim \mathcal{N}_r(L(\la))$ and $\dim \mathcal{N}_\ell(\wh{G}(\la))=\dim \mathcal{N}_\ell(L(\la))$. Thus, a rational matrix and any of its linearizations have the same number of right minimal indices and the same number of left minimal indices.

Proposition \ref{prop_UwhGV} relates right polynomial bases of $G(\la)$ and $\wh{G} (\la)$. An analogous result holds for left polynomial bases of $G(\la)$ and $\wh{G} (\la)$ as a consequence of Lemma \ref{lem_trans}. Such ``left'' result is omitted for brevity.

\begin{proposition}\label{prop_UwhGV}
 Let $G(\la) \in \FF(\la)^{p\times m}$ and let $\wh{G}(\la) \in \FF(\la)^{(p+s)\times (m+s)}$, $s\geq 0$. Let $U(\la)\in\efe[\la]^{(p+s)\times (p+s)}$ and $V(\la)\in\efe[\la]^{(m+s)\times (m+s)}$ be unimodular matrices such that
 $U(\la)\wh{G}(\la)V(\la)=\diag(G(\la),I_s).$
% Let $l=m-\rank G(\la)$.
 \begin{enumerate}
 \item[\rm{(a)}] If $H(\la)$ is a right polynomial basis of $G(\la)$ then $V(\la)\begin{bsmallmatrix}H(\la) \\ 0 \end{bsmallmatrix}$ is a right polynomial basis of $\wh{G}(\la)$.

 \item[\rm{(b)}] If $\wh{H}(\la)$ is a right polynomial basis of $\wh{G}(\la)$ then $V(\la)^{-1}\wh{H}(\la)=\begin{bsmallmatrix}H(\la) \\ 0 \end{bsmallmatrix}$ and $H(\la)$ is a right polynomial basis of $G(\la)$.
 \end{enumerate}
\end{proposition}

\textbf{Proof}.-
In order to prove (a) assume that $G(\la)H(\la)=0$. We obtain, via a direct multiplication, that
 	$$%U(\la)\wh{G}(\la)V(\la)\begin{bmatrix} H(\la) & 0\\ 0& I_s \end{bmatrix} = \begin{bmatrix} G(\la) & 0\\ 0& I_s \end{bmatrix}\begin{bmatrix} H(\la) & 0\\ 0& I_s \end{bmatrix} =\begin{bmatrix} 0& 0\\ 0& I_s \end{bmatrix}.
 	\wh{G}(\la)V(\la)\begin{bmatrix}H(\la) \\ 0 \end{bmatrix}=U(\la)^{-1}\begin{bmatrix}G(\la)&0\\0&I_s\end{bmatrix}
 	\begin{bmatrix}H(\la) \\ 0 \end{bmatrix}=0.
 	$$
 	%Thus, $\wh{G}(\la)V(\la)\begin{bmatrix}H(\la) \\ 0 \end{bmatrix}=0$
 	So, $V(\la)\begin{bsmallmatrix}H(\la) \\ 0 \end{bsmallmatrix}$ is a right polynomial basis of $\wh{G}(\la)$, because its columns are linearly independent and $\dim \mathcal{N}_r (\wh{G}(\la))=
 	%m+s- \rank \wh{G}(\la)= m+s-(s+\rank G(\la))=m-\rank G(\la)=
 	\dim \mathcal{N}_r (G(\la)).$
 	
For proving (b) assume that  $\wh{G}(\la) \wh{H}(\la)=0$. Therefore,
%$U(\la)^{-1}\begin{bmatrix} G(\la) & 0\\ 0& I_s \end{bmatrix}V(\la)^{-1}\wh{H}(\la)=0$ or
%\begin{equation}\label{whbasis}
$$
\begin{bmatrix} G(\la) & 0\\ 0& I_s \end{bmatrix}V(\la)^{-1}\wh{H}(\la)=0.
$$
%\end{equation}
Write $V(\la)^{-1}=\begin{bsmallmatrix}V_1(\la) \\ V_2(\la) \end{bsmallmatrix}$, where
$V_1(\la)\in \efe [\la]^{ m \times (m+s)}$ and $V_2(\la)\in \efe [\la]^{ s \times (m+s)}$. Thus, $G(\la)V_1(\la)\wh{H}(\la)=0$ and $V_2(\la)\wh{H}(\la)=0$. Set $H(\la)=V_1(\la)\wh{H}(\la)$.  It follows that $V(\la)^{-1}\wh{H}(\la)=\begin{bsmallmatrix}H(\la) \\ 0 \end{bsmallmatrix}$ and $G(\la)H(\la)=0$. Thus, the columns of $H(\la)$ form a right polynomial basis of $G(\la)$.
\hfill\halmos

\begin{rem}
\rm{Proposition \ref{prop_UwhGV} cannot be extented to right minimal bases, i.e., if $H(\la)$ is a right minimal basis of $G(\la)$,  $V(\la)\begin{bsmallmatrix}H(\la) \\ 0 \end{bsmallmatrix}$ may not be a right minimal basis of $\wh{G}(\la)$, and if $\wh{H} (\la)$ is a right minimal basis of $\wh{G}(\la)$, $V(\la)^{-1} \wh{H}(\la)$ may not contain in its first $m$ rows a minimal basis of $G(\la)$. {\color{black} Otherwise, if Proposition \ref{prop_UwhGV}  could be extended to right minimal bases, taking $G(\la)$ polynomial and $\wh{G}(\la)$ a linearization of $G(\la)$, then their right minimal bases and indices would be always related, which is in contradiction with \cite[Theorem 4.10 (b)]{DeDoMa14}.
}
}
\end{rem}

The next result relates the polynomial bases of a rational matrix and its linearizations through the unimodular matrices that connect the rational matrix and the transfer function matrix of the linearizations.

\begin{theorem}
Let $G(\la)\in\efe(\la)^{p\times m}$ and let
$L(\la)$ of \eqref{eq.deflinar}
be a linearization of $G(\la)$ with transfer function matrix $\wh{G}(\la)$. Let $U(\la)\in\efe[\la]^{(p+s)\times (p+s)}$, $V(\la)\in\efe[\la]^{(m+s)\times (m+s)}$ be unimodular matrices such that
 $U(\la)\wh{G}(\la)V(\la)=\diag(G(\la),I_s).$
 \begin{enumerate}
 \item[\rm{(a)}] $\begin{bsmallmatrix}H_1(\la)\\H_2(\la)\end{bsmallmatrix}$ is a right polynomial basis of $L(\la)$ if and only if $H_2(\la)=V(\la)\begin{bsmallmatrix}H(\la) \\ 0 \end{bsmallmatrix}$ for some right polynomial basis $H(\la)$ of $G(\la)$ and $H_1(\la)=-(A_1 \la +A_0)^{-1}(B_1 \la +B_0)H_2(\la)$.

 \item[\rm{(b)}] $\begin{bsmallmatrix}H_1(\la)\\H_2(\la)\end{bsmallmatrix}$ is a left polynomial basis of $L(\la)$ if and only if $H_2(\la)=U(\la)^{T}\begin{bsmallmatrix}H(\la) \\ 0 \end{bsmallmatrix}$ for some left polynomial basis $H(\la)$ of $G(\la)$ and $H_1(\la)=((C_1\la+C_0)(A_1 \la +A_0)^{-1})^T H_2(\la)$.
% \item[\rm{(a)}] If $H(\la)$ is a right polynomial basis of $G(\la)$ and $\wh{H}(\la)=V(\la)\begin{bmatrix}H(\la) \\ 0 %\end{bmatrix}$ then
% $\begin{bmatrix}(A_1 \la +A_0)^{-1}(B_1 \la +B_0)\wh{H}(\la) \\ \wh{H}(\la) \end{bmatrix}$ is a right polynomial basis of $L(\la)$.

% \item[\rm{(b)}] If $J(\la)$ is a left polynomial basis of $G(\la)$ and $\wh{J}(\la)=U(\la)^T\begin{bmatrix}J(\la) \\ 0 \end{bmatrix}$ then
%  $\begin{bmatrix}((C_1\la+C_0)(A_1 \la +A_0)^{-1})^T\wh{J}(\la) \\ \wh{J}(\la) \end{bmatrix}$ is a left polynomial basis of $L(\la)$.

% \item[\rm{(c)}] If $\begin{bmatrix}\wh{H}_1(\la)\\\wh{H}_2(\la)\end{bmatrix}$ is a right polynomial basis of $L(\la)$ then $\wh{H}_1(\la)=(A_1 \la +A_0)^{-1}(B_1 \la +B_0)\wh{H}_2(\la)$, $V(\la)^{-1}\wh{H}_2(\la)=\begin{bmatrix}H(\la) \\ 0 \end{bmatrix}$ and $H(\la)$ is a right polynomial basis of $G(\la)$.

% \item[\rm{(d)}] If $\begin{bmatrix}\wh{J}_1(\la)\\\wh{J}_2(\la)\end{bmatrix}$ is a left polynomial basis of $L(\la)$ then $\wh{J}_1(\la)=((C_1\la+C_0)(A_1 \la +A_0)^{-1})^T\wh{J}_2(\la)$, $U(\la)^{-T}\wh{J}_2(\la)=\begin{bmatrix}J(\la) \\ 0 \end{bmatrix}$ and $J(\la)$ is a left polynomial basis of $G(\la)$.
  \end{enumerate}
\end{theorem}
\textbf{Proof}.- As $L(\la)$ is a linearization of $G(\la)$, $L(\la)$ is a minimal polynomial system matrix and, therefore,  $A_1\la+A_0$ and $C_1\la+C_0$ are right coprime and $A_1\la+A_0$ and $B_1\la+B_0$ are left coprime. Thus we can apply Lemmas \ref{basisPG} and \ref{basisGP} and Proposition \ref{prop_UwhGV} to prove part (a). To prove part (b), use Lemma \ref{lem_trans} and part (a).
 \hfill\halmos

%\bigskip
%Although we can not relate the minimal bases of a rational matrix and its linearizations in general, we can say that the number of their minimal indices is the same.

%---------------------------------------------------------------------------------------------
\section{Minimal indices of strong linearizations of rational matrices}\label{Sec_minimalindicesstrong}
%-------------------------------------------------------------------------------------------
In this section we begin to study the relationship between  the
minimal indices of a rational matrix and the minimal indices
of its strong linearizations. As discussed in \cite[Remark 3.5]{AmDoMaZa18}, strong linearizations are particular cases of linearizations and, therefore, we know that the number of right (resp., left) minimal indices of a rational matrix and of its strong linearizations coincide. However, we will show in this section that it is not possible to obtain the right (resp., left) minimal indices of a rational matrix from those of its strong linearizations in general. Nevertheless,  we will prove in Theorem \ref{thm:summinindex} that the total sum of the right and left minimal indices of a rational matrix can be easily obtained from the total sum of the right and left minimal indices of any of its strong linearizations. { It has been recently shown in \cite{DasAl19laa,DaAl19_2,DaAl20} that for the families of Fiedler-like and affine spaces of strong linearizarions of {\em square} rational matrices it is possible to recover easily the minimal bases and indices of the rational matrix from these linerizations. We prove in Section \ref{Sec_blockminimalbases} that the same is possible for any rational matrix, i.e., possibly rectangular, from its strong block minimal bases linearizations. As corollaries, analogous results are proved in Section \ref{Sec_M12minimalbases} for $\mathbb{M}_1$ and $\mathbb{M}_2$-strong linearizations.}

We start by recalling the definition of strong linearization of a rational matrix.

\begin{definition}(\cite[Definition 3.4]{AmDoMaZa18})\label{def_stronglin}
Let $G(\la) \in\FF(\la)^{p\times m}$. Let $q_1$ be its first invariant order at infinity and $g=\min(0,q_1)$.
Let $n=\nu(G(\la))$. A strong linearization of $G(\la)$ is a linear polynomial matrix
$$%\begin{equation}\label{eq_lin_inf}
L(\la)=\begin{bmatrix}
A_1 \la +A_0 &B_1 \la +B_0\\-(C_1 \la +C_0)&D_1 \la +D_0
\end{bmatrix}\in\efe[\la]^{(n+q)\times (n+r)}
$$%\end{equation}
such that the following conditions hold:
\begin{itemize}
\item[\rm{(a)}] if $n>0$ then $\det(A_1\la+A_0)\neq 0$, and
\item[\rm{(b)}] if $\wh{G}(\la)=(D_1\la+D_0)+(C_1\la+C_0)(
A_1\la+A_0)^{-1}(B_1\la+B_0)$, $\wh{q}_{1}$ is its first invariant order at infinity and $\wh{g}=\min(0,\wh{q}_1)$ then:
\begin{itemize}
\item [\rm{(i)}] there are integers $s_1,s_2\geq 0$ and unimodular matrices $U_1(\la)\in\FF[\la]^{(p+s_1)\times (p+s_1)}$ and
$U_2(\la)\in\FF[\la]^{(m+s_1)\times (m+s_1)}$ so that $s_1-s_2=q-p=r-m$ and
\[
U_1(\la)\diag(G(\la),I_{s_1})U_2(\la)=\diag(\wh{G}(\la),I_{s_2})\text{, and}
\]
\item [\rm{(ii)}] there are biproper matrices $B_1(\la)\in\FF_{pr}(\la)^{(p+s_1)\times (p+s_1)}$
and $B_2(\la)\in\FF_{pr}(\la)^{(m+s_1)\times (m+s_1)}$ such that
\[
B_1(\la)\diag(\la^{g}G(\la),I_{s_1})B_2(\la)=\diag(\la^{\wh{g}}\wh{G}(\la),I_{s_2}).
\]
\end{itemize}
\end{itemize}
\end{definition}

As in the case of linearizations, we can also assume without loss of generality that $s_1=s$ and $s_2=0$ in the definition of strong linearizations. We will adopt such assumption in the rest of the paper.

\begin{rem} {\rm As commented in \cite[Remark 3.5]{AmDoMaZa18}, the requirement $n = \nu (G(\la))$ in Definition \ref{def_stronglin} might seem very restrictive. Thus, it is worth to emphasize that such requirement may be replaced by the assumptions that $L(\la)$ is a minimal polynomial system matrix and $A_1$ is invertible when $n>0$, as a consequence of the discussion in \cite[Remark 3.5]{AmDoMaZa18}, which are more direct requirements. We have decided to state Definition \ref{def_stronglin} exactly as in \cite{AmDoMaZa18} in order to avoid confusions.
}
\end{rem}

Recall that any rational matrix can be written uniquely as $G(\la)=D(\la)+G_{sp}(\la)$ with $D(\la)$ a polynomial matrix and $G_{sp}(\la)$ a strictly proper matrix. Moreover, if $D(\la)\neq 0$ then the first invariant order at infinity of $G(\la)$, $q_1$, is equal to $-\deg(D(\la))$; otherwise, if $G(\la)$ is strictly proper, $q_1>0$. We define
\begin{equation}\label{eq_d}
d=-\min(0,q_1)=\left\{\begin{array}{ll}\deg(D(\la))&\text{if }D(\la)\neq 0\\0&\text{if }D(\la)=0\end{array}\right..
\end{equation}
Notice that $g$ in Definition \ref{def_stronglin} is equal to $-d$.

We show now with Example \ref{ex} that the minimal indices of a strong linearization of a rational matrix may be arbitrarily different than the minimal indices of the rational matrix in general. In order to develop Example \ref{ex}, we present the following lemma first.

\begin{lemma} \label{lemm.auxxmin}
Let
$$
K_u(\la)=\begin{bmatrix}
1&\la&&&\\&1&\la&&\\&&\ddots&\ddots&\\&&&1&\la
\end{bmatrix}\in\efe[\la]^{u\times(u+1)}
$$
for any positive integer $u$ and let $0_{u,1}$ be the $u\times 1$ zero matrix. Then,
\begin{itemize}
\item [\rm{(i)}] $K_u(\la)$ is unimodularly equivalent to $\begin{bmatrix}
I_u&0_{u,1}
\end{bmatrix}$.
\item [\rm{(ii)}] $\la^{-1}K_u(\la)$ is equivalent at infinity to $\begin{bmatrix}
I_u&0_{u,1}
\end{bmatrix}$.
\end{itemize}
\end{lemma}

\textbf{Proof}.- In order to prove (i), multiply $K_u(\la)$ on the right by the unimodular matrix
$$
\begin{bmatrix}
1&-\la&\la^2&(-\la)^3&\cdots&(-\la)^u\\
&1&-\la&\la^2&\cdots&(-\la)^{u-1}\\
&&\ddots&\ddots&\ddots&\vdots\\
&&&1&-\la&\la^2\\
&&&&1&-\la\\
&&&&&1
\end{bmatrix}.
$$
To prove (ii), multiply $\la^{-1}K_u(\la)$ on the right by the biproper matrix
$$
\begin{bmatrix}
0&0&0&\cdots&0&1\\
1&0&0&\cdots&0&-1/\la\\
-1/\la&1&0&\cdots&0&(-1/\la)^{2}\\
(-1/\la)^{2}&-1/\la&1&\cdots&0&(-1/\la)^{3}\\
\vdots&\vdots&&\ddots&\vdots&\vdots\\
%&&&1&-\la&\la^2\\
(-1/\la)^{u-1}&(-1/\la)^{u-2}&&\cdots&1&(-1/\la)^{u}
\end{bmatrix}.
$$

\hfill\halmos

\begin{example}\label{ex}
{\rm Let $G(\la)=\begin{bsmallmatrix}
\la+\la^{-1}&0\\0&0
\end{bsmallmatrix}\in\efe(\la)^{2\times2}$.
We may consider infinitely many strong linearizations of $G(\la)$. Let
$$
L_{\epsilon,\eta}(\la)=\left[\begin{array}{c|ccc}
\la&1&&\\\hline-1&\la&&\\&&K_{\epsilon}(\la)&\\&&&K_{\eta}(\la)^T
\end{array}\right]\in\efe[\la]^{(1+(2+\epsilon+\eta))\times (1+(2+\epsilon+\eta))}.
$$
We prove now that for each pair of positive integers $\epsilon$ and $\eta$, $L_{\epsilon,\eta}(\la)$ is a strong linearization of $G(\la)$. First, notice that $L_{\epsilon,\eta}(\la)$ is a minimal polynomial system matrix with  transfer function matrix
$$
\begin{array}{ll}\widehat{G}_{\epsilon,\eta}(\la)&= %\begin{bmatrix}
%\la&&\\&K_{\epsilon}(\la)&\\&&K_{\eta}(\la)^T
%\end{bmatrix}+\begin{bmatrix}1\\0\\\vdots\\0\end{bmatrix}\la^{-1}\begin{bmatrix}1&0&\cdots&0\end{bmatrix}
%\\
%&=
\begin{bmatrix}
\la+\la^{-1}&&\\&K_{\epsilon}(\la)&\\&&K_{\eta}(\la)^T
\end{bmatrix}.
\end{array}
$$
Using Lemma \ref{lemm.auxxmin}, it is easy to prove that $\widehat{G}_{\epsilon,\eta}(\la)$ is unimodularly equivalent to
$$
\begin{bmatrix}
\la+\la^{-1}&&&\\&I_{\epsilon}&0_{\epsilon, 1}&\\&&&I_{\eta}\\&&&0_{1,\eta}
\end{bmatrix},
$$
which is unimodularly equivalent to
$
\begin{bsmallmatrix}
G(\la)&0\\0&I_{\epsilon+\eta}
\end{bsmallmatrix}.
$
Thus, $L_{\epsilon,\eta}(\la)$ is a linearization of $G(\la)$. Furthermore, $G(\la)$ can be written as $$G(\la)=\begin{bmatrix}
\la&0\\0&0
\end{bmatrix}+\begin{bmatrix}
\la^{-1}&0\\0&0
\end{bmatrix}
$$
and $\widehat{G}_{\epsilon,\eta}(\la)$ can be written as
$$
\widehat{G}_{\epsilon,\eta}(\la)=\begin{bmatrix}
\la&&\\&K_{\epsilon}(\la)&\\&&K_{\eta}(\la)^T
\end{bmatrix}+\begin{bmatrix}
\la^{-1}&&\\&0&\\&&0
\end{bmatrix}.
$$
Therefore, with the notation of Definition \ref{def_stronglin}, $g=\widehat{g}=-1$. The matrix $\la^{-1}\widehat{G}_{\epsilon,\eta}(\la)$ is
$$
\begin{bmatrix}
1+\la^{-2}&&\\&\la^{-1}K_{\epsilon}(\la)&\\&&\la^{-1}K_{\eta}(\la)^T
\end{bmatrix},
$$
which, by  Lemma \ref{lemm.auxxmin}, is equivalent at infinity to
$$
\begin{bmatrix}
1+\la^{-2}&&&\\&I_{\epsilon}&0_{\epsilon, 1}&\\&&&I_{\eta}\\&&&0_{1,\eta}
\end{bmatrix}
\text{ and to }
\begin{bmatrix}
\la^{-1}G(\la)&0\\0&I_{\epsilon+\eta}
\end{bmatrix}.
$$
Hence, $L_{\epsilon,\eta}(\la)$ is a strong linearization of $G(\la)$. Notice that the} unique right minimal index of $G(\la)$ is 0 and the unique left minimal index of $G(\la)$ is 0 as well, while the unique right minimal index of $L_{\epsilon,\eta}(\la)$ is $\epsilon$ and the unique left minimal index of $L_{\epsilon,\eta}(\la)$ is $\eta$. Thus, strong linearizations do not preserve minimal indices.

\end{example}

Denote by $\mu(G(\la))$ the sum of the right and left minimal indices of a rational matrix $G(\la)$.
Our next goal is to analyze how this is related with the sum of the right and left minimal indices of any of its strong linearizations.  In order to study this relationship, we will make use of Van Dooren's index sum theorem, proved for the first time in \cite[Theorem 3]{VeDoka79}, and that we rewrite in a way convenient for our purposes in Lemma \ref{lem_ist}. Interested readers are referred to the recent paper \cite{AnDoHoMac19} for more information on this fundamental result.

\begin{lemma}\label{lem_ist}
Let $G(\la)\in\efe(\la)^{p\times m}$ be any rational matrix with finite Smith--McMillan form $\diag\left(\frac{\epsilon_1(\la)}{\psi_1(\la)},\ldots,\frac{\epsilon_r(\la)}{\psi_r(\la)},0_{p-r,m-r}\right)$. Let $q_1\leq\cdots\leq q_r$ be its invariant orders at infinity. Then
\begin{equation}\label{eq_muG}
\mu(G(\la))=\sum_{i=1}^r\deg(\psi_i(\la))-\sum_{i=1}^r\deg(\epsilon_i(\la))-\sum_{i=1}^r q_i.
\end{equation}
\end{lemma}

\textbf{Proof}.- By the index sum theorem (see \cite[Theorem 3]{VeDoka79} or \cite[Theorem 6.5-11]{Kailath80}) $\mu(G(\la))$ is equal to the total number of poles (finite and at infinity) of $G(\la)$ minus the total number of zeros (finite and at infinity) of $G(\la)$. The total number of finite zeros of $G(\la)$ is the sum of all partial multiplicities of all finite zeros of $G(\la)$, that is, $\sum_{i=1}^r\deg(\epsilon_i(\la))$. In the same way, the total number of finite poles of $G(\la)$ is the sum of all partial multiplicities of all finite poles of $G(\la)$, i.e., $\sum_{i=1}^r\deg(\psi_i(\la))$. Therefore, the total number of finite poles minus the total number of finite zeros is $\sum_{i=1}^r\deg(\psi_i(\la))-\sum_{i=1}^r\deg(\epsilon_i(\la))$. On the other hand, the total number of infinite poles minus the total number of infinite zeros is $-\sum_{i=1}^r q_i$ since the positive $q_i$ are the orders of the infinite zeros while minus the negative $q_i$ are the orders of the infinite poles. Thus, equation (\ref{eq_muG}) is obtained.
\hfill\halmos

\bigskip
Let $G(\la)\in\efe(\la)^{p\times m}$ be any rational matrix, let $d$ be defined as in (\ref{eq_d}) and let
\begin{equation}\label{eq_lin_inf}
L(\la)=\begin{bmatrix}
A_1 \la +A_0 &B_1 \la +B_0\\-(C_1 \la +C_0)&D_1 \la +D_0
\end{bmatrix}\in\efe[\la]^{(n+(p+s))\times (n+(m+s))}
\end{equation}
be a linear minimal polynomial system matrix
with $A_1$ invertible if $n>0$. We say that $L(\la)$ preserves the finite and infinite structures of poles and zeros of $G(\la)$ if
the following conditions simultaneously hold:
\begin{itemize}
\item[(i)] the finite poles of $G(\la)$ are the finite zeros of $A_1\la+A_0$, with the same partial
multiplicities in both matrices,
\item[(ii)] the finite zeros of $G(\la)$ are the finite zeros of $L(\la)$, with the
same partial multiplicities, and
\item[(iii)]  the number and orders of the infinite zeros of $\la^{-1}L(\la)$
are the same as the number and orders of the infinite zeros of $\la^{-d}G(\la)$ if
$D_1+C_1A_1^{-1}B_1\neq 0$ or of $\diag(\la^{-1}I_s,\la^{-d-1}G(\la))$ otherwise.
\end{itemize}

\begin{theorem}(\cite[Theorem 3.10]{AmDoMaZa18})\label{thm_spcharstlin}
Let $G(\la)\in\FF(\la)^{p\times m}$ and $n=\nu(G(\la))$. Let $L(\la)$ be the pencil of \eqref{eq_lin_inf}.
%\[
%L(\la)=\begin{bmatrix}
%A_1 \la +A_0 &B_1 \la +B_0\\-(C_1 \la +C_0)&D_1 \la +D_0
%\end{bmatrix}\in\efe[\la]^{(n+(p+s))\times (n+(m+s))}.
%\]
Then $L(\la)$ is a strong linearization of $G(\la)$ if and only if the following two conditions hold:
\begin{enumerate}
\item[\rm (I)] $\mbox{\rm dim} \, \mathcal{N}_r (G(\la)) = \mbox{\rm dim} \,  \mathcal{N}_r ( L(\la) )$ (which is equivalent to $\mbox{\rm dim} \, \mathcal{N}_{\ell} (G(\la)) = \mbox{\rm dim} \,  \mathcal{N}_{\ell} ( L(\la) )$), and

\item[\rm (II)] $L(\la)$ preserves the finite and infinite structures of poles and zeros of $G(\la)$.
\end{enumerate}
\end{theorem}

The following result relates the invariant orders at infinity of  a rational matrix and its a strong linearization.  We remark that although Lemma \ref{lem_qL} was not explicitly stated in \cite{AmDoMaZa18}, it is related to discussions in \cite[pp. 1682--1683]{AmDoMaZa18}.

\begin{lemma}\label{lem_qL}
Let $G(\la)\in\efe(\la)^{p\times m}$ be any rational matrix with invariant orders at infinity $q_1\leq\cdots\leq q_r$ and $d$ be defined as in (\ref{eq_d}). Let $L(\la)$ of \eqref{eq_lin_inf}
be any strong linearization of $G(\la)$ and $q_1^L\leq\ldots\leq q_{\ell}^L$ be the invariant orders at infinity of $L(\la)$. Then $\ell=n+s+r$ and
\begin{itemize}
\item[\rm{(i)}] If $D_1+C_1A_1^{-1}B_1\neq 0$ then  $q_i^L=-1$ for $i=1,\ldots, n+s$, and $q_{n+s+i}^L=q_i+d-1$ for $i=1,\ldots,r$.

\item[\rm{(ii)}] If $n>0$ and $D_1+C_1A_1^{-1}B_1=0$ then
$q_i^L=-1$ for $i=1,\ldots,n$, $q_{n+i}^L=0$ for $i=1,\ldots,s$, and $q_{n+s+i}^L=q_i+d$ for $i=1,\ldots,r$.

\item[\rm{(iii)}] If $n=0$ and $D_1=0$ then $L(\la)=D_0$, $q_i^L=0$ for $i=1,\ldots,s+r$, and $q_i=-d$ for $i=1,\ldots,r$.
\end{itemize}
\end{lemma}

\textbf{Proof}.- By Theorem \ref{thm_spcharstlin} (I) and the rank-nullity theorem, $\ell=n+s+r$ is the rank of $L(\la)$. As $q_1\leq\cdots\leq q_r$ are the invariant orders at infinity of $G(\la)$, there exist two biproper matrices $B_1(\la)\in\efe_{pr}(\la)^{p\times p}$ and $B_2(\la)\in\efe_{pr}(\la)^{m\times m}$ such that
\begin{equation}\label{eq_Ginf}
G(\la)=B_1(\la)\diag\left(\left(\frac{1}{\la}\right)^{q_1},\ldots,
\left(\frac{1}{\la}\right)^{q_r},0_{p-r,m-r}\right)B_2(\la).
\end{equation}
%Moreover, we know by the previous theorem and the rank-nullity theorem that $\rank L(\la)=
%n+m+s-\mbox{\rm dim} \,  \mathcal{N}_r ( L(\la) )=n+m+s-\mbox{\rm dim} \,  \mathcal{N}_r ( G(\la) )=
%n+s+r$.
We distinguish two cases:

Suppose first that $D_1+C_1A_1^{-1}B_1\neq 0$. By  Theorem \ref{thm_spcharstlin} again, the number and orders of the infinite zeros of $\la^{-1}L(\la)$
are the same as the number and orders of the infinite zeros of $\la^{-d}G(\la)$. Since $\la^{-1}L(\la)$ and $\la^{-d}G(\la)$ are both proper rational matrices and $\rank L(\la)-\rank G(\la)=n+s$, $\la^{-1}L(\la)$ must be equivalent at infinity to $\begin{bsmallmatrix}
\la^{-d}G(\la)&0\\0&I_{n+s}\end{bsmallmatrix}$. Thus $L(\la)$ is equivalent at infinity to $\begin{bmatrix}
\la^{-d+1}G(\la)&0\\0&\la I_{n+s}\end{bmatrix}$, that is, there exist two biproper matrices $B_3(\la)\in\efe_{pr}(\la)^{(p+n+s)\times (p+n+s)}$ and $B_4(\la)\in\efe_{pr}(\la)^{(m+n+s)\times (m+n+s)}$ such that
$$
L(\la)=B_3(\la)\begin{bmatrix}
\la^{-d+1}G(\la)&0\\0&\la I_{n+s}\end{bmatrix}B_4(\la)=B_3(\la)\la^{-d+1}\begin{bmatrix}
G(\la)&0\\0&\la^d I_{n+s}\end{bmatrix}B_4(\la).
$$
Put $\overline{B}_1(\la)=B_3(\la)\begin{bsmallmatrix}
B_1(\la)&0\\0&I_{n+s}\end{bsmallmatrix}$ and $\overline{B}_2(\la)=\begin{bsmallmatrix}
B_2(\la)&0\\0&I_{n+s}\end{bsmallmatrix}B_4(\la)$, which are biproper matrices.
Using (\ref{eq_Ginf}),
\begin{align*}
L(\la)& =\overline{B}_1(\la)\la^{-d+1}\begin{bmatrix}
\diag\left((\frac{1}{\la})^{q_1},\ldots,
(\frac{1}{\la})^{q_r},0\right)&0\\0&\la^d I_{n+s}\end{bmatrix}
%\begin{bmatrix}
%B_2(\la)&0\\0&I_{n+s}\end{bmatrix}B_4(\la)=
\overline{B}_2(\la) \\
%B_3(\la)\begin{bmatrix}
%B_1(\la)&0\\0&I_{n+s}\end{bmatrix}
& = \overline{B}_1(\la)
\begin{bmatrix}
\diag\left((\frac{1}{\la})^{q_1+d-1},\ldots,
(\frac{1}{\la})^{q_r+d-1},0\right)&0\\0&(\frac{1}{\la})^{-1} I_{n+s}\end{bmatrix}
%\begin{bmatrix}
%B_2(\la)&0\\0&I_{n+s}\end{bmatrix}B_4(\la).
\overline{B}_2(\la).
\end{align*}
Notice, by (\ref{eq_d}), that $q_1+d\geq 0$. Therefore $-1\leq q_1+d-1\leq\cdots\leq q_r+d-1$. Thus, $q_i^L=-1$ for $i=1,\ldots, n+s$, and $q_{n+s+i}^L=q_i+d-1$ for $i=1,\ldots,r$.

Suppose now that $D_1+C_1A_1^{-1}B_1=0$. By Theorem \ref{thm_spcharstlin}, the number and orders of the infinite zeros of $\la^{-1}L(\la)$
are the same as those of $\diag(\la^{-1}I_s,\la^{-d-1}G(\la))$. As both matrices are proper and their rank difference is $n$, $\la^{-1}L(\la)$ must be equivalent at infinity to $\begin{bsmallmatrix}
\la^{-d-1}G(\la)&0&0\\0&\la^{-1}I_{s}&0\\0&0&I_n\end{bsmallmatrix}$. Thus $L(\la)$ is equivalent at infinity to $\begin{bsmallmatrix}
\la^{-d}G(\la)&0&0\\0&I_{s}&0\\0&0&\la I_n\end{bsmallmatrix}$, that is, there exist two biproper matrices $B_5(\la)\in\efe_{pr}(\la)^{(p+n+s)\times (p+n+s)}$ and $B_6(\la)\in\efe_{pr}(\la)^{(m+n+s)\times (m+n+s)}$ such that %$L(\la)=$
\begin{align*}
L(\la) &=
%B_5(\la)\begin{bmatrix}
%\la^{-d}G(\la)&0&0\\0&I_{s}&0\\0&0&\la I_n\end{bmatrix}B_6(\la)
%\\
%& =
B_5(\la)\la^{-d}\begin{bmatrix}
G(\la)&0&0\\0&\la^{d}I_{s}&0\\0&0&\la^{d+1} I_n\end{bmatrix}B_6(\la).
\end{align*}
By using (\ref{eq_Ginf}) and proceeding as in the previous case,  if $n>0$ then the invariant orders at infinity of $L(\la)$ are $q_i^L=-1$ for $i=1,\ldots,n$, $q_{n+i}^L=0$ for $i=1,\ldots,s$, and $q_{n+s+i}^L=q_i+d$ for $i=1,\ldots,r$. Otherwise, if $n=0$ then $D_1=0$, $L(\la)=D_0$ and, therefore, $q_i^L=0$ for $i=1,\ldots,s+r$. Moreover, since $D_0=B_5(\la)\begin{bsmallmatrix}
\la^{-d}G(\la)&0\\0&I_{s}\end{bsmallmatrix}B_6(\la)$, the invariant orders at infinity of $\la^{-d}G(\la)$ must be 0 and, in consequence, $q_i=-d$ for $i=1,\ldots,r$.
\hfill\halmos

\bigskip
The following lemma gives $\mu(L(\la))$, the sum of the right and left minimal indices of a strong linearization $L(\la)$ of a rational matrix $G(\la)$, in terms of the spectral invariants of $G(\la)$.

\begin{lemma}\label{lem_muL}
Let $G(\la)\in\efe(\la)^{p\times m}$ be any rational matrix with $\epsilon_1(\la),\ldots ,\epsilon_r(\la)$ as numerators in its finite Smith--McMillan form and with $q_1\leq\cdots\leq q_r$ as invariant orders at infinity. Let $d$ be defined as in (\ref{eq_d}). Let $L(\la)$ of \eqref{eq_lin_inf}
be any strong linearization of $G(\la)$.
\begin{itemize}
\item[\rm{(i)}] If $D_1+C_1A_1^{-1}B_1\neq 0$ then
$$
\begin{array}{ll}
\mu(L(\la))&=s+r(1-d)+n-\sum_{i=1}^r\deg(\epsilon_i(\la))-\sum_{i=1}^{r} q_i.
\end{array}
$$
\item[\rm{(ii)}] If $n>0$ and $D_1+C_1A_1^{-1}B_1=0$ then
$$
\begin{array}{ll}
\mu(L(\la))&=-dr+n-\sum_{i=1}^r\deg(\epsilon_i(\la))-\sum_{i=1}^{r} q_i.
\end{array}
$$
\item[\rm{(iii)}] If $n=0$ and $D_1=0$ then $L(\la)=D_0$, $\mu(L(\la))=0$, and $\epsilon_i(\la)=1$ for $i=1,\ldots,r$.
\end{itemize}
\end{lemma}

\textbf{Proof}.- We aim to apply Lemma \ref{lem_ist} to $L(\la)$. As seen in Lemma \ref{lem_qL}, $\rank L(\la)=n+s+r$. Since $L(\la)$ is a polynomial matrix it has no finite poles. Moreover, by Theorem \ref{thm_spcharstlin}, its total number of finite zeros is $\sum_{i=1}^r\deg(\epsilon_i(\la))$. Denote by $q_i^L$, $i=1,\ldots,n+s+r$, the invariant orders at infinity of $L(\la)$. By Lemma \ref{lem_ist}, $$\mu(L(\la))=-\sum_{i=1}^r\deg(\epsilon_i(\la))-\sum_{i=1}^{n+s+r} q_i^L.$$
By Lemma \ref{lem_qL}:
\begin{itemize}
\item[\rm{(i)}] If $D_1+C_1A_1^{-1}B_1\neq 0$ then
$$
\begin{array}{ll}
\mu(L(\la))&=-\sum_{i=1}^r\deg(\epsilon_i(\la))-(\sum_{i=1}^{n+s} (-1)+\sum_{i=1}^r (q_i+d-1))\\
&=-\sum_{i=1}^r\deg(\epsilon_i(\la))+n+s+r-dr-\sum_{i=1}^r q_i.
\end{array}
$$
\item[\rm{(ii)}] If $n>0$ and $D_1+C_1A_1^{-1}B_1=0$ then
$$
\begin{array}{ll}
\mu(L(\la))&=-\sum_{i=1}^r\deg(\epsilon_i(\la))-(\sum_{i=1}^{n} (-1) +\sum_{i=1}^{r} (q_i+d))\\
 &=-\sum_{i=1}^r\deg(\epsilon_i(\la))+n-dr-\sum_{i=1}^{r} q_i.
\end{array}
$$
\item[\rm{(iii)}] If $n=0$ and $D_1=0$ then $L(\la)=D_0$ and $\mu(L(\la))=-\sum_{i=1}^r\deg(\epsilon_i(\la))-0$. But since $L(\la)$ is constant its total number of finite zeros is 0 and, therefore, $\epsilon_i(\la)=1$ for $i=1,\ldots,r$. \hfill\halmos
\end{itemize}
%\hfill\halmos

\bigskip
Finally, the following result shows the relationship between the sum of the right and left minimal indices of a rational matrix and of its strong linearizations.

\begin{theorem} \label{thm:summinindex}
Let $G(\la)\in\efe(\la)^{p\times m}$ be any rational matrix of rank $r$. Let $d$ be defined as in (\ref{eq_d}). Let $$L(\la)=\begin{bmatrix}
A_1 \la +A_0 &B_1 \la +B_0\\-(C_1 \la +C_0)&D_1 \la +D_0
\end{bmatrix}\in\efe[\la]^{(n+(p+s))\times (n+(m+s))}$$
be any strong linearization of $G(\la)$. Then
$$
\mu(G(\la))=\left\{\begin{array}{ll}\mu(L(\la))+dr-(r+s),&\text{if }D_1+C_1A_1^{-1}B_1\ne 0\\\mu(L(\la))+dr,&\text{if }n>0\text{ and }D_1+C_1A_1^{-1}B_1=0\\dr, &\text{if }n=0\text{ and }D_1=0\end{array}\right..
$$
\end{theorem}

\textbf{Proof}.- Let $\diag\left(\frac{\epsilon_1(\la)}{\psi_1(\la)},\ldots,\frac{\epsilon_r(\la)}{\psi_r(\la)},0_{p-r,m-r}\right)$ be the finite Smith--McMillan form of $G(\la)$ and $q_1\leq\cdots\leq q_r$ be its invariant orders at infinity.
By definition of strong linearization, $n=\nu(G(\la))$. Moreover, $\nu(G(\la))=\sum_{i=1}^r\deg(\psi_i(\la))$ and, therefore, $n=\sum_{i=1}^r\deg(\psi_i(\la))$. By using Lemma \ref{lem_ist}, $\mu(G(\la))=n-\sum_{i=1}^r\deg(\epsilon_i(\la))-\sum_{i=1}^r q_i$. Now, by Lemma \ref{lem_muL}:
\begin{itemize}
\item[\rm{(i)}] If $D_1+C_1A_1^{-1}B_1\neq 0$ then
$
\mu(L(\la))=s+r(1-d)+\mu(G(\la)).
$
\item[\rm{(ii)}] If $n>0$ and $D_1+C_1A_1^{-1}B_1=0$ then
$
\mu(L(\la))=-dr+\mu(G(\la)).
$
\item[\rm{(iii)}] If $n=0$ and $D_1=0$ then $\mu(L(\la))=0$ and, by Lemmas \ref{lem_qL} and \ref{lem_muL}, $\mu(G(\la))=dr$. \hfill\halmos
\end{itemize}

%{\color{red} Esto nos permite sacar alguna conclusi�n o deducir algo m�s?}
\begin{example}
{\rm We show that, certainly, the previous result is satisfied for the matrices in Example \ref{ex}. It was proved that the matrices $L_{\epsilon,\eta}(\la)$ are strong linearizations of $G(\la)=\begin{bsmallmatrix}	
	\la+\la^{-1}&0\\0&0\end{bsmallmatrix}$. Notice that, under the same notation as above, $r=1$, $d=1$, $n=1$, $s=\epsilon+\eta$, $A_1=1, B_1=0, C_1=0$ and $D_1+C_1A_1^{-1}B_1\neq 0$. As we proved $\mu(G(\la))=0$ and $\mu(L_{\epsilon,\eta}(\la))=\epsilon+\eta$. Thus, $\mu(G(\la))=\mu(L_{\epsilon,\eta}(\la))+dr-(r+s)$, as claimed.
	}	
\end{example}

%---------------------------------------------------------------------------------------------
\section{Minimal bases and indices of strong block minimal bases linearizations of rational matrices}\label{Sec_blockminimalbases}
%-------------------------------------------------------------------------------------------

The aim of this section is to study the relationship between  the
minimal bases and indices of a rational matrix and the minimal bases and indices
of its strong block minimal bases linearizations. { This family of strong linearizations is a rather general family introduced in \cite[Theorem 5.11]{AmDoMaZa18}. It will be shown in Section \ref{sec.Fiedler} that the families of Fiedler-like linearizations of rational matrices introduced in \cite{AlBe16,AlBe18,DasAl19laa,DaAl20} are, modulo permutations, particular instances of strong block minimal bases linearizations. Actually, this is a consequence of the corresponding results for polynomial matrices in \cite{BuDoPe18} and \cite[Lemma 2.7]{DoMaQu19}. Moreover, the strong block minimal bases linearizations are closely connected to those introduced in \cite{DoMaQu19}. In contrast to Fiedler-like linearizations, affine spaces of linearizations \cite{DaAl19_2} and the linearizations in \cite{DoMaQu19}, which are only defined for square rational matrices, strong block minimal bases linearizations are valid for general rectangular rational matrices.} Strong block minimal bases linearizations of rational matrices are built on strong block minimal bases linearizations of polynomial matrices, presented previously in \cite[Definition 3.1]{Dopico:2018:BlockKronecker} (see \cite{Dopico:2016:BlockKronecker} for an expanded version of this latter reference). In order to introduce these families of linearizations and prove the results in this section, we need to recall first a number of concepts in the next paragraphs.

A matrix polynomial $N(\la)\in\efe[\la]^{m\times l}$ with $m<l$ is a minimal basis if the columns of $N(\la)^T$ form a minimal basis of the subspace they span. Moreover, two matrix polynomials $K(\la)\in\efe[\la]^{m_1\times l}$ and $N(\la)\in\efe[\la]^{m_2\times l}$ are dual minimal bases  if they are both minimal bases satisfying $m_1+m_2=l$ and $K(\la)N(\la)^T=0$ (see \cite{Dopico:2018:BlockKronecker,Fo75}).

Let us recall the definition of strong block minimal bases pencils associated to a polynomial matrix (see \cite[Definition 3.1 and Theorem 3.3]{Dopico:2018:BlockKronecker} or \cite[Definition 5.2]{AmDoMaZa18}).  Let $P(\la) \in \FF[\la]^{p \times m}$ be a polynomial matrix. A strong block minimal bases pencil associated to $P(\la)$ is a linear polynomial matrix with the following structure
\begin{equation}
  \label{eq:minbaspencil}
  \begin{array}{cl}
  \mathcal{L}(\la) =
  \left[
    \begin{array}{cc}
      M(\la) & K_2 (\la)^T \\
      K_1 (\lambda) &0
      \end{array}
    \right]&
    \begin{array}{l}
      \left. \vphantom{K_2 (\la)^T} \right\} {\scriptstyle p + \widehat{p}}\\
      \left. \vphantom{K_1 (\la)} \right\} {\scriptstyle \widehat{m}}
    \end{array}\\
%    &
    \hphantom{\mathcal{L}(\la) =}
    \begin{array}{cc}
      \underbrace{\hphantom{K_1 (\lambda)}}_{\scriptstyle m + \widehat{m}} & \underbrace{\hphantom{K_2 (\la)^T}}_{\widehat{p}}
    \end{array}
  \end{array}
  \>,
\end{equation}
where $K_1(\la) \in \FF[\la]^{\widehat{m} \times (m + \widehat{m})}$ (respectively $K_2(\la) \in \FF[\la]^{\widehat{p} \times (p + \widehat{p})}$) is a minimal basis with all its row degrees equal to $1$ and with the row degrees of a minimal basis $N_1(\la) \in \FF[\la]^{m \times (m + \widehat{m})}$ (respectively $N_2(\la) \in \FF[\la]^{p \times (p + \widehat{p})}$) dual to $K_1(\la)$ (respectively $K_2(\la)$) all equal, and such that
\begin{equation} \label{eq:Dpolinminbaslin}
P(\la) = N_2(\la) M(\la) N_1(\la)^T.
\end{equation}
If, in addition, $\deg(P(\la)) = \deg(N_2(\la)) +  \deg(N_1(\la)) + 1$ then $\mathcal{L}(\la)$ is said to be a strong block minimal bases pencil associated to $P(\la)$ with sharp degree. The key property is that any strong block minimal bases pencil associated to $P(\la)$ is a strong linearization of $P(\la)$ \cite[Theorem 3.3]{Dopico:2018:BlockKronecker}.

Let $G(\la) = D(\la) + G_{sp} (\la)$ be the unique decomposition of $G(\la) \in \FF (\la)^{p\times m}$ into its polynomial part $D(\la) \in \FF[\la]^{p\times m}$ and its strictly proper part $G_{sp}(\la) \in \FF_{pr}(\la)^{p\times m}$, and let $G_{sp} (\la) = C (\la I_n - A)^{-1} B$ be a minimal order state-space realization of $G_{sp} (\la)$ with $n= \nu (G(\la))$. Assume\footnote{If $\deg (D(\la)) \leq 1$, then the polynomial system matrix $\begin{bsmallmatrix}
                \la I_n - A & B \\
                -C & D(\la)
              \end{bsmallmatrix}$ with transfer function matrix $G(\la)$ gives directly a strong linearization of $G(\la)$, as discussed in \cite{AmDoMaZa18}, and the idea of strong block minimal bases linearizations is of no interest.} that $\deg (D(\la)) >1$ and let (\ref{eq:minbaspencil})
be a strong block minimal bases pencil associated to $D(\la)$ with sharp degree, with $N_1(\la) \in \FF[\la]^{m \times (m + \widehat{m})}$ and $N_2(\la) \in \FF[\la]^{p \times (p + \widehat{p})}$ minimal bases dual to $K_1(\la)$ and $K_2(\la)$, respectively, such that $D(\la) = N_2 (\la) M(\la) N_1(\la)^T$. Let $\widehat{K}_1 \in \FF^{m \times (m + \widehat{m})}$, $\widehat{N}_1(\la) \in \FF[\la]^{\widehat{m} \times (m + \widehat{m})}$, $\widehat{K}_2 \in \FF^{p \times (p + \widehat{p})}$ and $\widehat{N}_2(\la) \in \FF[\la]^{\widehat{p} \times (p + \widehat{p})}$ be matrices such that for $i=1,2$
\begin{equation} \label{eq.unimodSBMBP}
U_i(\lambda) =
\begin{bmatrix}
K_i(\lambda) \\ \widehat{K}_i
\end{bmatrix}
\quad \mbox{and} \quad
U_i(\lambda)^{-1}=
\begin{bmatrix}
\widehat{N}_i(\lambda)^T & N_i(\lambda)^T
\end{bmatrix}
\end{equation}
are unimodular (see in \cite[Lemma 5.5]{AmDoMaZa18} the result that guaratees that all these matrices exist and are well-defined). Let $T, S \in\FF^{n\times n}$ be any nonsingular constant matrices. By \cite[Theorem 5.11]{AmDoMaZa18} the linear polynomial matrix
\begin{equation} \label{eq:ratstrongblockmin}
L(\la) = \left[
\begin{array}{c|cc}
T(\la I_n -A)S & \phantom{a} TB\widehat{K}_1 \phantom{a} & 0 \\ \hline \phantom{\Big|}
-\widehat{K}_2^T C S \phantom{\Big|}& M(\la) & K_2(\la)^T \\
0 & K_1 (\la) & 0
\end{array}
\right]
\end{equation}
is a strong linearization of $G(\la)$ and is called strong block minimal bases linearization of $G(\la)$.

Furthermore, by \cite[Theorem 5.7]{AmDoMaZa18}, there are matrices $X(\la)\in\efe[\la]^{\widehat{p}\times m}$ ($X(\la)=\wh{N}_2(\la)M(\la)N_1(\la)^T$), $Y(\la)\in\efe[\la]^{p\times\widehat{m}}$  ($Y(\la) = N_2 (\la) M(\la) \widehat{N}_1 (\la)^T$), and $Z(\la)\in\efe[\la]^{\widehat{p}\times\widehat{m}}$ ($Z(\la) = \widehat{N}_2 (\la) M(\la) \widehat{N}_1 (\la)^T$) such that
%Let us define the unimodular matrices
\begin{align} \label{eq.ffvu}
V (\la) &
%:= \begin{bmatrix}
%\widehat{N}_1(\lambda)^T & N_1(\lambda)^T & 0  \\
%0 & 0 & I_{\widehat{p}}
%\end{bmatrix}
%\begin{bmatrix}
%0 & I_{\widehat{m}} & 0 \\
%I_m & 0 & 0 \\
%-X(\lambda) & 0 & I_{\widehat{p}}
%\end{bmatrix}
=\begin{bmatrix}
N_1(\lambda)^T &\widehat{N}_1(\lambda)^T &  0  \\
-X(\la) & 0 & I_{\widehat{p}}
\end{bmatrix}\ \text{ and }\
%\\
U(\la) %&
%:=
%\begin{bmatrix}
%0 & I_p & -Y(\lambda)\\
%0 & 0 & I_{\widehat{m}} \\
%I_{\widehat{p}} & 0 & -Z(\lambda)
%\end{bmatrix}
%\begin{bmatrix}
%\widehat{N}_2(\lambda) & 0 \\ N_2(\lambda) & 0 \\0 & I_{\widehat{m}}
%\end{bmatrix}
=\begin{bmatrix}
N_2(\lambda) &-Y(\la)\\0&I_{\wh{m}}\\\widehat{N}_2(\lambda) &  -Z(\la)  \\
\end{bmatrix}
\end{align}
are unimodular matrices and
%where $X(\la) = \widehat{N}_2 (\la) M(\la) N_1 (\la)^T$, $Y(\la) = N_2 (\la) M(\la) \widehat{N}_1 (\la)^T$, and $Z(\la) = %\widehat{N}_2 (\la) M(\la) \widehat{N}_1 (\la)^T$. Then
%\begin{equation} \label{eq.unimodSBMBP}
$$
U(\la) \,\begin{bmatrix}M(\la)&K_2(\la)^T\\K_1(\la)&0\end{bmatrix}\, V(\la) = \diag( D(\la), I_{\widehat{m} + \widehat{p}}),
$$
 as can be easily checked through a direct matrix multiplication.
%\end{equation}
Moreover, $U(\la)\begin{bsmallmatrix}
-\wh{K}_2^TCS\\0
\end{bsmallmatrix}=%\begin{bmatrix}
%N_2(\lambda) &-Y(\la)\\0&I_{\wh{m}}\\\widehat{N}_2(\lambda) &  -Z(\la)  \\
%\end{bmatrix}\begin{bmatrix}
%-\wh{K}_2^TCS\\0
%\end{bmatrix}=
\begin{bsmallmatrix}
-CS\\0
\end{bsmallmatrix}
$ and
$\begin{bmatrix}TB\wh{K}_1&0\end{bmatrix}V(\la)=%\begin{bmatrix}TB\wh{K}_1&0\end{bmatrix}\begin{bmatrix}
%N_1(\lambda)^T &\widehat{N}_1(\lambda)^T &  0  \\
%-X(\la) & 0 & I_{\widehat{p}}
%\end{bmatrix}=
\begin{bmatrix}TB&0\end{bmatrix}$.
Thus,
$$
\begin{bmatrix}T^{-1}&0\\0&U(\la)\end{bmatrix}L(\la)\begin{bmatrix}S^{-1}&0\\0&V(\la)\end{bmatrix}=
\begin{bmatrix}\la I_n-A&B&0\\-C&D(\la)&0\\0&0&I_{\wh{m}+\wh{p}}\end{bmatrix}.
$$
Let $\wh{G}(\la)$ be the transfer function matrix of $L(\la)$, i.e.,
\begin{equation}\label{whG}
\wh{G}(\la)=\left[\begin{array}{cc}M(\la)+\wh{K}_2^TC(\la I_n-A)^{-1}B\wh{K}_1&K_2(\la)^T\\K_1(\la)&0\end{array}\right].
\end{equation}
%By \cite[Proposition 4.3]{AmDoMaZa18},
 Taking into account the developments above, a straightforward computation yields
\begin{equation} \label{eq.ffugvg}
U(\la)\wh{G}(\la)V(\la)=\diag(G(\la),I_{\wh{m}+\wh{p}}),
\end{equation}
which implies, among other properties, $\dim \mathcal{N}_r (\wh{G}(\la))=\dim \mathcal{N}_r (G(\la))$ and
$\dim \mathcal{N}_\ell (\wh{G}(\la))=\dim \mathcal{N}_\ell (G(\la))$, in agreement with the properties of any (strong) linearization of $G(\la)$.

%\begin{lemma}
%\begin{enumerate}
% \item If $H(\la)$ is a polynomial basis of $G(\la)$ then %$\left[\begin{array}{c}N_1(\lambda)^T\\-X(\la)\end{array}\right]H(\la)%$  is a polynomial basis of $\wh{G}(\la)$.

% \item If %$\wh{H}(\la)=\left[\begin{array}{c}\wh{H}_1(\la)\\\wh{H}_2(\la)\end{ar%ray}\right]$, with $\wh{H}_1(\la)\in\efe[\la]^{(m+\wh{m})\times l}$ %and $\wh{H}_2(\la)\in\efe[\la]^{\wh{p}\times l}$,  is a polynomial %basis of $\wh{G}(\la)$ then $\wh{K}_1\wh{H}_1(\la)$ is a polynomial %basis of $G(\la)$, $K_1(\la)\wh{H}_1(\la)=0$ and %$X(\la)\wh{K}_1\wh{H}_1(\la)+\wh{H}_2(\la)=0$.
% \end{enumerate}
%\end{lemma}

%\textbf{Proof}.-

%$V(\la)\left[\begin{array}{c}H(\la)\\0\end{array}\right]=\begin{bmatri%x}
%N_1(\lambda)^T &\widehat{N}_1(\lambda)^T &  0  \\
%-X(\la) & 0 & I_{\widehat{p}}
%\end{bmatrix}\left[\begin{array}{c}H(\la)\\0\end{array}\right]=
%\left[\begin{array}{c}N_1(\lambda)^T\\-X(\la)\end{array}\right]H(\la)$

%$V(\la)^{-1}\wh{H}(\la)=\begin{bmatrix}
%\wh{K}_1 &  0  \\K_1(\la)&0\\
%X(\la)\wh{K}_1 & I_{\widehat{p}}
%\end{bmatrix}\left[\begin{array}{c}\wh{H}_1(\la)\\\wh{H}_2(\la)\end{ar%ray}\right]=\left[\begin{array}{c}\wh{K}_1\wh{H}_1(\la)\\
%K_1(\la)\wh{H}_1(\la)\\X(\la)\wh{K}_1\wh{H}_1(\la)+\wh{H}_2(\la)\end{a%rray}\right]$. Thus $\wh{K}_1\wh{H}_1(\la)$ is a polynomial basis of %$G(\la)$, $K_1(\la)\wh{H}_1(\la)=0$ and %$X(\la)\wh{K}_1\wh{H}_1(\la)+\wh{H}_2(\la)=0$.
%\hfill\halmos

In order to investigate the relationship between the minimal bases and indices of a rational matrix and those of its strong block minimal bases linearizations, we prove Lemma \ref{minimalbasis}. This lemma first establishes the relationship between vectors in the right null-space of the rational matrix and in the right null-spaces of the transfer functions of any of its strong block minimal bases linearizations. Secondly, it relates the right minimal bases of the rational matrix and those of the transfer functions of its strong block minimal bases linearizations.  Lemma \ref{minimalbasis} is based on \cite[Lemma A.1]{Dopico:2018:BlockKronecker}, which is a similar result corresponding to strong block minimal bases pencils of polynomial matrices.

\begin{lemma}\label{minimalbasis}
Let $G(\la)\in\efe(\la)^{p\times m}$ and let $L(\la)$ as in (\ref{eq:ratstrongblockmin}) be a strong block minimal bases linearization of $G(\la)$. Let $\wh{G}(\la)$ be its transfer function matrix, as in (\ref{whG}). Let $N_1(\la)$ be a minimal basis dual to $K_1(\la)$ and let $\wh{N}_2(\la)$ be the matrix in (\ref{eq.unimodSBMBP}).
\begin{itemize}
 \item[\rm{(a)}]  If $h(\la)\in\mathcal{N}_r(G(\la))$ then
 $$z(\la)=\left[\begin{array}{c}N_1(\lambda)^T\\-\wh{N}_2(\la)M(\la)N_1(\la)^T\end{array}\right]h(\la)\in\mathcal{N}_r(\wh{G}(\la)).$$
Moreover, if $0\neq h(\la)\in\mathcal{N}_r(G(\la))$ is a vector polynomial then $z(\la)$ is also a vector polynomial and
\begin{equation}\label{degz}
\deg(z(\la))=\deg(N_1(\la)^T h(\la))=\deg(N_1(\la))+\deg(h(\la)).
\end{equation}
 \item[\rm{(b)}] If $\{h_1(\la),\ldots,h_l(\la)\}$ is a right minimal basis of $G(\la)$ then
 $$\left\{\left[\begin{array}{c}N_1(\lambda)^T\\-\wh{N}_2(\la)M(\la)N_1(\la)^T\end{array}\right]h_1(\la),\ldots,\left[\begin{array}{c}N_1(\lambda)^T\\-\wh{N}_2(\la)M(\la)N_1(\la)^T\end{array}\right]h_l(\la)\right\}$$
 is a right minimal basis of $\wh{G}(\la)$.
% \item[\rm{(c)}] Any right minimal basis of $\wh{G}(\la)$ has the form  $$\left\{\left[\begin{array}{c}N_1(\lambda)^T\\-\wh{N}_2(\la)M(\la)N_1(\la)^T\end{array}\right]h_1(\la),\ldots,\left[\begin{array}{c}N_1(\lambda)^T\\-\wh{N}_2(\la)M(\la)N_1(\la)^T\end{array}\right]h_l(\la)\right\}$$
% where $\{h_1(\la),\ldots,h_l(\la)\}$ is some right minimal basis of $G(\la)$.
\end{itemize}
\end{lemma}

\textbf{Proof}.-
%\begin{enumerate}
%\item
By Proposition \ref{prop_UwhGV}, equation \eqref{eq.ffugvg} and using the structure of $V(\la)$  in \eqref{eq.ffvu} (recall that $X(\la)=\wh{N}_2(\la)M(\la)N_1(\la)^T$) we obtain the first part of (a). Now, we are going to prove (\ref{degz}) following the ideas of
\cite[Lemma A.1]{Dopico:2018:BlockKronecker}.
%It only remains to prove (\ref{eq_degz}).
  It can be seen as in the proof of \cite[Lemma A.1]{Dopico:2018:BlockKronecker} that for any vector polynomial  $g(\la)\neq 0$
 \begin{equation}\label{deg1}
 \deg(N_1(\la)^T g(\la))=\deg(N_1(\la))+\deg(g(\la)),
 %\text{ for any vector polynomial } g(\la)\neq 0
\end{equation}
for any vector polynomial  $y(\la)\neq 0$
 \begin{equation}\label{deg2}
 \deg(K_2(\la)^Ty(\la))=\deg(K_2(\la))+\deg(y(\la))=1+\deg(y(\la)), %\text{ for any vector polynomial } y(\la)\neq 0
 \end{equation}
   and
 \begin{equation}\label{max}
 \deg(z(\la))=\max\{\deg(N_1(\la)^Th(\la)),\deg(X(\la)h(\la))\}.
 \end{equation}
 If $X(\la)h(\la)=0$ then (\ref{degz}) follows. Otherwise, use $0=\wh{G}(\la)z(\la)$ and consider the expression of $\wh{G}(\la)$ in (\ref{whG})
 \begin{align*}
 0&=\left[\begin{array}{cc}M(\la)+\wh{K}_2^TC(\la I_n-A)^{-1}B\wh{K}_1&K_2(\la)^T\\K_1(\la)&0\end{array}\right]\left[\begin{array}{c}N_1(\lambda)^T\\-X(\la)\end{array}\right]h(\la)
 \\
 &=\left[\begin{array}{c}M(\la)N_1(\lambda)^T+\wh{K}_2^TC(\la I_n-A)^{-1}B-K_2(\la)^T X(\la)\\0\end{array}\right]h(\la).
 \end{align*}
 Therefore, $M(\la)N_1(\lambda)^T h(\la)-K_2(\la)^T X(\la)h(\la)=-\wh{K}_2^TC(\la I_n-A)^{-1}Bh(\la)$. Since the expression on the left hand side of this equality is polynomial, the expression on the right hand side must be polynomial. Moreover, by Lemma \ref{coldeg}, $\deg(\wh{K}_2^TC(\la I_n-A)^{-1}Bh(\la))<\deg(h(\la))$ since $\wh{K}_2^TC(\la I_n-A)^{-1}B$ is strictly proper.
 Write the previous expression as
 $$
 K_2(\la)^T X(\la)h(\la)=M(\la)N_1(\lambda)^T h(\la)+\wh{K}_2^TC(\la I_n-A)^{-1}Bh(\la).
 $$
 Notice that (\ref{deg2}) implies that
 \[
 %\begin{array}{ll}
 1+\deg(X(\la)h(\la)) %& =\deg(K_2(\la)^T X(\la)h(\la))\\
% &
=\deg(M(\la)N_1(\lambda)^T h(\la)+\wh{K}_2^TC(\la I_n-A)^{-1}Bh(\la)).
% \end{array}
 \]
 Let us see now that, using the previous expression,
 \begin{equation}\label{eq_deg}
 \deg(X(\la)h(\la))\leq\deg(N_1(\la)^T h(\la)).
 \end{equation}
 If $\deg(\wh{K}_2^TC(\la I_n-A)^{-1}Bh(\la))\leq\deg(M(\la)N_1(\lambda)^T h(\la))$ then
 $$
 1+\deg(X(\la)h(\la))\leq\deg(M(\la)N_1(\lambda)^T h(\la))\leq 1+\deg(N_1(\la)^T h(\la)).
 $$
% Therefore, (\ref{deg1}) and (\ref{max}) prove that $\deg(z(\la))=\deg(N_1(\la))+\deg(h(\la))$.
 Otherwise, if $\deg(\wh{K}_2^TC(\la I_n-A)^{-1}Bh(\la))>\deg(M(\la)N_1(\lambda)^T h(\la))$ then
$$
 1+\deg(X(\la)h(\la))= \deg(\wh{K}_2^TC(\la I_n-A)^{-1}Bh(\la))<\deg(h(\la))\text{ and}
 $$
 $$
 \deg(X(\la)h(\la))<\deg(h(\la))-1<\deg(h(\la))+\deg(N_1(\la))=\deg(N_1(\la)^Th(\la)).
 $$
Therefore, (\ref{deg1}), (\ref{max}) and (\ref{eq_deg}) prove that $\deg(z(\la))=\deg(N_1(\la))+\deg(h(\la))$.
% Thus, $\deg(z(\la))=\deg(N_1(\la))+\deg(h(\la))$.

The proof of part (b) is similar to the proof of \cite[Lemma A.1]{Dopico:2018:BlockKronecker} taking into account that  $\dim \mathcal{N}_r (\wh{G}(\la))=\dim \mathcal{N}_r (G(\la))$.  Therefore, the details are omitted.
\hfill\halmos

\bigskip
As a corollary of  Lemma \ref{minimalbasis} we get the following result on the relationship between the minimal indices of a rational matrix and of the transfer function of any of its strong block minimal bases linearizations.
\begin{corollary}\label{minimalindices}
Let $G(\la)\in\efe(\la)^{p\times m}$ and let $L(\la)$ as in (\ref{eq:ratstrongblockmin}) be a strong block minimal bases linearization of $G(\la)$. Let $\wh{G}(\la)$ be its transfer function matrix, as in (\ref{whG}). Let $N_1(\la)$ be a minimal basis dual to $K_1(\la)$ and $N_2(\la)$ be a minimal basis dual to $K_2(\la)$.
\begin{enumerate}
\item[\rm{(a)}] If $\varepsilon_1\leq\cdots\leq\varepsilon_l$ are the right minimal indices of $G(\la)$ then $\varepsilon_1+\deg(N_1(\la))\leq\cdots\leq\varepsilon_l+\deg(N_1(\la))$ are the right minimal indices of $\wh{G}(\la)$.
\item[\rm{(b)}] If $\eta_1\leq\cdots\leq\eta_q$ are the left minimal indices of $G(\la)$ then $\eta_1+\deg(N_2(\la))\leq\cdots\leq\eta_q+\deg(N_2(\la))$ are the left minimal indices of $\wh{G}(\la)$.
\end{enumerate}
\end{corollary}

\textbf{Proof}.- Part (a) follows from part (b) of Lemma \ref{minimalbasis} and (\ref{degz}). Suppose now that $\eta_1\leq\cdots\leq\eta_q$ are the left minimal indices of $G(\la)$. By Lemma \ref{lem_trans}, $\eta_1\leq\cdots\leq\eta_q$ are the right minimal indices of $G(\la)^T$. Notice that $L(\la)^T$ is a strong block minimal bases linearization of $G(\la)^T$ with transfer function matrix $\wh{G}(\la)^T$. Observe that $S^T$, $A^T$, $T^T$, $B^T$, $C^T$, $\wh{K}_1$, $\wh{K}_2$, $M(\la)^T$, $K_1(\la)$, $K_2(\la)$ in $L(\la)^T$ play the role of  $T$, $A$, $S$, $-C$, $-B$, $\wh{K}_2$, $\wh{K}_1$, $M(\la)$, $K_2(\la)$, $K_1(\la)$ in $L(\la)$ respectively. In particular, $K_2(\la)$ in $L(\la)^T$ plays the role of $K_1(\la)$ in $L(\la)$. Thus, by part (a), $\eta_1+\deg(N_2(\la))\leq\cdots\leq\eta_q+\deg(N_2(\la))$ are the right minimal indices of $\wh{G}(\la)^T$. By Lemma \ref{lem_trans} again, $\eta_1+\deg(N_2(\la))\leq\cdots\leq\eta_q+\deg(N_2(\la))$ are the left minimal indices of $\wh{G}(\la)$.
\hfill\halmos

\bigskip
 Now, we provide a recovery result for the minimal bases of a rational matrix from the minimal bases of the transfer functions of any of its strong block minimal bases linearizations, i.e., the converse of Lemma \ref{minimalbasis}-(b).

\begin{lemma}\label{lem_minimalbasis}
Let $G(\la)\in\efe(\la)^{p\times m}$ and let $L(\la)$ as in (\ref{eq:ratstrongblockmin}) be a strong block minimal bases linearization of $G(\la)$. Let $\wh{G}(\la)$ be its transfer function matrix, as in (\ref{whG}). Let $N_1(\la)$ be a minimal basis dual to $K_1(\la)$, $N_2(\la)$ be a minimal basis dual to $K_2(\la)$ and $\wh{N}_1(\la)$ and $\wh{N}_2(\la)$ be the matrices appearing in (\ref{eq.unimodSBMBP}).
\begin{itemize}
 \item[\rm{(a)}] Any right minimal basis of $\wh{G}(\la)$ has the form  $$\left\{\left[\begin{array}{c}N_1(\lambda)^T\\-\wh{N}_2(\la)M(\la)N_1(\la)^T\end{array}\right]h_1(\la),\ldots,\left[\begin{array}{c}N_1(\lambda)^T\\-\wh{N}_2(\la)M(\la)N_1(\la)^T\end{array}\right]h_l(\la)\right\}$$
 where $\{h_1(\la),\ldots,h_l(\la)\}$ is some right minimal basis of $G(\la)$.

 \item[\rm{(b)}] Any left minimal basis of $\wh{G}(\la)$ has the form  $$\left\{\left[\begin{array}{c}N_2(\lambda)^T\\-\wh{N}_1(\la)M(\la)^TN_2(\la)^T\end{array}\right]j_1(\la),\ldots,\left[\begin{array}{c}N_2(\lambda)^T\\-\wh{N}_1(\la)M(\la)^TN_2(\la)^T\end{array}\right]j_q(\la)\right\}$$
 where $\{j_1(\la),\ldots,j_q(\la)\}$ is some left minimal basis of $G(\la)$.

\end{itemize}
\end{lemma}

\textbf{Proof}.- { The proof is like the one of \cite[Lemma 7.1]{Dopico:2016:BlockKronecker}. Therefore, it is omitted.
\hfill\halmos

\begin{rem} \label{rem.frorecov1}
\rm Lemma \ref{lem_minimalbasis} implies that a right (resp., left) minimal basis of $G(\la)$ can be obtained, or recovered, from any right (resp., left) minimal basis of $\wh{G} (\la)$, as it is described in this remark. Let us focus for brevity only on right minimal bases, since the procedure for left minimal bases is completely analogous. Note first that the vectors $\{\wh{h}_1 (\la), \ldots, \wh{h}_l (\la) \}$ obtained by taking the top $m + \widehat{m}$ entries of the vectors of any right minimal basis of $\wh{G} (\la)$ are always of the form
\begin{equation} \label{eq.hathath}
\{\wh{h}_1 (\la), \ldots, \wh{h}_l (\la) \} = \{N_1(\la)^T h_1 (\la), \ldots, N_1(\la)^T h_l (\la) \},
\end{equation}
with $\{h_1 (\la), \ldots, h_l (\la) \}$ a right minimal basis of $G(\la)$. Then, it is enough to multiply each $\wh{h}_j (\la)$ by a left inverse of $N_1(\la)^T$ in order to get the right minimal basis $\{h_1 (\la), \ldots, h_l (\la) \}$ of $G(\la)$. Such left inverse may be, for instance, the matrix $\wh{K}_1$ in \eqref{eq.unimodSBMBP}. Moreover, in some cases important in applications, the matrices $N_1 (\la)$ and $\wh{K}_1$ are very simple and allow us to recover a right minimal basis of $G(\la)$ without the need of performing any matrix multiplication. This happens, for instance, if $K_1 (\la) = L_\varepsilon (\la) \otimes I_m$ (and $K_2 (\la) = L_\eta (\la) \otimes I_p$) in \eqref{eq:ratstrongblockmin}, where
\begin{equation} \label{eq.bloquekronbasic}
L_k (\la) = \begin{bmatrix}
              -1 & \la &  &  &  \\
                 &  -1 & \la &  &  \\
                 &     &  \ddots & \ddots &  \\
                 &     &         & -1 & \la
            \end{bmatrix} \in \FF [\la]^{k \times (k+1)},
\end{equation}
which corresponds to the well-known block Kronecker linearizations of the polynomial part of $G(\la)$ \cite[Section 4]{Dopico:2018:BlockKronecker} (see also \cite[Examples 5.3 and 5.6]{AmDoMaZa18}). In this case,
$$
N_1 (\la)^T = \begin{bmatrix}
                \la^\varepsilon \\
                \vdots \\
                \la \\
                1
              \end{bmatrix} \otimes I_m \quad \mbox{and} \quad
              \wh{K}_1 = \begin{bmatrix}
                           0 & \cdots & 0 & 1
                         \end{bmatrix} \otimes I_m .
$$
Thus a minimal bases of $G(\la)$ can be obtained just by taking the last $m$ entries of the vectors $\{\wh{h}_1 (\la), \ldots, \wh{h}_l (\la) \}$ in \eqref{eq.hathath}.

\end{rem}

 The next Theorem \ref{th_bases_sbml} is the main result in this section, together with Theorem \ref{th_indices_sbml}, and one of the most relevant results in this paper. Theorem \ref{th_bases_sbml} describes the complete relationship between the minimal bases of a rational matrix and the minimal bases of its strong block minimal bases linearizations in both directions. It follows from combining results in Section \ref{Sec_PSM} with results previously obtained in this section.

\begin{theorem}\label{th_bases_sbml}
Let $G(\la)\in\efe(\la)^{p\times m}$ and let $L(\la)$ as in (\ref{eq:ratstrongblockmin}) be a strong block minimal bases linearization of $G(\la)$. Let $N_1(\la)$ be a minimal basis dual to $K_1(\la)$, $N_2(\la)$ be a minimal basis dual to $K_2(\la)$ and $\wh{N}_1(\la)$ and $\wh{N}_2(\la)$ be the matrices appearing in (\ref{eq.unimodSBMBP}).
\begin{enumerate}
\item[\rm{(a)}]  $\begin{bsmallmatrix}
H_1(\la)\\H_2(\la)\\H_3(\la)
\end{bsmallmatrix}$ is a right minimal basis of $L(\la)$ if and only if
$$\begin{array}{l}H_1(\la)=-S^{-1}(\la I_n-A)^{-1}BH(\la),\\
H_2(\la)=N_1(\la)^TH(\la),\\
H_3(\la)=-\wh{N}_2(\la)M(\la)N_1(\la)^TH(\la)\end{array}
$$
for some  right minimal basis $H(\la)$ of $G(\la)$.

\item[\rm{(b)}]  $\begin{bsmallmatrix}
H_1(\la)\\H_2(\la)\\H_3(\la)
\end{bsmallmatrix}$ is a left minimal basis of $L(\la)$ if and only if
$$\begin{array}{l}H_1(\la)=(C(\la I_n-A)^{-1}T^{-1})^TH(\la),\\
H_2(\la)=N_2(\la)^TH(\la),\\
H_3(\la)=-\wh{N}_1(\la)M(\la)^TN_2(\la)^TH(\la)\end{array}
$$
for some  left minimal basis $H(\la)$ of $G(\la)$.

% \item[\rm{(a)}] If $H(\la)$ is a right minimal basis of $G(\la)$ then
% $$
% \left[\begin{array}{c}-S^{-1}(\la I_n-A)^{-1}B\\N_1(\la)^T\\-\wh{N}_2(\la)M(\la)N_1(\la)^T\end{array}\right]H(\la)
% $$
%is a right minimal basis of $L(\la)$.

% \item[\rm{(b)}] If $J(\la)$ is a left minimal basis of $G(\la)$ then
% $$
% \left[\begin{array}{c}-(C(\la I_n-A)^{-1}T^{-1})^T\\N_2(\la)^T\\-\wh{N}_1(\la)M(\la)^TN_2(\la)^T\end{array}\right]J(\la)
% $$
%is a left minimal basis of $L(\la)$.

%\item[\rm{(c)}]  If $\begin{bmatrix}
%H_1(\la)\\H_2(\la)\\H_3(\la)
%\end{bmatrix}$ is a right minimal basis of $L(\la)$ then
%$$\begin{array}{l}H_1(\la)=-S^{-1}(\la I_n-A)^{-1}BH(\la),\\
%H_2(\la)=N_1(\la)^TH(\la),\\
%H_3(\la)=-\wh{N}_2(\la)M(\la)N_1(\la)^TH(\la)\end{array}
%$$
%where $H(\la)$ is some right minimal basis of $G(\la)$.

%\item[\rm{(d)}]  If $\begin{bmatrix}
%J_1(\la)\\J_2(\la)\\J_3(\la)
%\end{bmatrix}$ is a left minimal basis of $L(\la)$ then
%$$\begin{array}{l}J_1(\la)=-(C(\la I_n-A)^{-1}T^{-1})^TJ(\la),\\
%J_2(\la)=N_2(\la)^TJ(\la),\\
%J_3(\la)=-\wh{N}_1(\la)M(\la)^TN_2(\la)^TJ(\la)\end{array}
%$$
%where $J(\la)$ is some left minimal basis of $G(\la)$.
\end{enumerate}
\end{theorem}

\textbf{Proof}.-
Let $\wh{G}(\la)$ be the transfer function matrix of $L(\la)$. Notice that both $(T(\la I_n-A)S)^{-1}\begin{bmatrix} TB\wh{K}_1 & 0 \end{bmatrix}$ and $\begin{bsmallmatrix} -\wh{K}_2^TCS \\ 0 \end{bsmallmatrix}(T(\la I_n-A)S)^{-1}$
 are strictly proper matrices. By Corollary \ref{col_PGPmin}, $\begin{bsmallmatrix}
H_1(\la)\\H_2(\la)\\H_3(\la)
\end{bsmallmatrix}$ is a right minimal basis of $L(\la)$ if and only if $\begin{bsmallmatrix}
H_2(\la)\\H_3(\la)
\end{bsmallmatrix}$ is a right minimal basis of $\wh{G}(\la)$ and $H_1(\la)=-S^{-1}(\la I_n-A)^{-1}B\wh{K}_1H_2(\la)$. Now, by Lemma \ref{lem_minimalbasis}, $H_2(\la)=N_1(\la)^TH(\la)$ and $H_3(\la)=-\wh{N}_2(\la)M(\la)N_1(\la)^TH(\la)$ for some $H(\la)$ right minimal basis of $G(\la)$.

Part (b) is proved similarly.
\hfill\halmos

\begin{rem} \label{rem.frorecov2}
\rm
Theorem \ref{th_bases_sbml} implies that a right (resp., left) minimal basis of $G(\la)$ can be recovered from any right (resp., left) minimal basis of any of its strong block minimal bases linearizations. Such recovery procedure is completely analogous to the one described in Remark \ref{rem.frorecov1} except for the following minor variation: in the case of Theorem \ref{th_bases_sbml} the right (resp., left) minimal bases of $G(\la)$ have to be recovered from the entries $n+1, n+2, \ldots , n+ m + \widehat{m}$ (resp., $n+1, n+2, \ldots , n+ p + \widehat{p}$) of the vectors of the right (resp., left) minimal bases of its strong block minimal bases linearizations. As in Remark \ref{rem.frorecov1}, the recovery is extremely simple for strong block minimal bases linearizations of $G(\la)$ constructed from a block Kronecker linearization of its polynomial part.
\end{rem}

In the last result of this section, the relationship between the minimal indices of a rational matrix and those of its strong block minimal bases linearizations is established.
\begin{theorem}\label{th_indices_sbml}
Let $G(\la)\in\efe(\la)^{p\times m}$ and let $L(\la)$ as in (\ref{eq:ratstrongblockmin}) be a strong block minimal bases linearization of $G(\la)$. Let $N_1(\la)$ be a minimal basis dual to $K_1(\la)$ and $N_2(\la)$ be a minimal basis dual to $K_2(\la)$.
\begin{enumerate}
\item[\rm{(a)}] If $\varepsilon_1\leq\cdots\leq\varepsilon_l$ are the right minimal indices of $G(\la)$ then $\varepsilon_1+\deg(N_1(\la))\leq\cdots\leq\varepsilon_l+\deg(N_1(\la))$ are the right minimal indices of $L(\la)$.
\item[\rm{(b)}] If $\eta_1\leq\cdots\leq\eta_q$ are the left minimal indices of $G(\la)$ then $\eta_1+\deg(N_2(\la))\leq\cdots\leq\eta_q+\deg(N_2(\la))$ are the left minimal indices of $L(\la)$.
\end{enumerate}
\end{theorem}

\textbf{Proof}.- Let $\wh{G}(\la)$ be the transfer function matrix of $L(\la)$. If $\varepsilon_1\leq\cdots\leq\varepsilon_l$ are the right minimal indices of $G(\la)$ then, by Corollary \ref{minimalindices}, $\varepsilon_1+\deg(N_1(\la))\leq\cdots\leq\varepsilon_l+\deg(N_1(\la))$ are the right minimal indices of $\wh{G}(\la)$. Now, by Theorem \ref{minimalbasisPG}, these are the right minimal indices of $L(\la)$.

A similar proof can be done in order to prove (b).
\hfill\halmos

%---------------------------------------------------------------------------------------------
\section{Minimal bases and indices of $\M_{1}$ and $\M_{2}$-strong linearizations of rational matrices}\label{Sec_M12minimalbases}
%-------------------------------------------------------------------------------------------

$\M_{1}$ and $\M_{2}$-strong linearizations of square rational matrices have been recently introduced in \cite{DoMaQu19} by combining results from \cite{AmDoMaZa18} with the $\M_{1}$ and $\M_{2}$ ansatz spaces of linearizations of a polynomial matrix developed in \cite{FasSal17}, which in turn are inspired by the pioneer $\mathbb{L}_{1}$ and $\mathbb{L}_{2}$ vector spaces of linearizations of matrix polynomials introduced in \cite{MMMM06}. Among other properties, $\M_{1}$ and $\M_{2}$-strong linearizations of rational matrices allow us to deal very easily with rational matrices whose polynomial part is expressed in any orthogonal basis. In this section, we study the minimal bases and indices of $\M_{1}$ and $\M_{2}$-strong linearizations of rational matrices. Since these families of linearizations are closely connected to strong block minimal bases linearizations, it is not surprising that the results of this section are easily obtained from combining those in Section \ref{Sec_blockminimalbases} with  specific properties of $\M_{1}$ and $\M_{2}$-strong linearizations. {\color{black} Besides, we extend the results on minimal indices and bases to the strong linearizations of rational matrices whose polynomial part is expressed in other polynomial bases, not only orthogonal, that were developed in \cite[Section 9]{DoMaQu19}.} In order to proceed, we need to recap first some results and notations taken from \cite{DoMaQu19}.

The following lemma establishes a general result about the relationship  between the minimal bases and indices of two rational matrices connected by a nonsingular constant  matrix on the left. We will see that this simple result will
allow us to obtain the relationship between the minimal bases and indices of a rational matrix and its $\M_1$, $\M_2$-strong linearizations. The reason is that an  $\M_1$-strong linearization is a strong block minimal bases linearization premultiplied by a nonsingular constant matrix, and an $\M_2$-strong linearization is a strong block minimal bases linearization postmultiplied by a nonsingular constant matrix.

\begin{lemma}\label{lem_premult}
	Let $G_1(\la), G_2(\la)\in\efe(\la)^{p\times m}$ and $X\in\efe^{p\times p}$ be nonsingular such that $G_2(\la)=XG_1(\la)$. Then, $H(\la)$ is a right minimal basis of $G_1(\la)$ if and only if $H(\la)$ is a right minimal basis of $G_2(\la)$ and $\overline{H}(\la)$ is a left minimal basis of $G_1(\la)$ if and only if $X^{-T}\overline{H}(\la)$ is a left minimal basis of $G_2(\la)$. Moreover, $G_1(\la)$ and $G_2(\la)$ have the same right minimal indices and the same left minimal indices.
\end{lemma}

\textbf{Proof}.-
Notice that $G_1(\la)H(\la)=0$ if and only if $G_2(\la)H(\la)=0$. Moreover, by Lemma \ref{lem_trans}, $\overline{H}(\la)$ is a left minimal basis of $G_1(\lambda)$ if and only if $\overline{H}(\la)$ is a right minimal basis of $G_1(\lambda)^T$. Furthermore, $G_1(\lambda)^T\overline{H}(\la)=0$ if and only if $G_2(\lambda)^TX^{-T}\overline{H}(\la)=0$ and, by \cite[Lemma 2.16]{DeDoVa15}, $X^{-T}\overline{H}(\la)$ is a minimal basis with the same column degrees as $\overline{H}(\la)$. Therefore, $\overline{H}(\la)$ is a right minimal basis of $G_1(\la)^T$ if and only if $X^{-T}\overline{H}(\la)$ is a right minimal basis of $G_2(\lambda)^T$ and, by Lemma \ref{lem_trans} again, $X^{-T}\overline{H}(\la)$ is a left minimal basis of $G_2(\lambda)$.
\hfill\halmos

\bigskip
 The definitions of the $\M_1$ and  $\M_2$-strong linearizations introduced in Subsections \ref{subsec.7.1} and \ref{subsec.7.2} are based on the matrices and vectors presented in the next paragraphs. Consider a polynomial basis $\{\phi_{j}(\lambda)\}_{j=0}^{\infty}$ of $\efe[\lambda],$ viewed as an $\efe$-vector space, with $\phi_{j}(\lambda)$ a polynomial of degree $j,$  that satisfies the following three-term recurrence relation:

\begin{equation}\label{recu}
\alpha_{j}\phi_{j+1}(\lambda)=(\lambda-\beta_{j})\phi_{j}(\lambda)-\gamma_{j}\phi_{j-1}(\lambda) \quad j\geq 0
\end{equation}
where $\alpha_{j},\beta_{j},\gamma_{j}\in\efe,$ $ \alpha_{j}\neq 0,$ $\phi_{-1}(\lambda)=0,$ and $\phi_{0}(\lambda)=1.$

Let $G(\lambda)\in \efe(\lambda)^{m\times m}$ be a rational matrix, let $G(\lambda)=D(\lambda)+G_{sp}(\lambda)$ be its unique decomposition into its polynomial part $D(\lambda)\in\efe[\lambda]^{m\times m}$ and its strictly proper part $G_{sp}(\lambda)\in\efe_{pr}(\lambda)^{m\times m},$ and let $G_{sp}(\lambda)=C(\lambda I_{n}-A)^{-1}B$ be a minimal order state-space realization of $G_{sp}(\lambda),$ where $n=\nu(G(\lambda)).$ Assume that $\deg(D(\lambda))\geq 2.$ Write $D(\lambda)$ in terms of the polynomial basis $\{\phi_{j}(\lambda)\}_{j=0}^{\infty},$ as
	\begin{equation}\label{eq_D}
	D(\lambda)=D_{k}\phi_{k}(\lambda)+D_{k-1}\phi_{k-1}(\lambda)+\cdots + D_{1}\phi_{1}(\lambda)+D_{0}\phi_{0}(\lambda)
	\end{equation}
	with $D_{k}\neq 0$. Let
	\begin{equation}\label{eq_fik}
		\Phi_{k}(\lambda)=[\phi_{k-1}(\lambda)\cdots \phi_{1}(\lambda)\text{ }\phi_{0}(\lambda)]^{T},
	\end{equation}
	$$m_{\Phi}^{D}(\lambda)=\left[ \dfrac{(\lambda-\beta_{k-1})}{\alpha_{k-1}}D_{k}+D_{k-1}\quad D_{k-2}-\dfrac{\gamma_{k-1}}{\alpha_{k-1}}D_{k}\quad D_{k-3}\quad \cdots\quad D_{1}\quad D_{0}\right],
	$$
	$$M_{\Phi}(\lambda)=\left[ {\begin{array}{cccccc}
		-\alpha_{k-2} & (\lambda -\beta_{k-2}) & -\gamma_{k-2} &  \\
		& -\alpha_{k-3} & (\lambda - \beta_{k-3}) & -\gamma_{k-3} &  \\
		& &\ddots &\ddots & \ddots &  \\
		& & &-\alpha_{1}& (\lambda-\beta_{1}) & -\gamma_{1} \\
		& & & &-\alpha_{0}& (\lambda-\beta_{0})
		\end{array} } \right],
	$$
	\begin{equation}\label{eq_F}
	F_{\Phi}^{D}(\lambda)=\left[ {\begin{array}{cc}
		m_{\Phi}^{D}(\lambda) \\
		M_{\Phi}(\lambda)\otimes I_{m} \\
		\end{array} } \right].
	\end{equation}
	
\subsection{Minimal bases and indices of $\M_1$-strong linearizations of rational matrices} \label{subsec.7.1}
	We investigate first the relationship between the minimal bases and indices of a rational matrix and its $\M_1$-strong linearizations.
	
	For any nonsingular constant matrices $T,S\in\efe^{n\times n}$ the linear polynomial matrix
	\begin{equation}\label{eq_linsb}
	\begin{array}{ll}L(\lambda)& = \left[
	\begin{array}{c|c}
	
		T(\lambda I_{n}-A)S& 0_{n\times (k-1)m}\quad TB\\
		\hline
		
		-CS&m_{\Phi}^{D}(\lambda) \\
		0_{(k-1)m\times n} & M_{\Phi}(\lambda)\otimes I_{m}
		\end{array}
		\right]
		\\& =
		 \left[
		\begin{array}{c|c}
		
		T(\lambda I_{n}-A)S& 0_{n\times (k-1)m}\quad TB\\
		\hline
		
		\begin{array}{c}
		-CS\\
		0_{(k-1)m\times n}
		\end{array}& F_{\Phi}^{D}(\lambda)
		\end{array}\right]
\end{array}
		\end{equation}
	is a strong linearization of $G(\lambda)$ (see \cite[Theorem 3.8]{DoMaQu19}).  Taking into account also \cite[Lemmas 3.3 and 3.7]{DoMaQu19} and their proofs, note that $L(\la)$ is a strong block minimal bases linearization of $G(\lambda)$ as in (\ref{eq:ratstrongblockmin}) with $M(\la)=m_{\Phi}^{D}(\lambda)$, $K_1(\la)=M_{\Phi}(\lambda)\otimes I_m$, $K_2(\la)$ empty, $N_1(\la)=(\Phi_{k}(\lambda)\otimes I_m)^T=\Phi_{k}(\lambda)^T\otimes I_m$, $N_2(\la)=I_m$, $\wh{K}_1=e_k^T\otimes I_m$, $\wh{K}_2=I_m$ and $\wh{N}_2(\la)$ empty, where $e_k$ is the $k$th canonical vector of size $k\times 1$. Moreover, $\wh{N}_1(\la)$ is of the form $\wh{N}_1(\la)=Q(\la)^T\otimes I_m$ with $Q(\la)= \begin{bsmallmatrix}M_{\Phi}(\la)\\e_k^T\end{bsmallmatrix}^{-1}\begin{bsmallmatrix}
	I_{k-1}\\0\end{bsmallmatrix}$.
	
	Furthermore, let $v\in\efe^{k},$ $J\in\efe^{km\times (k-1)m}$ with $[v\otimes I_{m}\quad J]$ nonsingular and let $\mathcal{L}(\lambda)=[v\otimes I_{m}\quad J]F_{\Phi}^{D}(\lambda).$
	Then, the linear polynomial matrix
		\begin{align}\label{eq_linM1}
		    L_1(\lambda) & = \left[\begin{array}{c|c}
			I_{n} & 0 \\\hline
			0 & v\otimes I_{m}\quad J
			\end{array}\right]\left[
			\begin{array}{c|c}	
			T(\lambda I_{n}-A)S& 0_{n\times (k-1)m}\quad TB\\
			\hline			
%			\begin{array}{c}
-CS&m_{\Phi}^{D}(\lambda) \\
		0_{(k-1)m\times n} & M_{\Phi}(\lambda)\otimes I_{m}
%			-CS\\
%			0_{(k-1)m\times n}
%			\end{array}& F_{\Phi}^{D}(\lambda)
			\end{array}
			\right] \nonumber \\ & =
			\left[
		\begin{array}{c|c}
		T(\lambda I_{n}-A)S& 0_{n\times (k-1)m}\quad TB\\
		\hline
		-(v\otimes I_{m})CS& \mathcal{L}(\lambda)
		\end{array}
		\right]
		\end{align}
		is a strong linearization of $G(\lambda)$, which is called $\M_1$-strong linearization of $G(\la)$ (see \cite[Theorem 3.9]{DoMaQu19}). Put $X=\left[\begin{array}{c|c}I_{n} & 0 \\\hline
		0 & v\otimes I_{m}\quad J
		\end{array}\right]$, which is nonsingular. Thus, $L_1(\lambda)=X L(\lambda)$.

With all these results at hand, Theorem \ref{th.m1recovery} establishes the relationships between the minimal bases and indices of a rational matrix and its $\M_1$-strong linearizations.

%%%%%%%%%%%%%%%%%%%%%%%%%%%%%%%%%%%

\begin{theorem} \label{th.m1recovery}
Let $G(\la)\in\efe(\la)^{m\times m}$ and let $L_1(\la)$ as in (\ref{eq_linM1}) be an $\M_{1}$-strong linearization of $G(\la)$. Let $\Phi_k(\la)$ be as in (\ref{eq_fik}).

\begin{enumerate}
\item[\rm{(a)}]  $\begin{bsmallmatrix}
H_1(\la)\\H_2(\la)\\\vdots\\H_{k+1}(\la)
\end{bsmallmatrix}$ is a right minimal basis of $L_1(\la)$ if and only if  $H_{k+1}(\la)$ is a right minimal basis of $G(\la)$ and
$$\begin{array}{l}H_1(\la)=-S^{-1}(\la I_n-A)^{-1}BH_{k+1}(\la),\\
H_i(\la)=\phi_{k-i+1}(\la)H_{k+1}(\la),\  i=2,\ldots ,k.
\end{array}
$$

\item[\rm{(b)}]  If $\begin{bsmallmatrix}
H_1(\la)\\H_2(\la)
\end{bsmallmatrix}$ is a left minimal basis of $L_1(\la)$ then $(v^T\otimes I_m)H_2(\la)$ is a  left minimal basis of $G(\la)$ and $H_1(\la)=(C(\la I_n-A)^{-1}T^{-1})^T(v^T\otimes I_m)H_2(\la)$.

\item[\rm{(c)}] If $H(\la)$ is a left minimal basis  of $G(\la)$ then $\begin{bsmallmatrix}
H_1(\la)\\H_2(\la)
\end{bsmallmatrix}$ is a left minimal basis of $L_1(\la)$ where
$$\begin{array}{l}H_1(\la)=(C(\la I_n-A)^{-1}T^{-1})^TH(\la),\\
H_2(\la)=
\begin{bmatrix}v\otimes I_m & J\end{bmatrix}^{-T}\left[\begin{array}{c}H(\la)\\
-\wh{N}_1(\la)m_{\Phi}^{D}(\lambda)^TH(\la)\end{array}\right]\end{array}
$$
with $\wh{N}_1(\la)=Q(\la)^T\otimes I_m$ such that $Q(\la)= \begin{bsmallmatrix}M_{\Phi}(\la)\\e_k^T\end{bsmallmatrix}^{-1} \begin{bsmallmatrix}
	I_{k-1}\\0\end{bsmallmatrix}$.

\item[\rm{(d)}] If $\varepsilon_1\leq\cdots\leq\varepsilon_l$ are the right minimal indices of $G(\la)$ then $\varepsilon_1+k-1\leq\cdots\leq\varepsilon_l+k-1$ are the right minimal indices of $L_1(\la)$.

\item[\rm{(e)}] If $\eta_1\leq\cdots\leq\eta_l$ are the left minimal indices of $G(\la)$ then $\eta_1\leq\cdots\leq\eta_l$ are the left minimal indices of $L_1(\la)$.
\end{enumerate}
\end{theorem}

\textbf{Proof}.-
To prove (a), by using Lemma \ref{lem_premult}, we get that $\begin{bsmallmatrix}
H_1(\la)\\H_2(\la)\\\vdots\\H_{k+1}(\la)
\end{bsmallmatrix}$ is a right minimal basis of $L_1(\la)$ if and only if it is a right minimal basis of $L(\la)$ in (\ref{eq_linsb}). By the fact that $L(\la)$ is a strong block minimal basis linearization of $G(\la)$ and Theorem \ref{th_bases_sbml} (with $\wh{N}_2 (\la)$ empty), this occurs if and only if
 $H_1(\la)=-S^{-1}(\la I_n-A)^{-1}BH(\la)$ and $\begin{bsmallmatrix}
H_2(\la)\\\vdots\\H_{k+1}(\la)
\end{bsmallmatrix}=(\Phi_k(\la)\otimes I_m)H(\la)$ for some right minimal basis $H(\la)$ of $G(\la)$. But, since $\phi_0(\la)=1$, $H(\la)$ is uniquely determined as $H(\la)=H_{k+1}(\la)$, and
$H_i(\la)=\phi_{k-i+1}(\la)H_{k+1}(\la),\  i=2,\ldots ,k$.

The proof of the other parts can be done similarly by using Lemma \ref{lem_premult} and Theorems \ref{th_bases_sbml} or \ref{th_indices_sbml}. Observe that in this case $G(\la)$ is square and, therefore, it has a number of left minimal indices equal to the number of right minimal indices.
\hfill\halmos

\begin{rem} \label{rem.recoveryfro3}
 \rm
Part (a) of Theorem \ref{th.m1recovery}, together with the fact that $\phi_0 (\la) = 1$, provides a very simple recovery rule of a right minimal basis of $G(\la)$ from any right minimal basis of any of its $\M_{1}$-strong linearizations: simply take the last $m$ rows of the right minimal basis of the $\M_{1}$-strong linearization. Part (b) of Theorem \ref{th.m1recovery} also provides a simple recovery rule of a left minimal basis of $G(\la)$ from any left minimal basis of any of its $\M_{1}$-strong linearizations, though in this case some arithmetic operations are required unless $v$ is one of the canonical vectors of $\FF^k$.
\end{rem}

{
\subsection{Minimal bases and indices of extended $\mathbb{M}_1$-strong linearizations of rational matrices} \label{subsec.newextendedm1}
}

%\begin{rem}
 %\rm
{ Similar ideas and techniques to those in Subsection \ref{subsec.7.1} can be used to obtain a result similar to Theorem \ref{th.m1recovery} for the strong linearizations of rational matrices with polynomial part expressed in any degree-graded polynomial basis, i.e., one whose $j$th element has degree $j$, introduced in \cite[Section 9]{DoMaQu19}. We call such linearizations extended $\mathbb{M}_1$-strong linearizations. More precisely, suppose as in \cite[Section 9]{DoMaQu19} that $G(\la)=D(\la)+C(\la I_n-A)^{-1}B$ with $\deg(D(\la))=k\geq2$ and consider a polynomial basis $\{\psi_{j}(\lambda)\}_{j=0}^{\infty}$ of $\efe[\lambda],$ with $\psi_{j}(\lambda)$ a polynomial of degree $j,$  that satisfies the linear relation
$M_{\Psi}(\la) \Psi_{k}(\lambda) =0$,
where $M_{\Psi}(\la)\in\efe[\la]^{(k-1)\times k}$ is a minimal basis with all its row degrees equal to $1,$ and
\begin{equation}\label{eq_psik}\Psi_{k}(\lambda)=[\psi_{k-1}(\lambda)\;\cdots\; \psi_{1}(\lambda)\text{ }\psi_{0}(\lambda)]^{T}.\end{equation}
Note that  $\Psi_{k}(\lambda_0)\neq 0$ for all $\la_0\in\overline{\efe}$ because $\psi_{0}(\lambda)$ has degree $0$. Example 9.2 in \cite{DoMaQu19} shows how to construct $M_{\Psi}(\la)$ from a recurrence relation that holds under the assumption that the polynomials in the basis are monic. Let $m_{\Psi}^{D}(\la)$ be a pencil such that $m_{\Psi}^{D}(\la)(\Psi_{k}(\lambda)\otimes I_{m})=D(\la)$ (see again \cite[Example 9.2]{DoMaQu19} for an explicit construction). By \cite[Lemma 5.5]{AmDoMaZa18}, there exists a vector $w\in\efe^{k}$ such that
%\begin{equation}\label{unimodular}
%U(\lambda)=\left[ {\begin{array}{cc}
%	M_{\Psi}(\lambda)\otimes I_{m} \\
%	w^T\otimes I_{m}
%	\end{array} } \right]	
%\end{equation}
$U(\lambda)=\begin{bsmallmatrix}
	M_{\Psi}(\lambda)\\
	w^T
	\end{bsmallmatrix}$
is unimodular, and its inverse has the form $U(\lambda)^{-1}=[R(\lambda)\quad \Psi_{k}(\lambda)]$
with $R(\lambda)= \begin{bsmallmatrix}M_{\Psi}(\la)\\w^T\end{bsmallmatrix}^{-1}\begin{bsmallmatrix}
	I_{k-1}\\0\end{bsmallmatrix}$. { Note that in the construction of \cite[Example 9.2]{DoMaQu19}, one can take simply $w^T = e_k^T = [0\, \cdots \, 0 \, 1]^T$.}
Let $v\in\efe^{k}$, $[v\otimes I_{m}\quad J]$  nonsingular, and  $\mathcal{L}(\la)=[v\otimes I_{m}\quad J]\left[\begin{array}{c}
m_{\Psi}^{D}(\la)\\
M_{\Psi}(\la)\otimes I_{m}
\end{array}\right]
$. By \cite[Theorem 9.1]{DoMaQu19}, for any $T$ and $S$ nonsingular,
\begin{equation}\label{eq_linotherbases}
\begin{array}{ll}
L_1(\lambda)&=
\left[\begin{array}{c|c}
			I_{n} & 0 \\\hline
			0 & v\otimes I_{m}\quad J
			\end{array}\right]\left[
			\begin{array}{c|c}	
			T(\lambda I_{n}-A)S& TB(w^T\otimes I_m)\\
			\hline			
-CS&m_{\Psi}^{D}(\lambda) \\
		0_{(k-1)m\times n} & M_{\Psi}(\lambda)\otimes I_{m}
			\end{array}
			\right]\\&
	= \left[
	\begin{array}{c|c}
	
	T(\lambda I_{n}-A)S&  TB(w^T\otimes I_m)\\
	\hline \phantom{\Big|}
	
	\begin{array}{c}
	-(v\otimes I_m)CS\\
	\end{array}&\mathcal{L}(\lambda)
	\end{array}
	\right]
	\end{array}
\end{equation}
	is a strong linearization of $G(\lambda)$, { that we call extended $\mathbb{M}_1$-strong linearization}. Moreover, the matrix on the right of the first equality is a strong block minimal bases linearization of $G(\la)$ as in (\ref{eq:ratstrongblockmin}) with $M(\la)=m_{\Psi}^{D}(\lambda)$, $K_1(\la)=M_{\Psi}(\lambda)\otimes I_m$, $K_2(\la)$ empty, $\wh{K}_1=w^T\otimes I_m$, and $\wh{K}_2=I_m$. Furthermore, $N_1(\la)=\Psi_{k}(\lambda)^T\otimes I_m$, $N_2(\la)=I_m$, $\wh{N}_1(\la)=R(\la)^T\otimes I_m$, and $\wh{N}_2(\la)$ is empty. With all this in mind, by using the same techniques as for the proof of Theorem \ref{th.m1recovery}, the following result is obtained.

\smallskip

\begin{theorem} \label{thm.degreegradedbasis}
Let $G(\la)\in\efe(\la)^{m\times m}$ and let $L_1(\la)$ as in (\ref{eq_linotherbases}) be an { extended $\mathbb{M}_1$-strong linearization} of $G(\la)$. Let $\Psi_k(\la)$ be as in (\ref{eq_psik}).

\begin{enumerate}
\item[\rm{(a)}]  $\begin{bsmallmatrix}
H_1(\la)\\H_2(\la)\\\vdots\\H_{k+1}(\la)
\end{bsmallmatrix}$ is a right minimal basis of $L_1(\la)$ if and only if
$$\begin{array}{l}H_1(\la)=-S^{-1}(\la I_n-A)^{-1}BH(\la),\\
H_i(\la)=\psi_{k-i+1}(\la)H(\la),\  i=2,\ldots ,k+1.
\end{array}
$$
for some right minimal basis $H(\la)$ of $G(\la)$.

\item[\rm{(b)}]  If $\begin{bsmallmatrix}
H_1(\la)\\H_2(\la)
\end{bsmallmatrix}$ is a left minimal basis of $L_1(\la)$ then $(v^T\otimes I_m)H_2(\la)$ is a  left minimal basis of $G(\la)$ and $H_1(\la)=(C(\la I_n-A)^{-1}T^{-1})^T(v^T\otimes I_m)H_2(\la)$.

\item[\rm{(c)}] If $H(\la)$ is a left minimal basis  of $G(\la)$ then $\begin{bsmallmatrix}
H_1(\la)\\H_2(\la)
\end{bsmallmatrix}$ is a left minimal basis of $L_1(\la)$ where
$$\begin{array}{l}H_1(\la)=(C(\la I_n-A)^{-1}T^{-1})^TH(\la),\\
H_2(\la)=
\begin{bmatrix}v\otimes I_m & J\end{bmatrix}^{-T}\left[\begin{array}{c}H(\la)\\
-\wh{N}_1(\la)m_{\Psi}^{D}(\lambda)^TH(\la)\end{array}\right]\end{array}
$$
with $\wh{N}_1(\la)=R(\la)^T\otimes I_m$ such that $R(\la)= \begin{bsmallmatrix}M_{\psi}(\la)\\w^T\end{bsmallmatrix}^{-1} \begin{bsmallmatrix}
	I_{k-1}\\0\end{bsmallmatrix}$.

\item[\rm{(d)}] If $\varepsilon_1\leq\cdots\leq\varepsilon_l$ are the right minimal indices of $G(\la)$ then $\varepsilon_1+k-1\leq\cdots\leq\varepsilon_l+k-1$ are the right minimal indices of $L_1(\la)$.

\item[\rm{(e)}] If $\eta_1\leq\cdots\leq\eta_l$ are the left minimal indices of $G(\la)$ then $\eta_1\leq\cdots\leq\eta_l$ are the left minimal indices of $L_1(\la)$.
\end{enumerate}
\end{theorem}

{ Note that in Theorem \ref{thm.degreegradedbasis}(a), $H_{k+1} (\la)$ itself is a right minimal basis of $G(\la)$ since $\psi_0 (\la)$ has degree zero and is a constant.}

%\end{rem}
}

%%%%%%%%%%%%%%%%%%%%%%%%%%%%%%%%%%%%%%%%%%%%%%%%%

\subsection{Minimal bases and indices of $\M_2$-strong linearizations of rational matrices} \label{subsec.7.2}
%Consider a polynomial basis $\{\phi_{j}(\lambda)\}_{j=0}^{\infty}$ defined as in (\ref{recu}). Let $G(\lambda)\in \efe(\lambda)^{m\times m}$ be a rational matrix, let $G(\lambda)=D(\lambda)+G_{sp}(\lambda)$ be its unique decomposition into its polynomial part $D(\lambda)\in\efe[\lambda]^{m\times m}$ and its stricly proper part $G_{sp}(\lambda)\in\efe_{pr}(\lambda)^{m\times m},$ and let $G_{sp}(\lambda)=C(\lambda I_{n}-A)^{-1}B$ be a minimal order state-space realization of $G_{sp}(\lambda),$ where $n=\nu(G(\lambda))=\nu(G_{sp}(\lambda)).$ Assume that $\deg(D(\lambda))=k\geq 2.$ Write $D(\lambda)$ in terms of the polynomial basis $\{\phi_{k}(\lambda)\}_{k=0}^{\infty}$ as in (\ref{eq_D}).

We investigate now the relationship between the minimal bases and indices of a rational matrix and its $\M_2$-strong linearizations. The developments and results are very similar to those in Subsection \ref{subsec.7.1} and, therefore, are described briefly.

Let $Q(\la)$ be a $km\times lm$ pencil of the form
$Q(\lambda)=\displaystyle\sum_{i=1}^{k}\displaystyle\sum_{j=1}^{l}e_{i}e_{j}^{T}\otimes Q_{ij}(\lambda)$ for certain $m\times m$  pencils $Q_{ij}(\lambda)$, and where $e_i$ (resp., $e_j$) is the $i$th (resp., $j$th) canonical vector in $\FF^k$ (resp., $\FF^l$). The $lm \times km$ pencil $Q(\lambda)^{\mathcal{B}}=\displaystyle\sum_{i=1}^{k}\displaystyle\sum_{j=1}^{l}e_{j}e_{i}^{T}\otimes Q_{ij}(\lambda)$ is the block-transpose of $Q(\lambda).$
Notice that the block-transpose of $F_{\Phi}^{D}(\lambda)$ in (\ref{eq_F}) is $F_{\Phi}^{D}(\lambda)^{\mathcal{B}}=[
	{m_{\Phi}^{D}(\lambda)}^{\mathcal{B}} \quad
	{M_{\Phi}(\lambda)}^{T}\otimes I_{m} ].$

	For any nonsingular constant matrices $T,S\in\efe^{n\times n}$ the linear polynomial matrix
	\begin{equation}\label{eq_linsb2}
	\begin{array}{ll}\mathbb{L}(\lambda)&= \left[
	\begin{array}{c|c}
	
		T(\lambda I_{n}-A)S&TB\qquad 0_{n\times (k-1)m}\\
		\hline
		\begin{array}{c}
		0_{(k-1)m\times n}\\
		-CS\end{array}& m_{\Phi}^D(\la)^{\mathcal{B}}\quad M_{\Phi}(\lambda)^T\otimes I_{m}	
		\end{array}
		\right]\\[0.7cm] & =
		 \left[
		\begin{array}{c|c}
		
		T(\lambda I_{n}-A)S&  TB \quad 0_{n\times (k-1)m} \\
		\hline
		\begin{array}{c}
		0_{(k-1)m\times n}\\
		-CS\end{array}& F_{\Phi}^{D}(\lambda)^{\mathcal{B}}
		\end{array}\right]\end{array}
		\end{equation}
	is a strong linearization of $G(\lambda)$ (see \cite[Theorem 4.3]{DoMaQu19}). Notice that $\mathbb{L}(\la)$ is a strong block minimal bases linearization of $G(\lambda)$ as in (\ref{eq:ratstrongblockmin}) with $M(\la)=m_{\Phi}^{D}(\lambda)^{\mathcal{B}}$, $K_1(\la)$ empty, $K_2(\la)=M_{\Phi}(\lambda)\otimes I_m$, $N_1(\la)=I_m$, $N_2(\la)=(\Phi_{k}(\lambda)\otimes I_m)^T=\Phi_{k}(\lambda)^T\otimes I_m$, $\wh{K}_1=I_m$, $\wh{K}_2=e_k^T\otimes I_m$ and $\wh{N}_1(\la)$ empty. Moreover, $\wh{N}_2(\la)=Q(\la)^T\otimes I_m$ such that $Q(\la)= \begin{bsmallmatrix}M_{\Phi}(\la)\\e_k^T\end{bsmallmatrix}^{-1}\begin{bsmallmatrix}
	I_{k-1}\\0\end{bsmallmatrix}$.
	
Furthermore, let $w\in\efe^{k},$ $J\in\efe^{km\times (k-1)m}$ with
$\begin{bsmallmatrix}
w^T\otimes I_{m}\\ J^{\mathcal{B}}
\end{bsmallmatrix}$
nonsingular and
$\mathcal{L}(\lambda)=F_{\Phi}^{D}(\lambda)^{\mathcal{B}}\begin{bsmallmatrix}
w^T\otimes I_{m}\\ J^{\mathcal{B}}
\end{bsmallmatrix}$.
Then, the linear polynomial matrix
\begin{align} \nonumber
		\mathbb{L}_2(\lambda) &= \left[
			\begin{array}{c|c}			
				T(\lambda I_{n}-A)S&TB\qquad 0_{n\times (k-1)m}\\
				\hline
				\begin{array}{c}
				0_{(k-1)m\times n}\\
				-CS\end{array}& m_{\Phi}^D(\la)^{\mathcal{B}}\quad M_{\Phi}(\lambda)^T\otimes I_{m}	
				\end{array}
				\right]\left[\begin{array}{c|c}
						I_{n} & 0 \\\hline
						0 & w^T\otimes I_{m}\\ 0&J^{\mathcal{B}}
						\end{array}\right]\\ \label{eq_linM2}
		&=	\left[
		\begin{array}{c|c}		
		T(\lambda I_{n}-A)S& TB(w^T\otimes I_{m})\\
		\hline
		\begin{array}{c}
					0_{(k-1)m\times n}\\
					-CS\end{array}& \mathcal{L}(\lambda)
		\end{array}
		\right]
\end{align}
		is a strong linearization of $G(\lambda)$, which is called $\M_2$-strong linearization of $G(\la)$ (see \cite[Theorem 4.4]{DoMaQu19}). Put $Y=\left[\begin{array}{c|c}
								I_{n} & 0 \\\hline
								0 & w^T\otimes I_{m}\\ 0&J^{\mathcal{B}}
								\end{array}\right]$, which is nonsingular. Thus, $\mathbb{L}_2(\lambda)=\mathbb{L}(\lambda)Y$.

The relationship between the minimal bases and indices of a rational matrix and its $\M_2$-strong linearizations is given in Theorem \ref{thm.M2recovery}.

\begin{theorem} \label{thm.M2recovery}
Let $G(\la)\in\efe(\la)^{m\times m}$ and let $\mathbb{L}_2(\la)$ as in (\ref{eq_linM2}) be an $\M_{2}$-strong linearization of $G(\la)$. Let $\Phi_k(\la)$ be as in (\ref{eq_fik}).

\begin{enumerate}
\item[\rm{(a)}]  If $\begin{bsmallmatrix}
H_1(\la)\\H_2(\la)
\end{bsmallmatrix}$ is a right minimal basis of $\mathbb{L}_2(\la)$ then
$(w^T\otimes I_m)H_2(\la)$ is a right minimal basis of $G(\la)$ and
$H_1(\la)=-S^{-1}(\la I_n-A)^{-1}B(w^T\otimes I_m)H_2(\la)$.

\item[\rm{(b)}] If $H(\la)$  is a right minimal basis of $G(\la)$ then
 $\begin{bsmallmatrix}
H_1(\la)\\H_2(\la)
\end{bsmallmatrix}$ is a right minimal basis of $\mathbb{L}_2(\la)$ where
$$\begin{array}{l}H_1(\la)=-S^{-1}(\la I_n-A)^{-1}BH(\la),\\
H_2(\la)=
\begin{bmatrix}w^T\otimes I_m \\ J^{\mathcal{B}}\end{bmatrix}^{-1}\left[\begin{array}{c}H(\la)\\
-\wh{N}_2(\la)m_{\Phi}^{D}(\lambda)^{\mathcal{B}}H(\la)\end{array}\right]\end{array}
$$
with $\wh{N}_2(\la)=Q(\la)^T\otimes I_m$ such that $Q(\la)=\begin{bsmallmatrix}M_{\Phi}(\la)\\e_k^T\end{bsmallmatrix}^{-1} \begin{bsmallmatrix}
	I_{k-1}\\0\end{bsmallmatrix}$.

\item[\rm{(c)}]   $\begin{bsmallmatrix}
H_1(\la)\\H_2(\la)\\\vdots\\H_{k+1}(\la)
\end{bsmallmatrix}$ is a left minimal basis of $\mathbb{L}_2(\la)$ if and only if  $H_{k+1}(\la)$ is a left minimal basis of $G(\la)$ and
$$\begin{array}{l}H_1(\la)=(C(\la I_n-A)^{-1}T^{-1})^T H_{k+1}(\la),\\
H_i(\la)=\phi_{k-i+1}(\la)H_{k+1}(\la),\  i=2,\ldots ,k.
\end{array}
$$

\item[\rm{(d)}] If $\varepsilon_1\leq\cdots\leq\varepsilon_l$ are the right minimal indices of $G(\la)$ then $\varepsilon_1\leq\cdots\leq\varepsilon_l$ are the right minimal indices of $\mathbb{L}_2(\la)$.

\item[\rm{(e)}] If $\eta_1\leq\cdots\leq\eta_l$ are the left minimal indices of $G(\la)$ then $\eta_1+k-1\leq\cdots\leq\eta_l+k-1$ are the left minimal indices of $\mathbb{L}_2(\la)$.
\end{enumerate}
\end{theorem}	

\textbf{Proof}.- The proof can be done by using Lemmas \ref{lem_trans}(a) and \ref{lem_premult}, and Theorems \ref{th_bases_sbml} and \ref{th_indices_sbml}, and by following the same pattern as in the proof of Theorem \ref{th.m1recovery}.
\hfill\halmos

\begin{rem}
 \rm Comments similar to those in Remark \ref{rem.recoveryfro3} can be done in order to apply Theorem \ref{thm.M2recovery} to recover minimal bases of $G(\la)$ from those of any of its $\M_2$-strong linearizations. The only difference to be emphasized is that the roles of left and right minimal bases are interchanged in Theorems \ref{th.m1recovery} and \ref{thm.M2recovery}.
\end{rem}

{ Results completely analogous after obvious modifications to those in Subsection \ref{subsec.newextendedm1} can be obtained in the $\mathbb{M}_2$-framework. They are omitted for brevity.}

{
\section{Fiedler-like linearizations of rational matrices are block permutations of strong block minimal bases linearizations. Consequences.} \label{sec.Fiedler}

Fiedler pencils (FP), generalized Fiedler pencils (GFP), Fiedler pencils with repetition (FPR) and
generalized Fiedler pencils with repetition (GFPR) of square rational matrices $G(\la)$ have been introduced in the references \cite{AlBe16,AlBe18,DasAl19laa,DaAl20}, respectively. In simple words, the main idea in those definitions is to combine the  corresponding class of Fiedler-like pencils of the polynomial part of $G(\la)$ (introduced originally in \cite{Fiedler03,AntVol04,VolAnt11,BuDoFuRych15}) with a minimal order state-space realization of its strictly proper part in order to construct minimal order linear polynomial system matrices, of other rational matrices $\widehat{G} (\la)$, which are strong linearizations of $G(\la)$. We emphasize that FPRs and GFPRs are strong linearizations of $G(\la)$ under certain non-singularity generic hypotheses that, for simplicity, are assumed to hold in this section. We will recall that such assumptions hold by using the wordings ``FPR linearization'' and ``GFPR linearization''. Moreover, among the GFPs we only consider the {\em proper} GFPs (see \cite{AlBe18} and the references therein), since nonproper ones involve the inversion of some coefficients of the polynomial part of $G(\la)$ and are less interesting in applications.

It was proved in \cite[Section 4]{Dopico:2016:BlockKronecker} that every FP of a square polynomial matrix $D(\la)$ is a (block) permutation of a (particular) strong block minimal bases linearization of $D(\la)$. This was extended to proper GFPs, FPR linearizations and GFPR linearizations of $D(\la)$ in \cite{BuDoPe18}. This result allows to unify many different classes of strong linearizations of polynomial matrices and derive many properties of Fiedler-like linearizations from strong block minimal bases linearizations.

We prove in this section that, not surprisingly, FPs, proper GFPs, FPR linearizations and GFPR linearizations of a square rational matrix $G(\la)$ are also (block) permutations of particular strong block minimal bases linearizations of $G(\la)$. Moreover, the involved permutations are direct sums of an identity matrix plus the permutations corresponding to the Fiedler-like pencils of the polynomial part of $G(\la)$. We will discuss how the structure of the permutations allows us to prove inmediately that FPs, proper GFPs, FPR linearizations and GFPR linearizations of rational matrices are strong linearizations as a consequence of \cite[Theorem 5.11]{AmDoMaZa18} and \cite[Lemma 2.7]{DoMaQu19}. In addition, we will obtain as corollaries of Theorems \ref{th_bases_sbml} and \ref{th_indices_sbml}, the recovery rules of the minimal bases and indices from FPs, proper GFPs, FPR linearizations and GFPR linearizations of square rational matrices previously deduced in \cite{DasAl19laa,DaAl20} (see \cite{AlBe16-2} for a related result for eigenvectors).

We remark that the permutation results in this section were  mentioned very briefly (without proof) in \cite{AmDoMaZa18,DoMaQu19} and that are simple consequences of results available in the literature. They set  most Fiedler-like linearizations of rational matrices into the unified framework of strong block minimal bases linearizations. For proving the results, we need to recall some definitions on strong block minimal bases linearizations and Fiedler-like linearizations. In this section we consider only square polynomial and rational matrices, since most Fiedler-like pencils have been defined only in the square case.

\subsection{Extended block Kronecker linearizations and the antidiagonal sum condition} Let $L_k (\la)$ be the pencil defined in \eqref{eq.bloquekronbasic}, and
\begin{equation} \label{eq.lambdak}
\Lambda_k (\la) :=
\begin{bmatrix}
\la^k & \la^{k-1} & \cdots & \la & 1
\end{bmatrix}.
\end{equation}
It is well-known and easy to see that $L_k (\la)$ and $\Lambda_k (\la)$ are dual minimal bases, as well as $L_k (\la)\otimes I_p$ and $\Lambda_k (\la) \otimes I_p$ for any integer $p>0$ \cite{Dopico:2018:BlockKronecker}. Observe that, for any nonsingular constant matrix $Y$, also $Y(L_k (\la)\otimes I_p)$ and $\Lambda_k (\la) \otimes I_p$ are dual minimal bases. Moreover the row degrees of $Y(L_k (\la)\otimes I_p)$ are all equal to $1$ and the ones of $\Lambda_k (\la) \otimes I_p$ are all equal to $k$.

Given $D(\la) = D_q \la^q + D_{q-1} \la^{q-1} + \cdots + D_1 \la + D_0 \in \FF[\la]^{p\times p}$, we consider the family of strong block minimal bases linearizations associated to $D(\la)$ defined by setting in $\mathcal{L} (\la)$ in \eqref{eq:minbaspencil}
\begin{equation}\label{eq.extendedblock1}
K_1 (\la) = Y (L_\varepsilon (\la) \otimes I_p), \quad \quad \quad
K_2 (\la) = Z (L_\eta (\la) \otimes I_p),
\end{equation}
where $Y$ and $Z$ are nonsingular matrices and $q = \varepsilon + \eta + 1$, and $M(\la)$ such that
\begin{equation} \label{eq.extendedblock2}
D(\la) = (\Lambda_\eta (\la) \otimes I_p) M(\la) (\Lambda_\varepsilon (\la)^T \otimes I_p) \, .
\end{equation}
For simplicity, we assume that $D_q \ne 0$, which implies that $\mathcal{L} (\la)$ has sharp degree. Moreover, if $M(\la) = M_1 \la + M_0$ is partitioned into $(\eta + 1) \times (\varepsilon + 1)$ blocks each of size $p \times p$ and we denote such blocks by $[M(\la)]_{ij} = [M_1]_{ij} \la + [M_0]_{ij}$, $1\leq i \leq \eta + 1$ and $1 \leq j \leq \varepsilon + 1$ , then \eqref{eq.extendedblock2} is equivalent to
\begin{equation}\label{eq:antidiagonalsumcondition}
   \sum_{i+j=q+2-k} [M_1]_{ij} + \sum_{i+j=q+1-k} [M_0]_{ij}=D_k, \quad \mbox{for $k=0,1,\hdots,q$.}
 \end{equation}
This condition follows from \cite[Theorem 4.4]{Dopico:2018:BlockKronecker} and, if \eqref{eq:antidiagonalsumcondition} is satisfied, it is said in \cite[Definition 3.9]{BuDoPe18} that ``$M(\la)$ satisfies the antidiagonal sum (AS) condition for $D(\la)$''. Thus, using the terminology in \cite[Section 3]{BuDoPe18}, we call the strong block minimal bases linearizations defined by \eqref{eq.extendedblock1} and \eqref{eq.extendedblock2} {\em extended $(\varepsilon,p,\eta,p)$-block Kronecker linearizations satisfying the AS condition for $D(\la)$}. All of them are strong linearizations of $D(\la)$ and, for brevity, we often refer to them simply as extended block Kronecker linearizations\footnote{We remark that we are simplifying a bit the terminology used in \cite[Section 3]{BuDoPe18}. First, in \cite{BuDoPe18} the term ``extended block Kronecker {\em pencil}'' is used since  $Y$ and $Z$ in \eqref{eq.extendedblock1} are allowed to be singular. However, see \cite[Theorem 3.8]{BuDoPe18}, such pencils are linearizations only when $Y$ and $Z$ are nonsingular. In addition, in \cite{BuDoPe18} expressions like ``extended block Kronecker pencil with body $M(\la)$ satisfying the AS condition for $D(\la)$'' are used instead of ``extended block Kronecker linearization satisfying the AS condition for $D(\la)$''. The reason of this is that other AS conditions are investigated in \cite{BuDoPe18}.} for $D(\la)$. Obviously, they include the well known block Kronecker linearizations introduced in \cite[Section 4]{Dopico:2018:BlockKronecker} (see also Remark \ref{rem.frorecov1}) just by taking $Y$ and $Z$ identity matrices. Extended block Kronecker linearizations of polynomial matrices have also been analyzed (with other name) in \cite{FasSal18}.

The following very simple corollary of the AS condition \eqref{eq:antidiagonalsumcondition} is key for obtaining the main results of this section. For block Kronecker pencils, it was described in \cite[eq. (4.5)]{Dopico:2018:BlockKronecker}.

\begin{corollary} \label{cor.thea0block1} Let $D(\la) = D_q \la^q +  \cdots + D_1 \la + D_0 \in \FF[\la]^{p\times p}$ and let
$$
\mathcal{L}(\la) =
  \left[
    \begin{array}{cc}
      M(\la) & (Z (L_\eta (\la) \otimes I_p))^T \\
      Y (L_\varepsilon (\la) \otimes I_p) &0
      \end{array}
    \right] \in \FF [\la]^{pq \times pq}
$$
be an extended $(\varepsilon,p,\eta,p)$-block Kronecker linearization satisfying the AS condition for $D(\la)$. If $M(\la) = M_1 \la + M_0$ is partitioned into $(\eta + 1) \times (\varepsilon + 1)$ blocks each of size $p \times p$, then
$$
[M_0]_{\eta +1 , \varepsilon + 1} = D_0.
$$
%i.e., the zero coefficient of the lower-right block of $M(\la)$ is precisely %$D_0$.
\end{corollary}

The AS condition \eqref{eq:antidiagonalsumcondition} and Corollary \ref{cor.thea0block1} make it convenient to view any extended $(\varepsilon,p,\eta,p)$-block Kronecker linearization satisfying the AS condition for any polynomial matrix $D(\la)\in \FF [\la]^{p \times p}$ of degree $q$ partitioned into $q \times q$ blocks each of size $p\times p$. For brevity we will refer to this partition as the {\em natural partition} of a extended block Kronecker linearization.

Based on the definitions above, we define extended $(\varepsilon,p,\eta,p)$-block Kronecker linearizations satisfying the AS condition for a rational matrix $G(\la) \in \FF (\la)^{p \times p}$, which are particular cases of strong block minimal bases linearizations of $G(\la)$. For this purpose, we first express $G(\la) = D(\la) + G_{sp} (\la) = D_q \la^q + \cdots + D_1 \la + D_0 + C(\la I_n - A)^{-1} B$ in terms of the coefficients in the monomial basis of the polynomial part $D(\la)$ of $G(\la)$ and of a minimal order state-space realization $G_{sp} (\la) = C(\la I_n - A)^{-1} B$ of the strictly proper part of $G(\la)$. Then we use \eqref{eq:ratstrongblockmin} with $K_1(\la)$ and $K_2(\la)$ as in \eqref{eq.extendedblock1}, with $M(\la)$ satisfying \eqref{eq.extendedblock2}, and $\widehat{K}_1 = e_{\varepsilon + 1}^T \otimes I_p$ and $\widehat{K}_2 = e_{\eta + 1}^T \otimes I_p$, where $e_k$ is the last column of $I_k$. The expressions for $\widehat{K}_1$ and $\widehat{K}_2$ follow from \cite[Example 5.6]{AmDoMaZa18} (see also Remark \ref{rem.frorecov1}). This leads (for any $T, S$ nonsingular matrices) to
\begin{equation} \label{eq:ratextendedblock}
L(\la) = \left[
\begin{array}{c|cc}
T(\la I_n -A)S & \phantom{a} e_{\varepsilon + 1}^T \otimes (TB) \phantom{a} & 0 \\ \hline \phantom{\Big|}
-e_{\eta + 1} \otimes (C S) \phantom{\Big|}& M(\la) & (Z (L_\eta (\la) \otimes I_p))^T \\
0 & Y (L_\varepsilon (\la) \otimes I_p) & 0
\end{array}
\right],
\end{equation}
which are the desired extended $(\varepsilon,p,\eta,p)$-block Kronecker linearizations satisfying the AS condition for a rational matrix $G(\la) \in \FF (\la)^{p \times p}$. The natural partition of $\mathcal{L} (\la)$ in Corollary \ref{cor.thea0block1} induces the {\em natural partition} of $L(\la)$ in \eqref{eq:ratextendedblock}, which has the $(1,1)$-block of size $n\times n$, the other blocks in the first block row (resp. column) of size $n\times p$ (resp. $p\times n$), and the remaining blocks of size $p\times p$. In the case $Y$ and $Z$ are identity matrices, we simply say that \eqref{eq:ratextendedblock} is an $(\varepsilon,p,\eta,p)$-block Kronecker linearization satisfying the AS condition for $G(\la) \in \FF (\la)^{p \times p}$.

\begin{rem} \label{rem.a0block2} {\rm
The simple facts that in the natural partition of $L(\la)$ in \eqref{eq:ratextendedblock} the blocks $-CS$ and $TB$ are, together with the $(1,1)$-block, the only non-identically zero blocks in the first block column and row, respectively, and that $-CS$ and $TB$ are located, respectively, in the rows and columns corresponding to the block $[M_0]_{\eta +1 , \varepsilon + 1} = D_0$ in Corollary \ref{cor.thea0block1} will be important for obtaining the results in this section.}
\end{rem}

\subsection{Fiedler-like linearizations} The definitions of the families of Fiedler-like linearizations of polynomial and rational matrices require a good number of concepts and notations. Several of them are omitted for brevity. For polynomial matrices, we refer to the summary in \cite[Section 4]{BuDoPe18}, since we follow similar notations and definitions.

Given a polynomial matrix $D(\la) = D_q \la^q +  \cdots + D_1 \la + D_0 \in \FF[\la]^{p\times p}$, its associated Fiedler-like pencils are defined in terms of products of the $pq \times pq$ elementary matrices introduced in \cite[pp. 66-67]{BuDoPe18}, which are constructed from a $p\times p$ matrix $X$ and are denoted by $M_i (X)$, where $i \in \{ -q, \ldots, 0, \ldots , q \}$ is an index that determines $M_i (X)$. If the matrix $X$ is the coefficient of $D(\la)$ in \cite[p. 67]{BuDoPe18}, then $M_i (X)$ is denoted simply by $M_i^D$. Index tuples, i.e., finite ordered sequences of integers belonging to $\{ -q, \ldots, 0, \ldots , q\}$, and matrix assignments are used to denote in a concise way products of elementary matrices. For instance, let $\mathbf{t} = (t_1,t_2,t_3)$ and $\mathcal{X} = (X_1,X_2,X_3)$ be an index tuple and a matrix assignment for $\mathbf{t}$, respectively, then $M_\mathbf{t} (\mathcal{X}) := M_{t_1} (X_1) M_{t_2} (X_2) M_{t_3} (X_3)$. Moreover, if $\mathbf{t}$ is the empty tuple, then we define $M_\mathbf{t} (\mathcal{X}) := I_{pq}$. As in the case of extended block Kronecker linearizations, it is very convenient to view these elementary matrices and their products partitioned into $q\times q$ blocks each of size $p \times p$. We will call this partition again the {\em natural partition} of $M_\mathbf{t} (\mathcal{X})$ and $[M_\mathbf{t} (\mathcal{X})]_{jk}$, $1\leq j,k \leq q$, will denote the block of $M_\mathbf{t} (\mathcal{X})$ in the block position $(j , k)$.
A fundamental property is that the products $M_\mathbf{t} (\mathcal{X})$ of elementary matrices appearing in the definitions of Fiedler-like pencils are {\em operation-free} \cite{VolAnt11} (see also \cite[Definition 4.5]{BuDoFuRych15}), which means that their blocks in the natural partition are either the matrices in the matrix assignments, or $I_p$ or $0_p$ and that the positions of such blocks only depend on the index tuple $\mathbf{t}$, i.e., do not depend on the particular matrix assignment $\mathcal{X}$. In order to guarantee the operation-free property, we will require that some of the index tuples involved in the definitions of Fiedler-like pencils satisfy the Successor Infix Property (SIP) introduced in \cite[Definition 7]{VolAnt11} (see also \cite[Definition 4.6]{BuDoPe18}).

Next, we define the FPs, proper GFPs, FPR linearizations and GFPR linearizations of the polynomial matrix $D(\la)$ in terms of the following generic pencil (see \cite[Definition 4.29]{BuDoPe18}),
\begin{equation}\label{eq.fprgen}
\mathcal{L}_D(\lambda)=M_{{\boldsymbol\ell}_{t},{\boldsymbol\ell}_{z}}(\mathcal{X},\mathcal{ Z})(\lambda M^D_{\mathbf{z}}
-M^D_{\mathbf{t}})M_{\mathbf{r}_{z},\mathbf{r}_{t}}(\mathcal{W}, \mathcal{Y}).
\end{equation}

$\bullet$ If in \eqref{eq.fprgen} ${\boldsymbol\ell}_{t},{\boldsymbol\ell}_{z}, \mathbf{r}_{z},\mathbf{r}_{t}$ are empty tuples, $\mathbf{z} = -q$ and $\mathbf{t}$ is a permutation of $\{0,1,\ldots, q-1\}$, then $\mathcal{L}_D(\lambda)$ is a FP associated to $D(\la)$.

$\bullet$ If in \eqref{eq.fprgen} ${\boldsymbol\ell}_{t},{\boldsymbol\ell}_{z}, \mathbf{r}_{z},\mathbf{r}_{t}$ are empty tuples, $\mathbf{t}$ and $\mathbf{z}$ are permutations, respectively, of $C_0$ and $-C_1$, where $\{C_0, C_1 \}$ is a partition of $\{0,1,\ldots, q\}$ such that $0\in C_0$ and $q \in C_1$, then $\mathcal{L}_D(\lambda)$ is a proper GFP associated to $D(\la)$.

$\bullet$ Let $h\in\left\{  0, 1, \ldots , q-1\right\}$. If in \eqref{eq.fprgen} ${\mathbf{t}}$ and $\mathbf{z}$ are permutations of $\left\{  0,\ldots, h\right\}
$ and $\left\{  -q, -q+1, \ldots, -h-1\right\}$, respectively, ${\boldsymbol\ell}_{t}$ and
${\mathbf{r}}_{t}$ are tuples with indices from $\left\{  0, \ldots , h-1\right\}  $ such
that $({\boldsymbol\ell}_{t},{\mathbf{t,r}}_{t})$ satisfies the SIP, and
${\boldsymbol\ell}_{z}$ and ${\mathbf{r}}_{z}$ are tuples with indices from
$\left\{  -q, -q+1, \ldots , -h-2\right\}  $ such that $({\boldsymbol\ell}_{z},{\mathbf{z,r}}_{z})$ satisfies the SIP, then $\mathcal{L}_D(\lambda)$ is a GFPR associated to $D(\la)$. This pencil is a GFPR linearization whenever the matrices $M_{{\boldsymbol\ell}_{t},{\boldsymbol\ell}_{z}}(\mathcal{X},\mathcal{ Z})$ and
$M_{\mathbf{r}_{z},\mathbf{r}_{t}}(\mathcal{W}, \mathcal{Y})$ are nonsingular, or, equivalently, the matrix assignments $\mathcal{X},\mathcal{ Z},\mathcal{W}, \mathcal{Y}$ are nonsingular according to \cite[Definition 4.1]{BuDoFuRych15}. The FPRs associated to $D(\la)$ are those GFPRs such that $\mathcal{X},\mathcal{ Z},\mathcal{W}, \mathcal{Y}$ are the trivial matrix assignments for ${\boldsymbol\ell}_{t},{\boldsymbol\ell}_{z}, \mathbf{r}_{z},\mathbf{r}_{t}$, respectively (see \cite[p. 67]{BuDoPe18}).

The operation-free property of the products of elementary matrices in $\mathcal{L}_D(\lambda) := \la \mathcal{L}_1 + \mathcal{L}_0$ in \eqref{eq.fprgen} for FPs, proper GFPs and GFPRs and the definition of the involved index tuples imply that for every FP, for every proper GFP, and for every GFPR for which the matrix assignments $\mathcal{X},\mathcal{ Z},\mathcal{W}, \mathcal{Y}$ do not contain $-D_0$ (the matrix coefficient of degree zero of $D(\la)$), there is exactly one block position in the natural partition of $\mathcal{L}_0$ whose entry is identically equal\footnote{The expression ``identically equal'' means that this happens always in the same block entry for any value of $D_0$.} to $D_0$. This property is very easy to prove for FPs and proper GFPs; for GFPRs, it can be found in \cite[Theorem 5.3]{BuDoFuRych15}.  Moreover, this block position in the natural partition of $\mathcal{L}_0$ is uniquely determined by the index tuples in \eqref{eq.fprgen}, i.e., does not depend on the particular matrix assignment in the case of GFPRs, and we call it the {\em intrinsic block position of $D_0$ in $\mathcal{L}_0$}. We use the word ``intrinsic'' because for other GFPRs with the same index tuples and with some of the matrix assignments $\mathcal{X},\mathcal{ Z},\mathcal{W}, \mathcal{Y}$ containing $-D_0$, $D_0$ is the entry of the ``intrinsic block position'' and of other block positions of $\mathcal{L}_0$. However, under arbitrary changes of the matrix assignments (including symbolic changes of variables), $D_0$ only remains in the ``intrinsic position''. This allows us to define the intrinsic block position of $D_0$ in $\mathcal{L}_0$ for the family of all GFPRs sharing the same index tuples in \eqref{eq.fprgen}.

The intrinsic block position of $D_0$ in $\mathcal{L}_0$ can be easily determined by using the consecutive consecutions and the consecutive inversions at $0$ of two index tuples (see \cite[Definition 5.7]{DasAl19laa} or \cite[Definition 2.8]{DaAl20} for the definition of consecutive consecutions and inversions). This is stated in Lemma \ref{lem.D_0_pos}, which is the only new result so far in this section.

\begin{lemma} \label{lem.D_0_pos} Let $D(\la) = D_q \la^q +  \cdots + D_1 \la + D_0 \in \FF[\la]^{p\times p}$ be a polynomial matrix and let $\mathcal{L}_D(\lambda) := \la \mathcal{L}_1 + \mathcal{L}_0$ as in \eqref{eq.fprgen} be a FP or a proper GFP or a GFPR associated to $D(\la)$. Let
$i_0 (\boldsymbol\ell_t,\mathbf{t})$ be the number of consecutive inversions at $0$ of the index tuple $(\boldsymbol\ell_t,\mathbf{t})$ and $c_0 (\mathbf{t},\mathbf{r}_t)$ be the number of consecutive consecutions at $0$ of $(\mathbf{t},\mathbf{r}_t)$. Then the intrinsic block position of $D_0$ in $\mathcal{L}_0$ is $(q- i_0 (\boldsymbol\ell_t,\mathbf{t}), q - c_0 (\mathbf{t},\mathbf{r}_t))$.
\end{lemma}

\textbf{Proof}.- Once the (highly non-trivial) operation-free property is taken into account, the proof is a simple exercise of block matrix multiplication. For brevity, we only sketch the proof for GFPRs with matrix assignments $\mathcal{X}$, $\mathcal{ Z}$, $\mathcal{W}$, $\mathcal{Y}$ not containing $-D_0$, since the proofs of the other cases are similar. Note first that the commutativity relations of the elementary matrices  \cite[Remark 4.3]{BuDoPe18} imply that the elementary matrices in $M_{{\boldsymbol\ell}_{z}}(\mathcal{ Z})$ and $M_{\mathbf{r}_{z}}(\mathcal{W})$ commute with the other elementary matrices in $\mathcal{L}_0$ and, so, they do not affect the position of $D_0$ (see also \cite[Lemma 4.33 and p. 76]{BuDoPe18}). Next, note that $D_0$ is in the block position $(q,q)$ if $i_0 (\boldsymbol\ell_t,\mathbf{t}) = c_0 (\mathbf{t},\mathbf{r}_t) = 0$, as a consequence of the commutativity relations of the elementary matrices. Otherwise, note that each elementary matrix different from $M_0^D$ corresponding to the subtuple of $(\boldsymbol\ell_t,\mathbf{t})$ (resp. $(\mathbf{t}, \mathbf{r}_t)$) containing the index $0\in \mathbf{t}$ and defining $i_0 (\boldsymbol\ell_t,\mathbf{t})$ (resp. $c_0 (\mathbf{t}, \mathbf{r}_t)$) moves $D_0$ one position up (resp. left), while the rest of elementary matrices in $(\boldsymbol\ell_t,\mathbf{t}, \mathbf{r}_t)$  do not move $D_0$. \hfill\halmos

\begin{rem} \label{re.i0c0} {\rm We emphasize that $(i_0 (\boldsymbol\ell_t,\mathbf{t}), c_0 (\mathbf{t},\mathbf{r}_t))\allowbreak =  (i_0 (\mathbf{t}), c_0 (\mathbf{t}))$  for FPs and proper GFPs. For FPs, one of the elements in this pair is necessarily different from zero and the other one is necessarily zero. For proper GFPs, at most one is different from zero, but it may happen that both are zero.
}
\end{rem}

Based on the definitions of FPs, proper GFPs, FPRs and GFPRs of polynomial matrices, FPs, proper GFPs, FPRs and GFPRs of square rational matrices have been defined in \cite[Definition 3.2]{AlBe16}, \cite[Definition 2.2]{AlBe18}, \cite[Definition 5.4]{DasAl19laa} and \cite[Definition 3.2]{DaAl20}, respectively. As in the case of polynomial matrices, for rational matrices GFPRs include FPRs, and GFPRs are strong linearizations if the involved matrix assignments are nonsingular. Essentially, the strategy in \cite{AlBe16,AlBe18,DasAl19laa,DaAl20} is to consider the rational matrix $G(\la) = D(\la) + G_{sp} (\la) = D_q \la^q +  \cdots + D_1 \la + D_0 + C(\la E - A)^{-1} B \in \FF(\la)^{p\times p}$ expressed in terms of the coefficients in the monomial basis of the polynomial part $D(\la)$ of $G(\la)$ and of a minimal order state-space realization $C(\la E - A)^{-1} B$ of the strictly proper part, with $E \in \FF^{n\times n}$ nonsingular. Then, the elementary matrices $M_i(X)$ and $M_i^D$ of the polynomial part $D(\la)$ are carefully embedded into augmented elementary matrices $\mathbb{M}_i (X)$ and $\mathbb{M}_i^G$ of the rational matrix that incorporate the information of the state-space realization of the strictly proper part. Finally, the corresponding families of Fiedler-like pencils of $G(\la)$ are defined through products of these augmented elementary matrices using matrix assignments and the same index tuples as in the corresponding families of Fiedler-like pencils of $D(\la)$. Theorem \ref{thm.FiedlerRatvsPol} is a very important result in this context that summarizes in a concise way  \cite[Theorem 3.6]{AlBe16}, \cite[Theorem 2.7]{AlBe18}, \cite[Theorem 5.12]{DasAl19laa} and \cite[Theorem 3.6]{DaAl20}.

\begin{theorem} \label{thm.FiedlerRatvsPol} Let $G(\la) = D(\la) + G_{sp} (\la)\in \FF(\la)^{p\times p}$ be the unique decomposition of $G(\la)$ into its polynomial part $D(\la) = D_q \la^q +  \cdots + D_0$ and its strictly proper part $G_{sp} (\la)$, and let $C(\la E - A)^{-1} B$ be a minimal order state-space realization of $G_{sp}(\la)$ (with $E\in \FF^{n\times n}$ nonsingular). Let
$L_G(\lambda)= \mathbb{M}_{{\boldsymbol\ell}_{t},{\boldsymbol\ell}_{z}}(\mathcal{X},\mathcal{ Z})(\lambda \mathbb{M}^G_{\mathbf{z}}
-\mathbb{M}^G_{\mathbf{t}})\mathbb{M}_{\mathbf{r}_{z},\mathbf{r}_{t}}(\mathcal{W}, \mathcal{Y})$ and $\mathcal{L}_D(\lambda)=M_{{\boldsymbol\ell}_{t},{\boldsymbol\ell}_{z}}(\mathcal{X},\mathcal{ Z})(\lambda M^D_{\mathbf{z}}
-M^D_{\mathbf{t}})M_{\mathbf{r}_{z},\mathbf{r}_{t}}(\mathcal{W}, \mathcal{Y})$ be
FPs, or proper GFPs or GFPRs associated to $G(\la)$ and to $D(\la)$, respectively. Then
\begin{equation}\label{eq.FiedlerRatvsPol}
L_G (\la) =
\left[
\begin{array}{c|c}
A - \la E & e^T_{q - c_0 (\mathbf{t},\mathbf{r}_t)} \otimes B \\ \hline
e_{q- i_0 (\boldsymbol\ell_t,\mathbf{t})} \otimes C & \mathcal{L}_D(\lambda)
\end{array}
\right],
\end{equation}
where $e_j$ is the $j$th column of $I_q$.
\end{theorem}

\begin{rem} \label{rem.diffalam} {\rm We have followed in \eqref{eq.FiedlerRatvsPol} the classical notation in \cite{Rosen70} for polynomial system matrices and we have set $A-\la E$ in the upper-left corner. In contrast in
\cite{AlBe16,AlBe18,DasAl19laa,DaAl20}, $A-\la E$ is set in the lower-right corner and the two block rows and the two block columns in \eqref{eq.FiedlerRatvsPol} appear interchanged. Note also that in previous sections of this paper we have written the strictly proper part of $G(\la)$ as $G_{sp} (\la) = C(\la I_n - A)^{-1} B$, which is the same as $G_{sp} (\la) = (CS) (\la TS - TAS)^{-1} (TB)$ for any nonsingular $T$ and $S$ matrices. With an obvious change of notation, it can be written as $G_{sp} (\la) = C (\la E - A)^{-1} B$, with $E$ nonsingular. This allows to compare \eqref{eq.FiedlerRatvsPol} and \eqref{eq:ratextendedblock} in Subsection \ref{subsec.perm}.
}
\end{rem}

\begin{rem} \label{rem.a0block3} {\rm If we consider for \eqref{eq.FiedlerRatvsPol} the same natural partition as in \eqref{eq:ratextendedblock} and we take into account Lemma \ref{lem.D_0_pos}, then we obtain from Theorem \ref{thm.FiedlerRatvsPol} a simple recipe to construct the FPs, proper GFPs, FPRs and GFPRs of a rational matrix when the explicit expression of the corresponding pencil of its polynomial part is known: (1) construct any FP, proper GFP, FPR or GFPR of the polynomial part, $\mathcal{L}_D(\lambda) = \la \mathcal{L}_1 + \mathcal{L}_0$; (2) identify in $\mathcal{L}_0$ the intrinsic block position of $D_0$ (this is particularly simple in FPs, proper GFPs, and GFPRs with matrix assignments not containing $-D_0$, since there is only one block entry identically equal to $D_0$); (3) augment  $\mathcal{L}_D(\lambda)$ to a larger matrix partitioned in the natural way by adding one block column and one block row as follows: the $(1,1)$-block is $(A - \la E)$ and $C$ (resp. $B$) is the only remaining nonzero block in the first block column (resp. row) and is located in the block row (resp. column) of the intrinsic block position of $D_0$ in $\mathcal{L}_0$. This remark is related to Remark \ref{rem.a0block2} for extended block Kronecker linearizations of rational matrices and stresses the importance of the intrinsic position of $D_0$ in these families of linearizations.
}
\end{rem}

\subsection{Block permutations of Fiedler-like pencils} \label{subsec.perm}
Given two positive integers $q$ and $p$, we say that a matrix $\Pi$ is a  $(q,p)$-block permutation matrix if $\Pi = \Sigma \otimes I_p$, where $\Sigma$ is a $q\times q$ permutation matrix. The following result follows from Theorems 6.3, 7.1 and 8.1 of \cite{BuDoPe18}.

\begin{theorem} \label{thm.permut1}
Let $D(\la) = D_q \la^q + \cdots + D_1 \la + D_0 \in \FF[\la]^{p\times p}$ be a polynomial matrix and $\mathcal{L}_D (\la)$ be a FP, or a proper GFP or a GFPR with nonsingular matrix assignments associated to $D(\la)$. Then, there exist two $(q,p)$-block permutation matrices $\Pi_1$ and $\Pi_2$ such that
\begin{equation}\label{eq.permut1}
\Pi_1 \mathcal{L}_D (\la) \Pi_2 = \left[
    \begin{array}{cc}
      M(\la) & (Z (L_\eta (\la) \otimes I_p))^T \\
      Y (L_\varepsilon (\la) \otimes I_p) &0
      \end{array}
    \right] \in \FF [\la]^{pq \times pq}
\end{equation}
is an extended $(\varepsilon,p,\eta,p)$-block Kronecker linearization satisfying the AS condition for $D(\la)$. Moreover the parameters $\varepsilon$ and $\eta$ are determined uniquely from the index tuples defining $\mathcal{L}_D (\la)$. If $\mathcal{L}_D (\la)$ is a FP or a proper GFP, then $Y$ and $Z$ in \eqref{eq.permut1} are identity matrices and $\Pi_1 \mathcal{L}_D (\la) \Pi_2$ is a $(\varepsilon,p,\eta,p)$-block Kronecker linearization.
\end{theorem}

Next, we extend Theorem \ref{thm.permut1} to Fiedler-like and extended block Kronecker linearizations of rational matrices. We warn the reader that for writing \eqref{eq.permut2} in a compact way, the vectors $e_{\varepsilon +1}$ and $e_{\eta +1}$ in \eqref{eq.permut2} are different from those in \eqref{eq:ratextendedblock}: in \eqref{eq.permut2} they are the corresponding columns of $I_q$, while in \eqref{eq:ratextendedblock} they are the last columns of $I_{\varepsilon +1}$ and $I_{\eta +1}$, respectively.

\begin{theorem} \label{thm.permut2} Let $G(\la) = D(\la) + G_{sp} (\la)\in \FF(\la)^{p\times p}$ be the unique decomposition of $G(\la)$ into its polynomial part $D(\la) = D_q \la^q +  \cdots + D_0$ and its strictly proper part $G_{sp} (\la)$, and let $C(\la E - A)^{-1} B$ be a minimal order state-space realization of $G_{sp}(\la)$ (with $E \in \FF^{n \times n}$ nonsingular). Let $L_G(\lambda)$ be a FP, or a proper GFP or a GFPR with nonsingular matrix assignments associated to $G(\la)$ as in \eqref{eq.FiedlerRatvsPol}, where $\mathcal{L}_D(\lambda)$ is the corresponding FP, or proper GFP or GFPR associated to $D(\la)$. Then, there exist two $(q,p)$-block permutation matrices $\Pi_1$ and $\Pi_2$ such that
\begin{equation}\label{eq.permut2}
\begin{bmatrix}
  I_n  &  \\
    & \Pi_1
\end{bmatrix}
L_G (\la)
\begin{bmatrix}
  I_n  &  \\
    & \Pi_2
\end{bmatrix}
= \left[
\begin{array}{c|c}
A - \la E &  e^T_{\varepsilon +1} \otimes B  \\ \hline
e_{\eta +1 } \otimes C & \Pi_1 \mathcal{L}_D(\lambda) \Pi_2
\end{array}
\right],
\end{equation}
is an extended $(\varepsilon,p,\eta,p)$-block Kronecker linearization satisfying the AS condition for $G(\la)$. Moreover the parameters $\varepsilon$ and $\eta$ are determined uniquely from the index tuples defining $L_G (\la)$. If $L_G (\la)$ is a FP or a proper GFP, then the pencil in the right-hand side of \eqref{eq.permut2} is an $(\varepsilon,p,\eta,p)$-block Kronecker linearization of $G(\la)$.
\end{theorem}

\textbf{Proof}.- Let $\Pi_1$ and $\Pi_2$ be the two block permutation matrices in Theorem \ref{thm.permut1}. Then, from \eqref{eq.FiedlerRatvsPol}, we get
\begin{equation}\label{eq.permut3}
\begin{bmatrix}
  I_n  &  \\
    & \Pi_1
\end{bmatrix}
L_G (\la)
\begin{bmatrix}
  I_n  &  \\
    & \Pi_2
\end{bmatrix}
= \left[
\begin{array}{c|c}
A - \la E &  (e^T_{q - c_0 (\mathbf{t},\mathbf{r}_t)} \otimes B) \Pi_2  \\ \hline
\Pi_1 (e_{q- i_0 (\boldsymbol\ell_t,\mathbf{t})} \otimes C) & \Pi_1 \mathcal{L}_D(\lambda) \Pi_2
\end{array}
\right],
\end{equation}
where $\Pi_1 \mathcal{L}_D(\lambda) \Pi_2$ is the extended $(\varepsilon,p,\eta,p)$-block Kronecker linearization satisfying the AS condition for $D(\la)$ in Theorem \ref{thm.permut1}. According to \eqref{eq:ratextendedblock}, it only remains to prove that
\begin{equation} \label{eq.permut4}
(e^T_{q - c_0 (\mathbf{t},\mathbf{r}_t)} \otimes B) \Pi_2 = e^T_{\varepsilon +1} \otimes B \quad \mbox{and} \quad \Pi_1 (e_{q- i_0 (\boldsymbol\ell_t,\mathbf{t})} \otimes C) = e_{\eta +1 } \otimes C.
\end{equation}
The proof of these two equalities is simple in the case $L_G(\la)$ is a FP, a proper GFP, or a GFPR whose matrix assignments do not contain $-D_0$. The key fact in these cases is that if $\mathcal{L}_D(\lambda) := \la \mathcal{L}_1 + \mathcal{L}_0$ then there is exactly one block identically equal to $D_0$ in the natural partition of $\mathcal{L}_0$ located in the block position $(q- i_0 (\boldsymbol\ell_t,\mathbf{t}), q - c_0 (\mathbf{t},\mathbf{r}_t))$, according to Lemma \ref{lem.D_0_pos}. Therefore, $\Pi_1 \mathcal{L}_D(\lambda) \Pi_2 := \la \mathcal{T}_1 + \mathcal{T}_0$ has also exactly one block identically equal to $D_0$ in the natural partition of $\mathcal{T}_0$, since $\Pi_1$ and $\Pi_2$ are $(q,p)$-block permutations. Moreover, Corollary \ref{cor.thea0block1} implies that $(\eta+1, \varepsilon +1)$ is the block position of $D_0$ in $\mathcal{T}_0$. This implies that $\Pi_1$ moves the block row $q- i_0 (\boldsymbol\ell_t,\mathbf{t})$ to the block row $\eta +1$ and that $\Pi_2$ moves the block column  $q - c_0 (\mathbf{t},\mathbf{r}_t)$ to the block column $\varepsilon + 1$, which imply \eqref{eq.permut4}.

In the case of GFPRs whose matrix assignments contain $-D_0$, there are more than one blocks equal to $D_0$ in $\mathcal{L}_0$ (and, so, in $\mathcal{T}_0$), and the proof of \eqref{eq.permut4} requires to use the concept of the intrinsic block position of $D_0$ in $\mathcal{L}_0$, introduced before Lemma \ref{lem.D_0_pos}. Note that Theorem \ref{thm.permut1} implies that the parameters $\varepsilon$ and $\eta$ are determined uniquely by the index tuples defining $L_G (\la)$, which are the same as those defining $\mathcal{L}_D (\la)$, and, so, are the same for all the GFPRs with the same index tuples (independently of the matrix assignments). Thus, the intrinsic block position of $D_0$ in $\mathcal{L}_0$ in $(q- i_0 (\boldsymbol\ell_t,\mathbf{t}), q - c_0 (\mathbf{t},\mathbf{r}_t))$ is mapped by $\Pi_1$ and $\Pi_2$ to $(\eta + 1, \varepsilon + 1)$ in $\mathcal{T}_0$, since by Corollary \ref{cor.thea0block1} is the only block entry of $\mathcal{T}_0$ that contains $D_0$ when the matrix assignments change arbitrarily but the defining tuples do not change. This implies \eqref{eq.permut4}.
\hfill\halmos

\subsection{Some consequences of Theorem \ref{thm.permut2}: Fiedler-like pencils are strong linearizations and recovery of minimal bases from Fiedler-like pencils} \label{subsec.consec} We study three easy consequences of Theorem \ref{thm.permut2} and its proof. More precisely, (1) we provide alternative proofs to those in \cite{AlBe16,AlBe18,DasAl19laa,DaAl20} for the facts that FPs, proper GFPs, FPRs and GFPRs (with nonsingular matrix assignments) of rational matrices are strong linearizations of rational matrices; (2) we provide alternative proofs to those in \cite{DasAl19laa,DaAl20} for the recovery rules of the minimal bases of a rational matrix from those of its Fiedler-like pencils; (3) we discuss briefly how to recover minimal indices of a rational matrix from those of its Fiedler-like pencils, a problem that has been solved in \cite{DasAl19laa,DaAl20}.

\begin{corollary} \label{cor.consequence1}
Let $G(\la) \in \FF (\la)^{p \times p}$ be a rational matrix. If $L_G (\la)$ is a FP, or a proper GFP or a GFPR with nonsingular matrix assignments associated with $G(\la)$, then $L_G (\la)$ is a strong linearization of $G(\la)$.
\end{corollary}

\textbf{Proof}.- $L_G(\la)$ satisfies \eqref{eq.permut2} and the right-hand side of \eqref{eq.permut2} is a strong linearization of $G(\la)$ by \cite[Theorem 5.11]{AmDoMaZa18}. Then, \cite[Lemma 2.7]{DoMaQu19} implies that $L_G (\la)$ is also a strong linearization of $G(\la)$.
\hfill\halmos

\smallskip

Corollary \ref{cor.consequence2} covers the recovery of minimal bases. In order to check that the results in Corollary \ref{cor.consequence2} are the same as those in \cite{DasAl19laa,DaAl20} recall that for FPs and proper GFPs ${\boldsymbol\ell}_{t},{\boldsymbol\ell}_{z}, \mathbf{r}_{z},\mathbf{r}_{t}$ are empty tuples and also Remark \ref{rem.diffalam}.
\begin{corollary} \label{cor.consequence2} Let $G(\la) \in \FF(\la)^{p\times p}$ be as in Theorem \ref{thm.permut2}. Let
$L_G(\lambda)= \mathbb{M}_{{\boldsymbol\ell}_{t},{\boldsymbol\ell}_{z}}(\mathcal{X},\mathcal{ Z})(\lambda \mathbb{M}^G_{\mathbf{z}}
-\mathbb{M}^G_{\mathbf{t}})\mathbb{M}_{\mathbf{r}_{z},\mathbf{r}_{t}}(\mathcal{W}, \mathcal{Y})$ be a FP, or a proper GFP or a GFPR linearization associated to $G(\la)$, let
$i_0 (\boldsymbol\ell_t,\mathbf{t})$ be the number of consecutive inversions at $0$ of the index tuple $(\boldsymbol\ell_t,\mathbf{t})$, $c_0 (\mathbf{t},\mathbf{r}_t)$ be the number of consecutive consecutions at $0$ of $(\mathbf{t},\mathbf{r}_t)$ and $e_j$ be the $j$th column of $I_{q}$.
\begin{itemize}
  \item [\rm (a)] If $Q(\la) \in \FF[\la]^{(n + pq) \times l}$ is a right minimal basis of $L_G(\la)$ then
      $$
      \begin{bmatrix}
        0_{p \times n}  & e^T_{q-c_0 (\mathbf{t},\mathbf{r}_t)} \otimes I_p
      \end{bmatrix} Q(\la)
      $$
      is a right minimal basis of $G(\la)$.
  \item [\rm (b)] If $Q(\la) \in \FF[\la]^{(n + pq) \times l}$ is a left minimal basis of $L_G(\la)$ then
      $$
      \begin{bmatrix}
        0_{p \times n}  & e^T_{q-i_0 (\boldsymbol\ell_t,\mathbf{t})} \otimes I_p
      \end{bmatrix} Q(\la)
      $$
      is a left minimal basis of $G(\la)$.
\end{itemize}
\end{corollary}

\textbf{Proof}.- We prove part (a). Let $L(\la)$ be the extended $(\varepsilon,p,\eta,p)$-block Kronecker linearization satisfying the AS condition for $G(\la)$ in the right-hand side of \eqref{eq.permut2}. Note that $Q(\la) = \mbox{diag} (I_n, \Pi_2) \, \widetilde{Q}(\la)$, where $\widetilde{Q}(\la) \in \FF[\la]^{(n + pq) \times l}$ is a right minimal basis of $L(\la)$. Next, we apply Theorem \ref{th_bases_sbml}(a) to $L (\la)$ and $\widetilde{Q}(\la)$, taking into account that in this case $N_1 (\la)^T = \Lambda_\varepsilon (\la)^T \otimes I_p$,  with $\Lambda_\varepsilon (\la)$ as in \eqref{eq.lambdak}. Therefore, Theorem \ref{th_bases_sbml} implies that $$
      \begin{bmatrix}
        0_{p \times n}  & e^T_{\varepsilon +1} \otimes I_p
      \end{bmatrix} \widetilde{Q}(\la) = \begin{bmatrix}
        0_{p \times n}  & (e^T_{\varepsilon +1} \otimes I_p) \Pi_2^T
      \end{bmatrix} Q(\la)
      $$
is a right minimal basis of $G(\la)$. Finally, note that $(e^T_{\varepsilon +1} \otimes I_p) \Pi_2^T = e^T_{q-c_0 (\mathbf{t},\mathbf{r}_t)} \otimes I_p$, since we have seen in the proof of Theorem \ref{thm.permut2} that the block permutation $\Pi_2$ moves the block column  $q - c_0 (\mathbf{t},\mathbf{r}_t)$ to the block column $\varepsilon + 1$. Part (b) is proved analogously  via Theorem \ref{th_bases_sbml}(b) and $N_2 (\la)^T = \Lambda_\eta (\la)^T \otimes I_p$.
\hfill\halmos

\begin{rem} \label{rem.a0block4} {\rm In the same spirit of Remark \ref{rem.a0block3} on Theorem \ref{thm.FiedlerRatvsPol}, Lemma \ref{lem.D_0_pos} allows to express the recovery rules in Corollary \ref{cor.consequence2} in simple words as follows: (1) identify in the explicit expression of $L_G(\la) = \la L_1 + L_0$ the intrinsic block position of $D_0$ in the natural partition of $L_0$ (which is immediate for FPs, proper GFPs and GFPRs with matrix assignments not containing $-D_0$); (2) the rows of any right (resp. left) minimal basis of $L_G (\la)$ corresponding to the columns (resp. rows) of the intrinsic block position of $D_0$ in $L_0$ are a right (resp. left) minimal basis of $G(\la)$.}
\end{rem}

Though Corollary \ref{cor.consequence2} is enough for recovering the minimal bases of $G(\la)$ from those of its Fiedler-like pencils in applications, we remark that it is weaker than Theorem \ref{th_bases_sbml} for the minimal bases of strong block minimal bases linearizations. In contrast to Theorem \ref{th_bases_sbml},  Corollary \ref{cor.consequence2} does not allow to construct the minimal bases of $L_G(\la)$ from those of $G(\la)$. In this context, note that for block Kronecker linearizations the matrices $\widehat{N}_1 (\la)$ and $\widehat{N}_2 (\la)$ in Theorem \ref{th_bases_sbml} are known and have simple expressions \cite[Remark 7.5]{Dopico:2016:BlockKronecker}. In the case of extended block Kronecker linearizations, expressions of $\widehat{N}_1 (\la)$ and $\widehat{N}_2 (\la)$ involving the inverses of the matrices $Y$ and $Z$ in \eqref{eq:ratextendedblock} can also be obtained.

Finally, note that \eqref{eq.permut2} implies that the minimal indices of $L_G (\la)$ are those of the extended $(\varepsilon,p,\eta,p)$-block Kronecker linearization of $G(\la)$ in the right hand-side. Then, Theorem \ref{th_indices_sbml} with $N_1 (\la) = \Lambda_\varepsilon (\la) \otimes I_p$ and $N_2 (\la) = \Lambda_\eta (\la) \otimes I_p$ imply that the right (resp. left) minimal indices of $G(\la)$ are those of $L_G (\la)$ minus $\varepsilon$ (resp. $\eta$). The parameters $\varepsilon$ and $\eta$ can be obtained from the index tuples defining $L_G (\la)$ as explained in Theorems 6.3, 7.1 and 8.1 of \cite{BuDoPe18}, which requires to use a number of definitions related to index tuples that we omit for brevity.

\section{Some remarks on eigenvectors} \label{sec.eigenvectors} If a rational matrix $G(\la)$ is regular, i.e., it is square and $\det G(\la)$ is not identically zero, then it does not have minimal bases nor minimal indices. In this case, $\la_0 \in \overline{\FF}$ is an eigenvalue of $G(\la)$ if $\la_0$ is a finite zero but not a pole of $G(\la)$. Then, there exist nonzero vectors $x,y$ such that $y^T G(\la_0) = 0$ and $G(\la_0) x = 0$. Such vectors are called, respectively, left and right eigenvectors of $G(\la)$ associated to $\la_0$ and are very interesting magnitudes in rational and nonlinear eigenvalue problems \cite{guttel-tisseur}. Equivalently, the left (resp. right) eigenvectors associated to $\la_0$ are the nonzero vectors of the left (resp. right) null-space of the constant matrix $G(\la_0)$. We denote such null-spaces by $\mathcal{N}_\ell (G(\lambda_0))$ and $\mathcal{N}_r (G(\lambda_0))$, respectively. A standard method for computing eigenvectors of $G(\la)$, or more precisely bases of $\mathcal{N}_\ell (G(\lambda_0))$ and $\mathcal{N}_r (G(\lambda_0))$, is to compute those of one of its linearizations and to recover from them the eigenvectors of $G(\la)$. This has led to recovery procedures for eigenvectors from Fiedler-like linearizations \cite{AlBe16-2,DasAl19laa,DaAl20}, from $\mathbb{M}_1$ and $\mathbb{M}_2$-strong linearizations \cite{DoMaQu19}, and from strong linearizations in the affine spaces defined in \cite{DaAl19_2}. For completeness, we describe in this section very briefly how to recover eigenvectors from strong block minimal bases linearizations and how such recovery method allows to obtain the ones for Fiedler-like linearizations.

It can be shown that if $\la_0$ is an eigenvalue of a regular rational matrix $G(\la)$ and $L(\la)$ is a strong block minimal bases linearization of $G(\la)$, then the bases of $\mathcal{N}_r (G(\lambda_0))$ and $\mathcal{N}_r (L(\lambda_0))$ are related as the right minimal bases in Theorem \ref{th_bases_sbml}(a) with the only change of replacing $\la$ by $\la_0$. Similarly, the bases of $\mathcal{N}_\ell (G(\lambda_0))$ and $\mathcal{N}_\ell (L(\lambda_0))$ are related as the left minimal bases in Theorem \ref{th_bases_sbml}(b) with $\la$ replaced by $\la_0$. The proofs of these facts follow a pattern similar to the proofs of Theorem \ref{th_bases_sbml} but are much simpler, since all the arguments concerning the degrees are not needed.

Once the relationships between the bases of $\mathcal{N}_r (G(\lambda_0))$ and $\mathcal{N}_\ell (G(\lambda_0))$ and the bases of the corresponding null-spaces of the strong block minimal bases linearizations of $G(\la)$ are known, the same argument as in the proof of Corollary \ref{cor.consequence2} proves that if $L_G (\la)$ is a FP, or a proper GFP or a GFPR linearization of $G(\la)$, then bases of $\mathcal{N}_r (G(\lambda_0))$ and $\mathcal{N}_\ell (G(\lambda_0))$ can be obtained from those of $\mathcal{N}_r (L_G(\lambda_0))$ and $\mathcal{N}_\ell (L_G(\lambda_0))$ as in Corollary \ref{cor.consequence2} with the only change of replacing again $\la$ by $\la_0$. This provides alternative proofs to those in \cite{AlBe16-2,DasAl19laa,DaAl20} for the recovery of eigenvectors of a regular rational matrix from those of its Fiedler-like linearizations.

%Remarks similar to the ones above hold for $\mathbb{M}_1$ and %$\mathbb{M}_2$-strong linearizations of $G(\la)$ and Theorems \ref{th.m1recovery} %and \ref{thm.M2recovery}. We mention in this context that the recovery of %eigenvectors of a regular rational matrix from those of its $\mathbb{M}_1$ and %$\mathbb{M}_2$-strong linearizations of $G(\la)$ was studied in \cite[Section %5]{DoMaQu19}.

}

\section{Conclusions} \label{Sec.Conclusions}  { In this paper a complete theory about the relationship between the minimal bases and indices of a rational matrix and those of its polynomial system matrices, as well as those of its strong linearizations, has been developed.}

The original contributions of this paper are organized into two clearly different parts. On the one hand those in Sections \ref{Sec_PSM}, \ref{Sec_linearizations} and \ref{Sec_minimalindicesstrong}, which deal with general polynomial system matrices, general linearizations and general strong linearizations of rational matrices. On the other hand those in Sections \ref{Sec_blockminimalbases}, \ref{Sec_M12minimalbases} { and \ref{sec.Fiedler}, which deal with specific (though large) families of strong linearizations and establish connections among them. More precisely that Fiedler-like linearizations are particular cases of strong block minimal bases linearizations modulo permutations}. In the case of polynomial system matrices, we have shown that, under the standard assumption of minimality and a certain additional condition of properness, the minimal indices of the polynomial system matrices and their transfer functions are exactly the same and their minimal bases are easily related to each other. These results are connected to pioneer results by { Verghese, Van Dooren and Kailath \cite{VeDoka79,Ver80}, who proved similar results under different and nonequivalent assumptions}. In contrast, we have shown that the minimal bases and indices of a rational matrix and those of its linearizations and strong linearizations are not related to each other in general, and that only the sums of the left and the right minimal indices are determined by each other in the case of strong linearizations. { This latter result is based on the fundamental index sum theorem obtained by Paul Van Dooren in \cite{VeDoka79}}.

In the case of the families of strong block minimal bases linearizations and $\M_1$ and $\M_2$-strong linearizations of rational matrices, we have proved that the minimal indices and bases of the linearizations and the rational matrices are easily related to each other and that any of them can be obtained from the others and vice versa. { The results for strong block minimal bases linearizations are obtained by using those for polynomial system matrices in Section 3, and they imply easily the results for $\M_1$ and $\M_2$-strong linearizations. In the case of Fiedler-like pencils, we have shown how they allow to recover the minimal indices and bases of a rational matrix as a consequence of the results for strong block minimal bases linearizations. This approach gives alternative proofs to the results in \cite{DasAl19laa,DaAl20}. In this context, it is worth to emphasize the important unifying role played by strong block minimal bases linearizations of rational matrices. In addition, to compute minimal bases and indices by applying algorithms for pencils to strong block minimal bases linearizations allow to deal with rectangular matrices. %This is an advantage with respect to the other families of linearizations recently developed in the literature, since rectangular rational matrices appear often in the applications involving minimal bases.
}

%------------------------------------------------------------------------------------------------------------------


\begin{thebibliography}{ }
%------------------------------------------------------------------------------------------------------------------


\bibitem{AlBe16} R. Alam, N. Behera,
Linearizations for rational matrix functions and Rosenbrock
system polynomials, \textit{SIAM J. Matrix Anal. Appl.}, 37 (1), 354--380, 2016.


\bibitem{AlBe16-2}
R. Alam, N. Behera, Recovery of eigenvectors of rational matrix functions from Fiedler-like linearizations, \textit{Linear Algebra Appl.}, 510, 373--394, 2016.


\bibitem{AlBe18} R. Alam, N. Behera, Generalized Fiedler pencils for rational matrix functions, \textit{SIAM J. Matrix Anal. Appl.}, 39 (2), 587--610, 2018.


\bibitem{AmMaZa16} A. Amparan, S. Marcaida, I. Zaballa,
On coprime rational function matrices,  \textit{Linear Algebra Appl.}, 507, 1--31, 2016.

 \bibitem{AmDoMaZa16} A. Amparan, F. M. Dopico, S. Marcaida, I. Zaballa,
Strong linearizations of rational matrices, \textit{Manchester Institute for Mathematical Sciences EPrints, The University of Manchester}, MIMS EPrint: 2016.51.

\bibitem{AmDoMaZa18} A. Amparan, F. M. Dopico, S. Marcaida, I. Zaballa,
Strong linearizations of rational matrices, \textit{SIAM J. Matrix Anal. Appl.}, 39 (4), 1670--1700, 2018.


\bibitem{AnDoHoMac19} L. M. Anguas, F. M. Dopico, R. Hollister, D. S. Mackey, Van Dooren's index sum theorem and rational matrices with prescribed structural data, \textit{SIAM J. Matrix Anal. Appl.}, 40 (2), 720--738, 2019.


\bibitem{AntVol04} E. N. Antoniou, S. Vologiannidis, A new family of companion forms of polynomial matrices, \textit{Electron. J. Linear Algebra}, 11, 78--87, 2004.


\bibitem{beelen1987}
Th. G. J. Beelen, G. W. Veltkamp, Numerical computation of a coprime factorization of a transfer function matrix, \textit{Syst.  Contr. Lett.}, 9 (4), 281--288, 1987.



\bibitem{BuDoFuRych15} M. I. Bueno, F. M. Dopico, S. Furtado, M. Rychnovsky, Large vector spaces of block-symmetric strong linearizations of matrix polynomials, \textit{Linear Algebra Appl.}, 477, 165--210, 2015.

 \bibitem{BuDoPe18} M. I. Bueno, F. M. Dopico, J. P\'erez, R. Saavedra, B. Zykoski, A simplified approach to Fiedler-like pencils via block minimal bases pencils, \textit{Linear Algebra Appl.}, 547, 45--104, 2018.


 \bibitem{DasAl19laa} R. K. Das, R. Alam, Recovery of minimal bases and minimal indices of rational matrices from Fiedler-like pencils, \textit{Linear Algebra Appl.}, 566, 34--60, 2019.


\bibitem{DaAl19_2} R. K. Das, R. Alam,
Affine spaces of strong linearizations for rational matrices and the recovery of eigenvectors and minimal indices, \textit{Linear Algebra Appl.}, 569, 335--368, 2019.


\bibitem{DaAl20} R. K. Das, R. Alam, Structured strong linearizations of structured rational matrices, {\textit arXiv:2008.00427v1}


\bibitem{DeDoMa14} F. De Ter\'an, F. M. Dopico, D. S. Mackey,
Spectral equivalence of matrix polynomials and the index sum
theorem, \textit{Linear Algebra Appl.}, 459, 264--333, 2014.

\bibitem{DeDoVa15} F. De Ter\'an, F. M. Dopico, P. Van Dooren,
Matrix polynomials with completely prescribed eigenstructure, \textit{SIAM J. Matrix Anal. Appl.}, 36, 302--328, 2015.

%\bibitem{DeDoMaVa16} F. De Ter\'an, F. M. Dopico, S. Mackey, P. Van Dooren,
%Polynomial zigzag matrices, dual minimal bases, and the realization of completely singular polynomials, \textit{Linear Algebra Appl.}, 488, 460--504, 2016.

\bibitem{Dopico:2016:BlockKronecker}
F. M. Dopico, P. W. Lawrence, J. P{\'{e}}rez, P. Van Dooren,
Block Kronecker linearizations of matrix polynomials and their backward
errors, MIMS EPrint 2016.34, Manchester Institute for Mathematical Sciences, The University of Manchester, UK, 2016.

\bibitem{Dopico:2018:BlockKronecker}
F. M. Dopico, P. W. Lawrence, J. P{\'{e}}rez, P. Van Dooren,
  Block Kronecker linearizations of matrix polynomials and their backward
  errors, \textit{Numer. Math.}, 140, 373--426, 2018.

\bibitem{DoMaQu19}
F. M. Dopico, S. Marcaida, M. C. Quintana,
Strong linearizations of rational matrices with polynomial part expressed in an orthogonal basis, \textit{Linear Algebra Appl.}, 570, 1--45, 2019.


\bibitem{DoMaQuVD19} F. M. Dopico, S. Marcaida, M. C. Quintana, P. Van Dooren, Local linearizations of rational matrices with application to rational approximations of nonlinear eigenvalue problems,  \textit{Linear Algebra Appl., 604, 441--475, 2020}.


\bibitem{FasSal17} H. Fa{\ss}bender, P. Saltenberger, On vector spaces of linearizations for matrix polynomials in orthogonal bases, \textit{Linear Algebra Appl.}, 525, 59--83, 2017.


\bibitem{FasSal18}
H. Fa{\ss}bender, P. Saltenberger, Block Kronecker ansatz spaces for matrix polynomials, \textit{Linear Algebra Appl.}, 542, 118--148, 2018.





\bibitem{Fiedler03} M. Fiedler, A note on companion matrices, \textit{Linear Algebra Appl.}, 372, 325--331, 2003.


\bibitem{guttel-tisseur}  S. G\"uttel, F. Tisseur, The nonlinear eigenvalue problem, \textit{Acta Numer.}, 26, 1--94, 2017.


\bibitem{nleigs}
 S. G\"{u}ttel, R. Van Beeumen, K. Meerbergen, W. Michiels, NLEIGS: A class of fully rational Krylov methods for nonlinear eigenvalue problems, \textit{SIAM J. Sci. Comput.}, 36 (6), A2842--A2864, 2014.


\bibitem{Fo75}
G. D. Forney, Minimal bases of rational vector spaces with applications to multivariable linear systems, \textit{SIAM J. Control},  13 (3), 143--520, 1975.


\bibitem{Kailath80} T. Kailath, {\em Linear Systems}, Prentice Hall, New
Jersey, 1980.


\bibitem{kung1980} S. Kung, T. Kailath, Fast projection methods for minimal design problems in linear system theory, \textit{Automatica}, 16 (4), 399--403, 1980.


\bibitem{automatic} P. Lietaert, J. P\'erez, B. Vandereycken, K. Meerbergen, Automatic  rational  approximation  and  linearization  of  nonlinear  eigenvalue  problems, submitted.  Available  as arXiv:1801.08622v2


\bibitem{MMMM06} D. S. Mackey, N. Mackey, C. Mehl, V. Mehrmann, Vector spaces of linearizations for matrix polynomials, \textit{SIAM J. Matrix Anal. Appl.}, 28 (4), 971--1004, 2006.


\bibitem{mehrmann-voss-04}
V. Mehrmann, H. Voss, Nonlinear eigenvalue problems:  A challenge for modern eigenvalue methods, \textit{GAMM-Mitt.}, 27, 121--152, 2004.


\bibitem{Rosen70} H. H. Rosenbrock,
{\em State-space and Multivariable Theory}, Thomas Nelson and Sons,
London, 1970.

 \bibitem{silverman1983system}
L. M. Silverman,  A. Kitap{\c{c}}i, System structure at infinity, \textit{Syst.  Contr. Lett.}, 3 (3), 123--131, 1983.



\bibitem{SuBai11} Y. Su, Z. Bai, Solving rational eigenvalue problems via
linearization, \textit{SIAM J. Matrix Anal. Appl.}, 32 (1), 201--216, 2011.


\bibitem{vd1979} P. Van Dooren, The computation of Kronecker's canonical form of a singular pencil, \textit{Linear Algebra Appl.}, 27, 103--140, 1979.



\bibitem{vd1981} P. Van Dooren, The generalized eigenstructure problem in linear system theory, \textit{IEEE Trans. Automat. Contr.}, 26 (1), 111--129, 1981.

\bibitem{Vard91}
A. I. G. Vardulakis, {\em Linear Multivariable Control}, John Wiley and Sons, New York, 1991.


\bibitem{Ver80} G. Verghese, Comments on `Properties of the system matrix of a generalized state-space system', \textit{Int. J. Control}, 31 (5), 1007--1009, 1980.

\bibitem{VeDoka79} G. Verghese, P. Van Dooren, T. Kailath,
Properties of the system  matrix of a generalized state-space system, \textit{Int. J. Control}, 30 (2), 235--243, 1979.


\bibitem{VolAnt11} S. Vologiannidis, E.N. Antoniou, A permuted factors approach for the linearization of polynomial matrices, \textit{Math. Control Signals Systems}, 22, 317--342, 2011.



\bibitem{wang1973} S. Wang, E. Davison, A minimization algorithm for the design of linear multivariable systems, \textit{IEEE Trans. Automat. Contr.}, 18 (3), 220--225, 1973.


%%%%%%%%%%%%%%%%%%%%%%%%%%%%%%%%%%%%%%%%%%%%%%%%%%%%%%%%%%%
%%%%%%%%%%%%%%%%%%%%%%%%%%%%%%%%%%%%%%%%%%%%%%%%%%%%%%%%%%%%%%%%%%%%%%%%%%%%%%%%%%%%%%%%%%%%%%%%%%%%%%%%%%%%%%%%%%%%%%%%%%%%%%%%%%%%

%\bibitem{AmMaZa04} Amparan A., Marcaida S., Zaballa I.,
%Assignment of infinite structure to an open-loop system. \textit{Linear
%Algebra Appl.}, 379 (2004) 249--266.


%\bibitem{AmMaZa15} A. Amparan, S. Marcaida, I. Zaballa,
%Finite and infinite structures of rational matrices: a local
%approach, \textit{Electronic Journal of Linear Algebra}, 30 (2015) 196--226.

%\bibitem{DeDoMa14} F. De Teran, F. Dopico S. Mackey,
%Spectral equivalence of matrix polynomials and the index sum
%theorem, \textit{Linear Algebra and Applications}, 459 (2014) 264--333.


%\bibitem{Lanc08} P. Lancaster, Linearization of regular matrix polynomials,
%\textit{Electronic Journal of Linear Algebra}, 17 (2008) 21--27.

%\bibitem{MMMM06} S. Mackey, N. Mackey, C. Mehl, V. Mehrmann,
%Vectors spaces of linearizations of matrix polynomials,
%\textit{SIAM J. Matrix Anal. Appl.}, 28 (4) (2006) 971--1004.

%\bibitem{Marques79} E. Marques de S\'{a}, Imbedding Conditions for
%$\lambda$-matrices, \textit{Linear Algebra and its Applications}, 24 (1979)
%33--50.


%\bibitem{Schutter00} B. de Schutter, Minimal state-space realization in linear system theory:
%an overview,
%\textit{Journal of Computation and Applied Mathematics}, 121 (2000) 331--354.



%\bibitem{Thom79} R. C. Thompson, Interlacing Inequalities for
%Invariant Factors, \textit{Linear Algebra and its Applications}, 24 (1979) 1--31.


\end{thebibliography}
\end{document}